\documentclass[reqno, 11pt, a4paper]{amsart}

\usepackage[a4paper, left=2.5cm, top=2.5cm, right=2.5cm, bottom=2.5cm]{geometry}
\numberwithin{equation}{section}

\usepackage{csquotes}
\usepackage[backend=biber, maxbibnames=5, style=numeric-comp,
            isbn=false, url=false, doi=true, giveninits=true]{biblatex}
\usepackage{amsmath, amsfonts, amssymb, bm, mathtools}
\usepackage{mathrsfs}
\usepackage[inline]{enumitem}
\usepackage{booktabs}
\usepackage{multirow}
\usepackage{etoolbox}
\AtBeginEnvironment{table}{%
  \setlength{\aboverulesep}{1pt}%
  \setlength{\belowrulesep}{1.5pt}%
  \setlength{\doublerulesep}{1pt}%
}
\usepackage{array}
\usepackage[table,dvipsnames]{xcolor}
\usepackage{graphicx}
\graphicspath{{figures/}}
\usepackage{float}

\usepackage{tikz}
\usetikzlibrary{shapes.misc,matrix,fit,positioning,arrows.meta,decorations.markings,decorations.shapes,decorations.pathreplacing,calc,shapes.geometric,arrows,shapes,patterns,patterns.meta,math}
\tikzset{Matrix/.style={matrix of nodes, font=\footnotesize,text height=1pt, text depth=0.5pt, text width=8.5pt, align=center, column sep=0pt, row sep=0pt, nodes in empty cells}}
\usepackage{pgf}
\usepgflibrary{fpu}

\newtheorem{theorem}{Theorem}
\numberwithin{theorem}{section}

\newtheorem{remark}[theorem]{Remark}

\usepackage{acronym}
\acrodef{fe}[FE]{finite element}
\acrodef{dof}[DOF]{degree of freedom}
\acrodefplural{dof}[DOFs]{degrees of freedom}
\acrodef{gmg}[GMG]{geometric multigrid}

\usepackage{hyperref}
\hypersetup{breaklinks=true, colorlinks=true,
            linkcolor=black, citecolor=blue, urlcolor=blue}
\newcommand{\tab}[1]{Table~\ref{#1}}
\newcommand{\fig}[1]{Figure~\ref{#1}}
\newcommand{\sect}[1]{Section~\ref{#1}}

\newcommand{\ol}[1]{\overline{#1}}
\newcommand{\norm}[2]{\left\|#2\right\|_{#1}}   
\newcommand{\inner}[3]{\left(#1,\,#2\right)_{#3}}   

\newcommand{\Div}{\operatorname{div}}
\newcommand{\Grad}{\operatorname{grad}}
\newcommand{\Curl}{\operatorname{\bm{curl}}}

\newcommand{\Qk}[1]{\mathbb{Q}_{#1}}     
\newcommand{\RT}[1]{\mathbb{RT}_{#1}}    

\newcommand{\dom}{\Omega}                       
\newcommand{\bnd}{\Gamma}                        
\newcommand{\bnddir}{\Gamma_{\rm D}}                 
\newcommand{\bndneu}{\Gamma_{\rm N}}                 
\newcommand{\bgdom}{\Omega^{\rm bg}}             
\newcommand{\actdom}{\Omega_h}                   
\newcommand{\indom}{\Omega_h^{\rm in}}
\newcommand{\outdom}{\Omega_h^{\rm out}}
\newcommand{\aggdom}{\Omega_h^{\rm ag}}       
\newcommand{\gdom}{G}                           
\newcommand{\normal}{n}                          
\newcommand{\jump}[1]{[\![ #1 ]\!]}                   

\newcommand{\Tbg}{\mathcal{T}_h^{\rm bg}}        
\newcommand{\Th}{\mathcal{T}_h}                  
\newcommand{\Tin}{\mathcal{T}_h^{\rm in}}        
\newcommand{\Tout}{\mathcal{T}_h^{\rm out}}       
\newcommand{\Tcut}{\mathcal{T}_h^{\rm cut}}      
\newcommand{\Tagg}{\mathcal{T}_h^{\rm ag}}       
\newcommand{\Taggall}{\mathcal{T}_h^{\rm ag, all}}       
\newcommand{\agg}{\omega}                        

\newcommand{\meshsize}{h}                        

\newcommand{\Hgrad}[1]{H(\Grad,#1)}
\newcommand{\Hdiv}[1]{\bm{H}(\Div,#1)}
\newcommand{\Hcurl}[1]{\bm{H}(\Curl,#1)}

\newcommand{\Lgrad}{\Lambda^{\rm g}}
\newcommand{\Ldiv}{\Lambda^{\rm d}}
\newcommand{\Lcurl}{\Lambda^{\rm c}}
\newcommand{\Lx}{\Lambda^{\rm x}}
\newcommand{\Lxh}{\Lx_h}                           
\newcommand{\Lxact}{\Lx_{\rm act}}                 

\newcommand{\Lgradh}{\Lgrad_h}
\newcommand{\Ldivh}{\Ldiv_h}
\newcommand{\Lcurlh}{\Lcurl_h}

\newcommand{\Vgrad}{V^{\rm g}_h}
\newcommand{\Vdiv}{\bm{V}^{\rm d}_h}
\newcommand{\Vcurl}{\bm{V}^{\rm c}_h}

\newcommand{\Vzeroloc}{\mathcal{V}^{0}}
\newcommand{\Vgradloc}{\mathcal{V}^{\rm g}}
\newcommand{\Vdivloc}{\bm{\mathcal{V}}^{\rm d}}
\newcommand{\Vcurlloc}{\bm{\mathcal{V}}^{\rm c}}

\newcommand{\bh}{b_h}                             

\newcommand{\eps}{\varepsilon}                    

\newcommand{\Ah}{A_h^{\rm x}}           
\newcommand{\Aact}{A_{\rm act}^{\rm x}}    
\newcommand{\Dh}{D_h}           
\newcommand{\Bx}[1]{B_{\rm #1}^{\rm x}}   

\newcommand{\Tl}{\mathcal{T}_{\ell}}
\newcommand{\Tlbg}{\mathcal{T}_{\ell}^{\rm bg}}
\newcommand{\doml}{\Omega_{\ell}}
\newcommand{\domlmo}{\Omega_{\ell-1}}

\newcommand{\Vl}{V_{\ell}^{\rm x}}
\newcommand{\Vlmo}{V_{\ell-1}^{\rm x}}
\newcommand{\Vld}{(\Vl)'}                           
\newcommand{\Vlmod}{(\Vlmo)'}
\newcommand{\Vhd}{(\Vh)'}                           
\newcommand{\Llact}{\Lambda_{\ell, \rm act}^{\rm x}}

\newcommand{\Alact}{A_{\ell, \rm act}^{\rm x}}

\newcommand{\Ilhat}{\widehat{I}_{\ell}}
\newcommand{\Il}{I_{\ell}}

\newcommand{\Dl}{D_{\ell}}
\newcommand{\Bl}{B_{\ell}^{\rm x}}

\newcommand{\Star}[1]{{\rm star}(#1)}                

\newcommand{\actdomH}{\Omega_H}               

\newcommand{\Vh}{V_h^{\rm x}}                

\newcommand{\Id}{{\rm I}}

\newcommand{\Tlmo}{\mathcal{T}_{\ell-1}}          
\newcommand{\Tlagg}{\mathcal{T}_{\ell}^{\rm ag}}   
\newcommand{\Tlaggall}{\mathcal{T}_{\ell}^{\rm ag, all}}   

\newcommand{\Ll}{\Lambda^{\rm x}_{\ell}}           
\newcommand{\Al}{A_{\ell}^{\rm x}}                         
\newcommand{\Rl}{R_{\ell}}                         
\newcommand{\Dlact}{D_{\ell, \rm act}}             
\newcommand{\Rh}{R_h}                              
\newcommand{\VCycle}{\operatorname{VCycle}}        
\newcommand{\Bgmg}{\Bx{hom}}                  
\newcommand{\Bmg}{\Bx{aux}}                   
\newcommand{\Bagmg}{B_{\text{a-gmg}}^{\rm x}}  

\newcommand{\gridap}{\href{https://github.com/gridap/Gridap.jl}{\texttt{Gridap}}}
\newcommand{\gridapembedded}{\href{https://github.com/gridap/GridapEmbedded.jl}{\texttt{GridapEmbedded}}}
\newcommand{\gridapsolvers}{\href{https://github.com/gridap/GridapSolvers.jl}{\texttt{GridapSolvers}}}

\usepackage[normalem]{ulem}

\title[Multigrid preconditioners for unfitted de Rham discretisations]{Geometric multigrid preconditioners for unfitted discretisations of the de~Rham complex}
\date{\today}

\address{School of Mathematics, Monash University, Clayton, Victoria 3800, Australia}
\author[S. Badia]{Santiago Badia}
\email{santiago.badia@monash.edu}
\author[J. Manyer]{Jordi Manyer}
\email{jordi.manyer@monash.edu}

\begin{document}

\begin{abstract}
We propose two families of multigrid preconditioners for unfitted finite element discretisations of the Riesz maps of the de~Rham complex, in $H(\Grad)$, $H(\Curl)$ and $H(\Div)$, stabilised by cell aggregation. The first is a geometric multigrid method applied directly to the unfitted operators, on a hierarchy of background meshes that share the same boundary representation. Its smoothers are star patches built on the aggregated mesh, so that each patch captures the kernel of the stabilisation. The second is a multilevel auxiliary-space preconditioner that combines a patch smoother localised around the cut cells with the body-fitted problem on the active mesh, which is in turn solved by a non-nested geometric multigrid method. Both are used as building blocks of block-diagonal preconditioners for the mixed Hodge Laplacians. Numerical experiments in two and three dimensions, on geometries with a controlled cut fraction and on general geometries, show iteration counts that are robust with respect to the cut configuration, the mesh size and the number of levels. We also report results for higher-order discretisations, weakly imposed Dirichlet boundary conditions and a face-based stabilisation.
\end{abstract}

\maketitle

\section{Introduction}

Body-fitted finite element methods rely on a mesh that conforms to the boundary of the physical domain. For complex geometries, and even more so for evolving domains, generating such a mesh is a costly and fragile step of the simulation pipeline, in particular in three dimensions. Unfitted \ac{fe} methods remove this bottleneck by embedding the physical domain within a simple background domain, which can then be meshed by a structured grid. The physical domain boundary is then represented implicitly (e.g.\ by a level-set function) or explicitly (e.g.\ by a surface mesh in STL format) and intersected with the background mesh. Quadrature rules, using for example high-order quadrature techniques for implicitly defined domains~\cite{Saye2015}, then need to be designed on these intersections.

While efficient, this approach introduces the small cut-cell problem, in which arbitrarily small intersections between the mesh and the boundary lead to severe cut-dependent numerical instabilities and ill-conditioning \cite{dePrenter-ArchComputMethodsEng-2023j}. To restore robustness, several stabilisation strategies have been developed.
Ghost-penalty methods (also known as CutFEM)~\cite{Burman2010ghost} are among the most popular in this family. Robustness with respect to the cut configuration is recovered by penalising the jumps of the derivatives across the facets of the cut cells, or, in bulk form, the deviation of the discrete function from a global polynomial on patches of cells. Alternatively, cell-aggregation techniques, such as the aggregated unfitted \ac{fe} method (AgFEM) \cite{badia2018agfem}, achieve the same by aggregating badly cut cells with interior neighbours and constraining their \acp{dof}. The two approaches are closely related, and the aggregates can also be used to define locking-free stabilisation terms, either through aggregate-wise projections or through jumps across the facets interior to each aggregate~\cite{BadiaNeivaVerdugo2022}.
Both approaches were originally designed for grad-conforming approximations. Considerable effort has recently gone into extending them to the whole de~Rham complex \cite{badia2026darcy,DiPietroDroniouNilsson2025,Lehrenfeld-MathComput-2024j}, so that mixed and structure-preserving formulations of more complex problems can also be discretised on unfitted meshes.

As for any \ac{fe} method, the practical applicability of unfitted discretisations depends on the efficient solution of the resulting linear systems. Efficient and robust solvers are therefore crucial, especially for high-order discretisations, where the cost of a direct solve becomes prohibitive.
For fitted conforming discretisations, \ac{gmg} is a popular choice. Combined with star-patch (subspace correction) smoothers \cite{Arnold2000,farrell2021}, it delivers mesh-independent convergence across the whole de~Rham complex, and degree-independent convergence as well \cite{Pavarino1994,Schoberl2008,BrubeckFarrell2024}. Its main practical limitation is the mesh hierarchy. A sequence of nested, body-fitted meshes is impossible to generate for a general geometry, and the alternatives (non-nested hierarchies or adaptive refinement of a coarse fitted mesh) require additional machinery. Multigrid is a much more natural fit for unfitted methods, because the background mesh is structured, so coarsening it is trivial, and the geometry is carried by the level set rather than by the mesh.

Despite this, \ac{gmg} solvers for unfitted discretisations are a relatively unexplored area.
A \ac{gmg} solver for immersed \ac{fe} methods was designed in \cite{deprenter2020mg}, building on the cell-based additive Schwarz preconditioner of \cite{deprenter2019as}. The resulting algorithm was shown numerically to be robust with respect to both the cut configuration and the mesh size. A similar \ac{gmg} algorithm was proposed in \cite{ludescher2020mg} for unfitted interface problems, with a tailored prolongation and an interface-correction smoother that solves a global system on the interface unknowns.
For discontinuous Galerkin discretisations, Lu, Rupp and Kanschat~\cite{lu2022} introduced a homogeneous \ac{gmg}, in which the unfitted discretisation is used on every level of a hierarchy that inherits the same boundary representation, together with vertex-star smoothers. Cui and Kanschat~\cite{cui2025cutfemgpu} applied the same ideas to unfitted discretisations with ghost penalty and Nitsche boundary conditions, using multiplicative vertex-star patches on a GPU.
To the best of our knowledge, the only convergence theory for vertex-star multigrid on ghost-penalty discretisations is the two-level analysis of Wichrowski~\cite{wichrowski2026cutfemmg}, which is uniform in the mesh size and the cut and extends to the W-cycle under a smallness assumption on the two-level contraction rate. The analysis of the V-cycle is, however, left open.

Another interesting approach is the one in \cite{gross2023cutfemprec}, where the authors analysed a two-block additive Schwarz preconditioner for unfitted $H(\Grad)$ discretisations with ghost penalty. The solver relies on a space decomposition that splits the global problem between an interior Dirichlet problem on the uncut cells and a strip block made of the cut basis functions with diagonal scaling. The splitting is proved stable uniformly in the mesh size and the cut. The interior Dirichlet problem is then solved using an algebraic multigrid solver. All of these works are restricted to the $H(\Grad)$-conforming setting.

 The extension of \ac{gmg} algorithms to the complete de~Rham complex for unfitted discretisations has not been addressed in the literature.
In this work we design and test multilevel preconditioners for stabilised unfitted discretisations of the $H(\Grad)$-, $H(\Curl)$- and $H(\Div)$-conforming Riesz maps. Our contributions are the following.
\begin{itemize}
  \item We extend the homogeneous unfitted multigrid of \cite{lu2022,cui2025cutfemgpu} from $H(\Grad)$ to $H(\Curl)$- and $H(\Div)$-conforming discretisations. We refer to it as the \emph{homogeneous} multigrid preconditioner. To this end, we introduce \emph{aggregate-based} star patches, built on the aggregates of the stabilisation, that are designed to satisfy the required kernel decomposition property and thereby ensure robustness of the smoother.
  \item We also design an alternative multigrid solver that relies almost exclusively on body-fitted problems set on the active domain. The design connects two levels of preconditioning. First, we build an \emph{auxiliary-space} preconditioner~\cite{Xu1996} for the original unfitted problem, in the spirit of \cite{gross2023cutfemprec}, which pairs a patch smoother with a global body-fitted problem on the active mesh. Then, we solve the auxiliary problem with a non-nested \ac{gmg} preconditioner on a hierarchy of active meshes.
  \item We use the resulting solvers as the diagonal blocks of operator preconditioners for the mixed Hodge Laplacians of the complex.
\end{itemize}
Numerical experiments in two and three dimensions show that the proposed preconditioners are robust with respect to the mesh size, the number of \ac{gmg} levels, the cut configuration and the weight $\kappa$ of the differential term. This holds both for a stabilisation based on aggregate-wise projections and for a face-based ghost-penalty stabilisation restricted to the aggregates \cite{BadiaNeivaVerdugo2022}. We also explore robustness with respect to both higher-order discretisations and the presence of essential (Dirichlet) boundary conditions.

The analysis of the proposed multigrid methods is not addressed here, and is left for future work. Due to the non-nested nature of the hierarchies, such analysis requires a substantially different approach from the body-fitted case and thus warrants its own paper.

\section{Geometrical discretisation}\label{sec:geometry}

Let $\dom \subset \mathbb{R}^d$, with $d \in \{2,3\}$ the space dimension, be the physical domain with boundary $\bnd = \partial\dom$. In unfitted discretisation techniques, we assume that the physical domain is embedded in a background domain $\bgdom$ such that $\dom \subset \bgdom$. The background domain $\bgdom$ is typically chosen to be a geometrically simple domain, such as a rectangle or a box, that can be easily meshed.

Let $\Tbg$ be a conforming, quasi-uniform and shape-regular partition of $\bgdom$. We denote by $h_T$ the diameter of a cell $T \in \Tbg$, and by $\meshsize \doteq \max_{T \in \Tbg} h_T$ the characteristic mesh size of $\Tbg$.

We now introduce some geometrical notation that will be used throughout the paper. We partition the background mesh $\Tbg$ into $\{\Tin, \Tout, \Tcut\}$, where
\[
  \Tin \coloneqq \{\, T \in \Tbg : T \subset \dom \,\}, \qquad
  \Tout \coloneqq \{\, T \in \Tbg : T \cap \dom = \emptyset \,\}, \qquad
  \Tcut \coloneqq \Tbg \setminus (\Tin \cup \Tout).
\]
The cells in $\Tin$, $\Tout$ and $\Tcut$ are referred to as \emph{interior}, \emph{exterior} and \emph{cut} cells, respectively. The cut cells are those that intersect the boundary $\bnd$, while the interior and exterior cells are those that do not. We refer to the union of all interior and cut cells as the \emph{active mesh} $\Th \coloneqq \Tin \cup \Tcut$, as they are the only cells that play a role in the discretisation.

We define the active domain $\actdom$ as the interior of the closure of $\bigcup_{T\in\Th} T$ and similarly denote by $\indom$ and $\outdom$ the unions of cells in $\Tin$ and $\Tout$, respectively.

\subsection{Aggregation of small cut cells}\label{subsec:aggregates}

One of the main challenges of unfitted discretisations is the so-called \emph{small-cut cell problem}, which arises when the intersection of cut cells with the physical domain has arbitrarily small measure. This leads to ill-conditioned global systems and a loss of convergence for the underlying discretisation.
Over the years, several strategies have been proposed to address this issue. Within the context of stabilised unfitted methods, the two main approaches are (i) face-based stabilisation and (ii) bulk stabilisation.

In this work, we consider unfitted discretisations that address the small-cut cell problem via bulk-stabilisation terms defined by means of aggregate-wise $L^2$-projections. To this end, we introduce an aggregated mesh $\Taggall$, obtained by associating each cell in $\Tcut$ to a single interior \emph{root} cell in $\Tin$ through a chain of cells, consecutive cells sharing a facet that intersects $\dom$ (an interior or cut facet). Each cell in $\Th$ belongs to exactly one aggregate, while interior cells that absorb no cut cell form singleton aggregates. We denote by $\Tagg$ the set of non-trivial aggregates, that is the aggregates that contain at least one cut cell, and by $\aggdom$ the union of all non-trivial aggregates in $\Tagg$. We assume that the number of cells per aggregate is uniformly bounded, which holds for $\meshsize$ small enough with respect to the features of $\bnd$ \cite{badia2018agfem}. Since every aggregate contains an interior cell, this implies that the diameter $h_\agg$ of $\agg$ satisfies $h_\agg \le C \meshsize$ and $|\agg \cap \dom| \ge c\, |\agg|$, with constants $C$, $c$ independent of $\meshsize$ and of the position of $\bnd$ relative to the mesh. \fig{fig:geometry} illustrates these geometric definitions.

\begin{figure}
  \centering
  \begin{minipage}[b]{0.32\textwidth}
    \centering
    \input{figures/meshes/mesh_inoutcut.tex}\\[2pt]
    {\small (a) Cell partition}
  \end{minipage}\hfill
  \begin{minipage}[b]{0.32\textwidth}
    \centering
    \input{figures/meshes/mesh_active.tex}\\[2pt]
    {\small (b) Active mesh $\Th$}
  \end{minipage}\hfill
  \begin{minipage}[b]{0.32\textwidth}
    \centering
    \input{figures/meshes/mesh_aggregates.tex}\\[2pt]
    {\small (c) Aggregates $\Tagg$}
  \end{minipage}
  \caption{The embedded geometric setting on an $8\times8$ background mesh $\Tbg$ cut by a
  domain $\dom$ (boundary $\bnd$ drawn in black). \emph{(a)} The partition of $\Tbg$ into
  interior cells $\Tin$ (yellow), cut cells $\Tcut$ (grey) and exterior cells $\Tout$ (white).
  \emph{(b)} The active mesh $\Th = \Tin \cup \Tcut$. \emph{(c)} The non-trivial aggregates
  $\Tagg$, each drawn in a different colour; every cut cell is absorbed into an aggregate
  rooted at an interior cell, while interior cells that absorb no cut cell (uncoloured)
  remain non-agglomerated.}
  \label{fig:geometry}
\end{figure}

 Given a polynomial space $\mathcal{V}$, we define by $\mathcal{V}^{\rm ag}(\Tagg)$ the space of discontinuous polynomials over $\Tagg$ whose restriction to each aggregate $\agg$ belongs to $\mathcal{V}(\agg)$\footnote{In the implementation, the polynomial basis of $\mathcal{V}(\agg)$ is defined on a bounding box of $\agg$ of diameter $h_\agg \sim h$, which has no effect on the space but keeps the local problems well conditioned.}, that is
\[
  \mathcal{V}^{\rm ag}(\Tagg) \coloneqq \{ v \in L^2(\aggdom) : v|_\agg \in \mathcal{V}(\agg),\ \forall\, \agg \in \Tagg \}.
\]

We can then define the aggregate-wise $L^2$-projection $\Pi^{\rm ag} : L^2(\aggdom) \to \mathcal{V}^{\rm ag}(\Tagg)$, defined by
\begin{equation} \label{eq:agg-proj}
  \inner{\Pi^{\rm ag} u - u}{v}{\Tagg} = 0
  \qquad \forall\, v \in \mathcal{V}(\agg).
\end{equation}
This projection can be computed by solving a local problem on each aggregate $\agg \in \Tagg$. Using these projections, we define a stabilisation bilinear form $S^{\rm ag}$ as
\begin{equation}\label{eq:stabilisation}
  S^{\rm ag}(u,v) \coloneqq \inner{u - \Pi^{\rm ag} u}{v - \Pi^{\rm ag} v}{\Tagg} \quad \forall u,v \in L^2(\aggdom).
\end{equation}

\section{Discrete formulations}\label{sec:formulations}

This section introduces the discrete formulations for the unfitted \ac{fe} method, including the definition of global bilinear forms, stabilisation terms, and the imposition of Dirichlet boundary conditions through Nitsche's method.

We consider the Riesz bilinear forms on a generic domain $\gdom$,
\begin{equation}
  \begin{aligned}
    \Lgrad(u,v; \gdom) &\coloneqq \inner{u}{v}{\gdom} + \kappa^2 \inner{\Grad u}{\Grad v}{\gdom}, \\
    \Lcurl(\bm{u},\bm{v}; \gdom) &\coloneqq \inner{\bm{u}}{\bm{v}}{\gdom} + \kappa^2 \inner{\Curl \bm{u}}{\Curl \bm{v}}{\gdom}, \\
    \Ldiv(\bm{u},\bm{v}; \gdom) &\coloneqq \inner{\bm{u}}{\bm{v}}{\gdom} + \kappa^2 \inner{\Div \bm{u}}{\Div \bm{v}}{\gdom}, \\
    \Lambda^0(u,v; \gdom) &\coloneqq \inner{u}{v}{\gdom}.
  \end{aligned}
\end{equation}

In two dimensions $\Curl \bm{u} = \partial_1 u_2 - \partial_2 u_1$ is scalar-valued, and in three dimensions it is vector-valued. Let $\Vgrad \subset \Hgrad{\actdom}$, $\Vcurl \subset \Hcurl{\actdom}$, $\Vdiv \subset \Hdiv{\actdom}$ and $V_h^0 \subset L^2(\actdom)$ be global \ac{fe} spaces defined in terms of local polynomial spaces $\Vgradloc$, $\Vcurlloc$, $\Vdivloc$ and $\Vzeroloc$ respectively as follows:

\begin{equation}\label{eq:spaces}
  \begin{aligned}
  \Vgrad &\coloneqq \{ v \in \Hgrad{\actdom} : v|_T \in \Vgradloc(T),\ \forall T \in \Th \}, \\
  \Vcurl &\coloneqq \{ \bm{v} \in \Hcurl{\actdom} : \bm{v}|_T \in \Vcurlloc(T),\ \forall T \in \Th \}, \\
  \Vdiv &\coloneqq \{ \bm{v} \in \Hdiv{\actdom} : \bm{v}|_T \in \Vdivloc(T),\ \forall T \in \Th \}, \\
  V_h^0 &\coloneqq \{ v \in L^2(\actdom) : v|_T \in \Vzeroloc(T),\ \forall T \in \Th \}.
  \end{aligned}
\end{equation}

On the Cartesian cells used in this work, the local spaces are the tensor-product Lagrangian space $\Vgradloc = \Qk{k}$, the N\'ed\'elec space of the first kind $\Vcurlloc = \mathbb{ND}_k$ and the Raviart--Thomas space $\Vdivloc = \RT{k}$, where $k \ge 1$ for $\Vgradloc$ and $k \ge 0$ otherwise, and $\Vzeroloc = \Qk{k}$ is the discontinuous space. The \emph{lowest order} thus refers to $\Qk{1}$, $\mathbb{ND}_0$, $\RT{0}$ and discontinuous $\Qk{0}$.

We assume that these spaces form a discrete de~Rham complex, \[ \Vgrad \xrightarrow{\ \Grad\ } \Vcurl \xrightarrow{\ \Curl\ } \Vdiv \xrightarrow{\ \Div\ } V_h^0 \] in three dimensions, and $\Vgrad \xrightarrow{\Grad} \Vcurl \xrightarrow{\Curl} V_h^0$ together with $\Vgrad \xrightarrow{\Grad^\perp} \Vdiv \xrightarrow{\Div} V_h^0$, $\Grad^\perp \coloneqq (\partial_2, -\partial_1)$, in two dimensions, which the families $\Qk{k}$, $\mathbb{ND}_{k-1}$, $\RT{k-1}$ and discontinuous $\Qk{k-1}$ satisfy. 
In what follows, the superscript ${\rm x} \in \{{\rm g},{\rm c},{\rm d},0\}$ refers to any of the four cases, so that $\Lx$ and $\Vh$ denote the corresponding bilinear form and \ac{fe} space, and $\mathrm{D}^{\rm x} \in \{\Grad, \Curl, \Div, 0\}$ the corresponding differential operator, with $\mathrm{D}^0 \coloneqq 0$.
For ${\rm x} \in \{{\rm g},{\rm c},{\rm d}\}$ we denote by ${\rm x}'$ the conformity of the space that $\mathrm{D}^{\rm x}$ maps into in the complex above, namely
\[
  {\rm g}' = {\rm c}, \qquad
  {\rm c}' = \begin{cases} {\rm d} & \text{in three dimensions}, \\ 0 & \text{in two dimensions}, \end{cases} \qquad
  {\rm d}' = 0,
\]
so that $\mathrm{D}^{\rm x} V_h^{\rm x} \subset V_h^{{\rm x}'}$. We use the construction of \sect{subsec:aggregates} for each of the local spaces $\mathcal{V}^{\rm x}$, ${\rm x} \in \{{\rm g},{\rm c},{\rm d},0\}$, and write $S^{\rm ag}_{\rm x}$ for the corresponding stabilisation \eqref{eq:stabilisation}. On each space we then define two global bilinear forms, that we refer to as the \emph{active} and \emph{unfitted} bilinear forms
\begin{equation} \label{eq:global-biforms}
  \Lxact(u,v) \coloneqq \Lx(u,v; \actdom), \qquad
  \Lxh(u,v) \coloneqq \Lx(u,v; \dom) + \tau^2 S^{\rm x}(u,v),
\end{equation}
where $\tau > 0$ is a stabilisation parameter.

We denote by $\Ah, \Aact : \Vh \to \Vhd$ the operators induced by $\Lxh$ and $\Lxact$, respectively, where $\Vhd$ is the dual of $\Vh$.
The stabilisation terms $S^{\rm x}$ are assumed to fulfil the following property:
\begin{equation}\label{eq:stab-property}
  \Lxh(u,u) \simeq \Lxact(u,u) \qquad \forall u \in \Vh,
\end{equation}
where $a \simeq b$ means $c\, b \le a \le C\, b$ with constants $c, C > 0$ that are independent of the mesh size $\meshsize$ and of the position of $\bnd$ relative to the background mesh (the \emph{cut configuration}), and of $\kappa$. This is the cut-robustness property of the stabilisation, see \cite[Assumption~3.2]{badia2026darcy}.

A possible choice for these terms that satisfies \eqref{eq:stab-property} is defined in terms of the aggregate-wise $L^2$-projections in \eqref{eq:stabilisation} as
\begin{equation} \label{eq:stabilisation-biforms}
  \begin{aligned}
    S^{\rm g}(u,v) &\coloneqq S_{\rm g}^{\rm ag}(u,v) + \kappa^2 S_{{\rm g}'}^{\rm ag}(\Grad u, \Grad v), \\
    S^{\rm c}(\bm{u},\bm{v}) &\coloneqq S_{\rm c}^{\rm ag}(\bm{u},\bm{v}) + \kappa^2 S_{{\rm c}'}^{\rm ag}(\Curl \bm{u}, \Curl \bm{v}), \\
    S^{\rm d}(\bm{u},\bm{v}) &\coloneqq S_{\rm d}^{\rm ag}(\bm{u},\bm{v}) + \kappa^2 S_{0}^{\rm ag}(\Div \bm{u}, \Div \bm{v}), \\
    S^0(u,v) &\coloneqq S_{0}^{\rm ag}(u,v),
  \end{aligned}
\end{equation}
where $S^0$ has no second term since $\mathrm{D}^0 = 0$.

We also consider Dirichlet boundary conditions imposed weakly through a Nitsche-type bilinear form $\bh$ defined on the boundary $\bnd$ of the physical domain $\dom$, for ${\rm x} \in \{{\rm g},{\rm c},{\rm d}\}$. The space $V_h^0$ does not admit boundary conditions, and we set $b^0 \coloneqq 0$. Let $\bnd = \bnddir \cup \bndneu$ with $\bnddir \neq \emptyset$ and define
\begin{equation}
\begin{aligned}
  b^{\rm g}(u,v) &\coloneqq \kappa^2 \Big( \frac{\gamma}{\meshsize}\, \inner{u}{v}{\bnddir} - \inner{\Grad u \cdot \normal }{v}{\bnddir} - \inner{u}{\Grad v \cdot \normal}{\bnddir} \Big), \\
  b^{\rm c}(\bm{u},\bm{v}) &\coloneqq \kappa^2 \Big( \frac{\gamma}{\meshsize}\, \inner{\bm{u}\times\normal}{\bm{v}\times\normal}{\bnddir} + \inner{\Curl \bm{u}}{\bm{v}\times\normal}{\bnddir} + \inner{\bm{u}\times\normal}{\Curl \bm{v}}{\bnddir} \Big), \\
  b^{\rm d}(\bm{u},\bm{v}) &\coloneqq \kappa^2 \Big( \frac{\gamma}{\meshsize}\, \inner{\bm{u}\cdot\normal}{\bm{v}\cdot\normal}{\bnddir} - \inner{\Div \bm{u}}{\bm{v}\cdot\normal}{\bnddir} - \inner{\bm{u}\cdot\normal}{\Div \bm{v}}{\bnddir} \Big),
\end{aligned}
\end{equation}
where $\bm{u}\times\normal$ is the tangential trace\footnote{In two dimensions, $\bm{u}\times\normal \coloneqq u_1 n_2 - u_2 n_1$ is a scalar; with this convention, the signs of $b^{\rm c}$ are those of the integration by parts of $\Curl\Curl$ in both two and three dimensions.}. Then the associated Nitsche-type bilinear forms are defined as
\begin{equation} \label{eq:nitsche-weakforms}
  \widehat{\Lxact}(u,v) \coloneqq \Lxact(u,v) + b^{\rm x}(u,v), \qquad
  \widehat{\Lxh}(u,v) \coloneqq \Lxh(u,v) + b^{\rm x}(u,v).
\end{equation}

With these boundary conditions, the purely diffusive Poisson problem is recovered from $\widehat{\Lgradh}$ in the limit $\kappa \to \infty$, and is well posed due to the Dirichlet boundary conditions on $\bnddir \neq \emptyset$.

\begin{remark}[Essential boundary conditions]\label{rem:essential-bcs}
  Imposing Dirichlet boundary conditions weakly through Nitsche's method introduces additional terms whose analysis is delicate. The convergence and stability of $\widehat{\Lgradh}$ are well established, see e.g.\ \cite{badia2018agfem} for the aggregated setting and \cite{BurmanHansbo2012} for ghost-penalty stabilisation, and those of $\widehat{\Ldivh}$ have recently been proven in \cite{badia2026darcy}. To our knowledge, the convergence and stability of $\widehat{\Lcurlh}$ remain open problems. The difficulty is not specific to the unfitted setting: for body-fitted meshes, Nitsche's method for the $H(\Div)$ Riesz map and the Darcy problem is analysed in \cite{BurmanPuppi2022}, and for Maxwell's problem in $H(\Curl)$ in \cite{BoffiCodinaTurk2026}.
\end{remark}

\section{Abstract multigrid framework for unfitted discretisations}\label{sec:gmg}

This section introduces the multigrid machinery used by both our multilevel preconditioners.

\subsection{Mesh hierarchy}\label{subsec:hierarchy}

Let $\{\Tlbg\}_{\ell=0,\dots,L}$ be a hierarchy of nested background meshes of $\bgdom$, with $\mathcal{T}^{\rm bg}_L = \Tbg$ being the finest mesh and $\mathcal{T}^{\rm bg}_0$ the coarsest. Since the background meshes are not fitted to the physical domain, such a hierarchy is always easy to construct. One can for instance choose a hierarchy of nested Cartesian meshes, built by uniform refinement, where each cell of $\mathcal{T}^{\rm bg}_{\ell-1}$ is split into $2^d$ equal children to obtain $\Tlbg$, so that
\begin{equation}
  h_{\ell-1} = 2 h_\ell, \qquad \ell = 1, \ldots, L,
\end{equation}
where $h_\ell$ denotes the mesh size of $\Tlbg$. This is the choice that will be used for the numerical experiments in this work.

The physical boundary $\bnd$ is represented identically on all levels, that is $\bnd_0 = \dots = \bnd_L = \bnd$: the physical domain $\dom$ and its boundary are the \emph{same} geometric objects at every level $\ell = 0, \ldots, L$, and no level-dependent approximation $\dom_\ell$ or $\bnd_\ell$ is introduced. In practice, this is achieved by performing the cuts on the finest mesh and passing this boundary representation on to the coarser meshes. We refer to this construction as the \emph{inherited boundary}.

\begin{remark}
  The inherited boundary contrasts with methods that approximate $\bnd$ at each level (e.g.\ by a piecewise-linear level-set approximation $\bnd_\ell$). It is the standard choice in \ac{gmg} for unfitted methods~\cite{lu2022}, and it is what makes the level operators below well defined on a common physical domain.
The inherited boundary does not compromise the cost of the coarse levels. The integrals over the cut cells of a coarse level can be computed with moment-fitted quadratures \cite{Badia-JComputPhys-2022k}, whose number of points per cut cell is set by the polynomial degree and not by the fine geometry they resolve, so that integration and assembly on each level scale with its number of cells.
\end{remark}

\begin{figure}
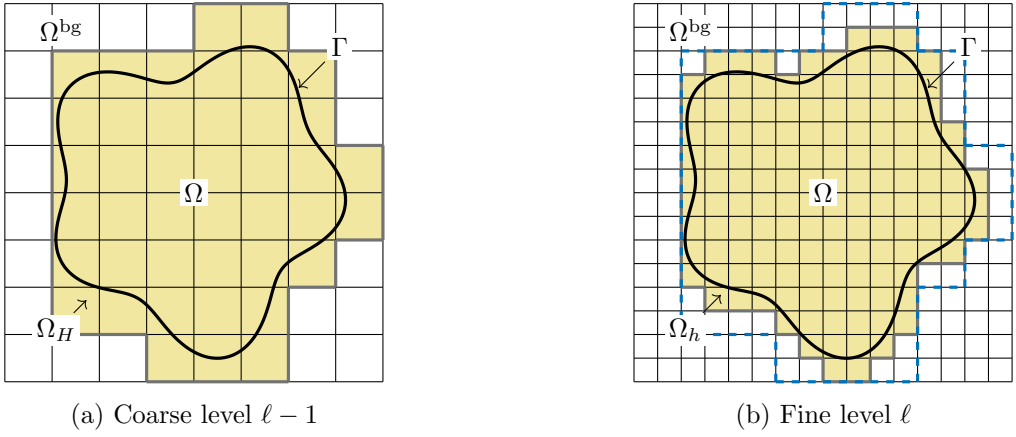

  \centering
  \begin{minipage}[b]{0.48\textwidth}
    \centering
    \input{figures/meshes/mesh_hierarchy_coarse.tex}\\[2pt]
    {\small (a) Coarse level $\ell-1$}
  \end{minipage}\hfill
  \begin{minipage}[b]{0.48\textwidth}
    \centering
    \input{figures/meshes/mesh_hierarchy_fine.tex}\\[2pt]
    {\small (b) Fine level $\ell$}
  \end{minipage}
  \caption{Two consecutive levels of the hierarchy. Both levels share the same background domain $\bgdom$ (the
  unit square) and the same physical boundary $\bnd$ (thick black), and are equipped with the
  nested background meshes $\mathcal{T}^{\rm bg}_{\ell-1}$ and $\Tlbg$ (thin black). The active meshes $\Tlmo$ and $\Tl$,
  i.e.\ the cells meeting the physical domain $\dom$, are shaded, and their active domains
  $\domlmo$ and $\doml$ are outlined in grey. The dashed blue line in \emph{(b)} is
  $\partial\domlmo$, which illustrates the inclusion $\doml \subset \domlmo$ of the active
  domains even though the active meshes themselves are non-nested.}
  \label{fig:two-level-setup}
\end{figure}

For each level $\ell$, we repeat the construction of \sect{sec:geometry} to obtain a hierarchy of non-nested active meshes $\Tl$ and their associated active domains $\doml$. Although the active meshes are non-nested, the inheritance of the boundary representation guarantees that the active domains form a decreasing sequence, $\doml \subset \domlmo$ for all $\ell > 0$. Indeed, a cell of $\mathcal{T}^{\rm bg}_{\ell-1}$ has a non-empty intersection with $\dom$ if and only if at least one of its children in $\Tlbg$ has. Hence every active cell of $\Tl$ is contained in an active cell of $\Tlmo$, whereas an active cell of $\Tlmo$ may contain exterior cells of $\Tl$. This is illustrated in \fig{fig:two-level-setup}.

On each level we also repeat the aggregation of \sect{subsec:aggregates}, which yields a hierarchy of aggregated meshes $\Tlaggall$ and $\Tlagg$, as well as the associated aggregate-wise projections and stabilisation forms \eqref{eq:agg-proj}--\eqref{eq:stabilisation}. The aggregates are built independently on each level. No relation between the aggregates of consecutive levels is required by the algorithms below.

\subsection{Level spaces and transfer operators}\label{subsec:spaces-transfer}

For each level $\ell$ and for ${\rm x} \in \{{\rm g, c, d}, 0\}$, we repeat the construction in \sect{sec:formulations} on the active mesh $\Tl$ and define the \ac{fe} space $\Vl$ of \eqref{eq:spaces}, so that $V_L^{\rm x} = \Vh$.
Since the active meshes are non-nested, the spaces $\Vl$ are also non-nested. However, we do have that for $\ell > 0$
\[
  \Vlmo|_{\doml} \subset \Vl.
\]
Then, let $\Ilhat : \Vlmo|_{\doml} \to \Vl$ be the natural injection operator. We define the intergrid prolongation operator $\Il : \Vlmo \to \Vl$ as
\begin{equation}\label{eq:prolongation}
  \Il(v_{\ell-1}) \coloneqq \Ilhat(v_{\ell-1}|_{\doml}) \quad \forall v_{\ell-1} \in \Vlmo,
\end{equation}
and the restriction operator $\Il^T : \Vld \to \Vlmod$ as its adjoint.

\subsection{The abstract V-cycle algorithm}\label{subsec:vcycle}

Let $\{\Al: \Vl \to \Vld\}_{\ell=0}^{L}$ and $\{\Dl : \Vld \to \Vl\}_{\ell=1}^{L}$ be two arbitrary families of symmetric positive definite operators, where $\Vld$ denotes the dual of $\Vl$. For the smoothing operators $\Dl$, positive semidefiniteness suffices for the definitions below, which is needed in \sect{sec:auxspace}. Let $\Il$ be the intergrid prolongation operator defined in \eqref{eq:prolongation}. On each level, for $n_\ell \geq 1$ and a damping parameter $\theta_\ell > 0$, we define the smoother $\Rl : \Vld \to \Vl$ as $n_\ell$ steps of the damped Richardson iteration with $\Dl$,
\begin{equation}\label{eq:effective-smoother}
  \Id - \Rl \Al \coloneqq (\Id - \theta_\ell \Dl \Al)^{n_\ell}.
\end{equation}
With these ingredients, the V-cycle preconditioner $\Bl$ for $\Al$ on level $\ell$ is defined recursively as follows. On the coarsest level $\ell = 0$, $B_0^{\rm x} \coloneqq (A_0^{\rm x})^{-1}$ is the exact inverse of the coarse grid operator. For every other level $\ell = 1, \dots, L$, $\Bl : \Vld \to \Vl$ is the operator that given $b_\ell \in \Vld$ performs the following steps:
\begin{enumerate}[label=\roman*]
  \item Sets $x^{(0)} = 0 \in \Vl$
  \item Pre-smoothing: $x^{(1)} = x^{(0)} + \Rl(b_\ell - \Al x^{(0)})$
  \item Coarse correction: $x^{(2)} = x^{(1)} + \Il B_{\ell-1}^{\rm x} \Il^T (b_\ell - \Al x^{(1)})$
  \item Post-smoothing: $x^{(3)} = x^{(2)} + \Rl^T(b_\ell - \Al x^{(2)})$
  \item Returns $\Bl b_\ell = x^{(3)}$
\end{enumerate}

This multigrid cycle can be expressed in operator form as
\begin{equation}\label{eq:error-propagation}
  \Id - \Bl \Al = (\Id - \Rl^T \Al)(\Id - \Il B_{\ell-1}^{\rm x} \Il^T \Al) (\Id - \Rl \Al) \, \quad \ell = 1, \dots, L.
\end{equation}
Using the adjoint smoother $\Rl^T$ in the post-smoothing step ensures that $\Bl$ is a symmetric operator. Since $\Dl$ is symmetric, $\Rl = \big(\Id - (\Id - \theta_\ell \Dl \Al)^{n_\ell}\big)(\Al)^{-1}$ is itself symmetric, so that $\Rl^T = \Rl$ for all the smoothers of this work. We keep the transpose to cover non-symmetric smoothers. For $\theta_\ell$ sufficiently small, $\Bl$ is symmetric positive definite and can therefore be used directly as a preconditioner within a preconditioned conjugate gradient iteration. In the experiments of \sect{sec:results} it is used within the flexible GMRES method (FGMRES).

We write
\begin{equation}\label{eq:vcycle-notation}
  B_L^{\rm x} = \VCycle\big(\{\Al\}_{\ell}, \{\Dl\}_{\ell}\big)
\end{equation}
for the finest-level preconditioner. The two multigrid preconditioners of this paper, in \sect{sec:embedded-gmg} and \sect{sec:auxspace}, are two instances of \eqref{eq:vcycle-notation} that differ in the families of operators and smoothers.

\subsection{Aggregate-based star-patch smoothers}\label{subsec:smoothers}

In this section, we generalise the definition of star-patch smoothers to meshes obtained by agglomeration. Star-patch smoothers are a class of additive Schwarz smoothers based on patches of cells around mesh entities \cite{farrell2021}. They are popular in multilevel methods for their robustness and local support, and have been shown to be effective across the whole de~Rham complex \cite{Arnold2000}.

The key condition for an effective and robust subspace-correction-based relaxation is that the subspaces associated with the patches of the smoother satisfy the so-called \emph{kernel decomposition property}, namely that the kernel of the differential term of the Riesz map\footnote{The operator $\Al$ itself is positive definite. The kernel is that of its semidefinite part, which dominates $\Al$ for large $\kappa$.}  together with its stabilisation (the second term in the expressions in \eqref{eq:stabilisation-biforms}), restricted to $\Vl$, can be decomposed as the sum of its intersections with the local patch spaces. In other words, we require that any kernel function can be written as the sum of kernel functions drawn from the subspaces. This is a non-trivial condition to fulfil when the operator has a large kernel, such as divergence-free and curl-free vector fields for $\Div$- and $\Curl$-conforming discretisations, respectively.
On fitted meshes, the classical vertex-star smoother satisfies this property for the whole de~Rham complex \cite{Arnold2000}. This is no longer the case in the unfitted setting, where the kernel of the stabilisation must also be taken into account.

The idea behind the stabilisations used for unfitted methods is to (weakly) force the piecewise-polynomial functions on each aggregate to behave as a single polynomial on the entire aggregate. Such functions cannot be supported on part of an aggregate only. Thus, we need to enlarge the support of the star patches to contain every aggregate they intersect. This is the idea behind aggregate-based star-patch smoothers that we propose in this work.

Let $\mathcal{T}$ be a fitted mesh of a bounded polytopal domain in $\mathbb{R}^d$. Given a colouring of its cells into connected components, we define $\mathcal{P}$ as the polytopal mesh obtained by agglomerating the cells of each connected component into a single polytope. For a trivial colouring, that is, when each cell forms its own connected component, $\mathcal{P}$ coincides with $\mathcal{T}$. By construction, $\mathcal{T}$ and $\mathcal{P}$ cover the same domain, and we denote by $\partial \mathcal{P}$ its boundary.

Let $\mathcal{F}(\mathcal{P})$ be the set of faces of $\mathcal{P}$, defined from the cells of $\mathcal{P}$ as
\[
  \mathcal{F}(\mathcal{P}) \coloneqq \{\, \ol{P} \cap \ol{P'} : P, P' \in \mathcal{P},\ P \ne P' \,\} \cup \{\, \ol{P} \cap \partial \mathcal{P} : P \in \mathcal{P} \,\},
\]
where the intersections are understood in the sense of the $(d-1)$-dimensional measure, i.e.\ only intersections of positive $(d-1)$-dimensional measure are faces.
The first set contains the faces shared by two polytopes, while the second contains the \emph{boundary faces}, that is, the intersection of each polytope with the domain boundary. These faces need not be planar due to the polytopal nature of the mesh. In particular, they are piecewise linear for $d = 2$ and piecewise planar for $d = 3$. Applying the same construction to the faces, now retaining the $(d-2)$-dimensional intersections of pairs of faces, gives the edges of the mesh for $d = 3$ (the vertices for $d = 2$), and applying it once more to the edges gives the vertices for $d = 3$. We denote by $\mathcal{N}(\mathcal{P})$ the set of vertices of the polytopal mesh $\mathcal{P}$.

With the above definitions, we can define the notion of a star patch around a mesh entity. Let $\mathcal{X}(\mathcal{P})$ be a set of mesh entities of $\mathcal{P}$, typically $\mathcal{F}(\mathcal{P})$, $\mathcal{N}(\mathcal{P})$, or a subset thereof.
For each $\nu \in \mathcal{X}(\mathcal{P})$, we define the \emph{star patch} around $\nu$ as
\begin{equation}
  \Star{\nu} \coloneqq \bigcup \big\{\, P \in \mathcal{P} : \nu \subset \ol{P} \,\big\},
\end{equation}
that is, the union of the polytopes whose closure contains $\nu$. Since every polytope is a union of cells of $\mathcal{T}$, $\Star{\nu}$ is also a union of cells of $\mathcal{T}$.

Given a \ac{fe} space $V$ on $\mathcal{T}$, let $\Gamma_\nu \coloneqq \bigcup \{\, f \in \mathcal{F}(\mathcal{P}) : f \subset \partial \mathcal{P},\ \nu \subset f \,\}$ be the union of the boundary faces incident to $\nu$. We define the local space associated with the star patch $\Star{\nu}$ as
\begin{equation}
  V_\nu \coloneqq \{\, v \in V : \operatorname{supp} v \subset \ol{\Star{\nu}},\ \operatorname{tr} v = 0 \text{ on } \partial \Star{\nu} \setminus \Gamma_\nu \,\},
\end{equation}
where $\operatorname{tr} v$ denotes the trace relevant to the conformity of $V$, that is, $v$, $\bm{v} \times \normal$ or $\bm{v} \cdot \normal$ for $H(\Grad)$, $H(\Curl)$ and $H(\Div)$, respectively, and no condition for the discontinuous space $V_h^0$, whose basis functions are cell-interior. Equivalently, $V_\nu$ is spanned by the basis functions of $V$ whose \acp{dof} are attached to mesh entities lying in the interior of $\Star{\nu}$ or on $\Gamma_\nu$.
In other words, these local spaces have homogeneous Dirichlet boundary conditions (on the relevant trace) on the whole boundary of the star patch, except on the boundary faces incident to $\nu$, where the \acp{dof} remain free, as in the global space $V$.

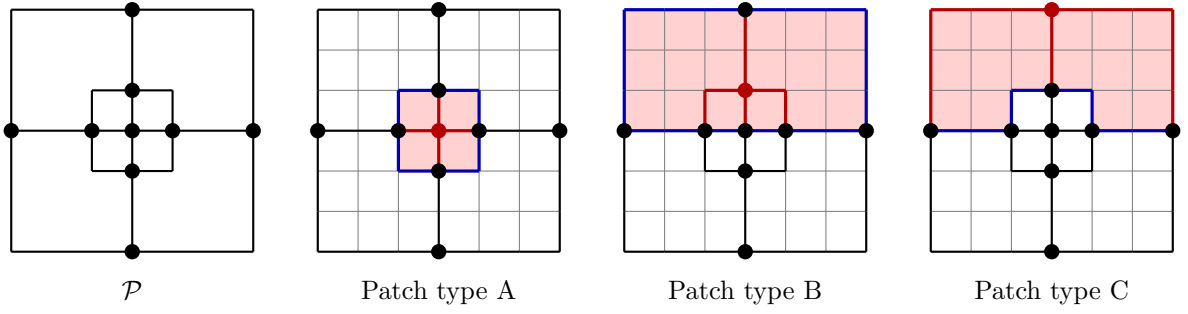
\begin{figure}
  \centering
  \begin{minipage}[b]{0.24\textwidth}
    \centering
    \begin{tikzpicture}[scale=3.2]
\draw[thick, black] (0.0000,0.0000) -- (1.0000,0.0000);
\draw[thick, black] (0.0000,0.5000) -- (1.0000,0.5000);
\draw[thick, black] (0.0000,1.0000) -- (1.0000,1.0000);
\draw[thick, black] (0.3333,0.3333) -- (0.6667,0.3333);
\draw[thick, black] (0.3333,0.6667) -- (0.6667,0.6667);
\draw[thick, black] (0.0000,0.0000) -- (0.0000,1.0000);
\draw[thick, black] (0.5000,0.0000) -- (0.5000,1.0000);
\draw[thick, black] (1.0000,0.0000) -- (1.0000,1.0000);
\draw[thick, black] (0.3333,0.3333) -- (0.3333,0.6667);
\draw[thick, black] (0.6667,0.3333) -- (0.6667,0.6667);
\fill[black] (0.0000,0.5000) circle (0.9pt);
\fill[black] (0.3333,0.5000) circle (0.9pt);
\fill[black] (0.5000,0.0000) circle (0.9pt);
\fill[black] (0.5000,0.3333) circle (0.9pt);
\fill[black] (0.5000,0.5000) circle (0.9pt);
\fill[black] (0.5000,0.6667) circle (0.9pt);
\fill[black] (0.5000,1.0000) circle (0.9pt);
\fill[black] (0.6667,0.5000) circle (0.9pt);
\fill[black] (1.0000,0.5000) circle (0.9pt);
\end{tikzpicture}\\[2pt]
    {\small $\mathcal{P}$}
  \end{minipage}\hfill
  \begin{minipage}[b]{0.24\textwidth}
    \centering
    \begin{tikzpicture}[scale=3.2]
\fill[red!18] (0.3333,0.5000) rectangle (0.5000,0.6667);
\fill[red!18] (0.5000,0.3333) rectangle (0.6667,0.5000);
\fill[red!18] (0.5000,0.5000) rectangle (0.6667,0.6667);
\fill[red!18] (0.3333,0.3333) rectangle (0.5000,0.5000);
\draw[thin, black!50] (0,0.0000) -- (1,0.0000);
\draw[thin, black!50] (0.0000,0) -- (0.0000,1);
\draw[thin, black!50] (0,0.1667) -- (1,0.1667);
\draw[thin, black!50] (0.1667,0) -- (0.1667,1);
\draw[thin, black!50] (0,0.3333) -- (1,0.3333);
\draw[thin, black!50] (0.3333,0) -- (0.3333,1);
\draw[thin, black!50] (0,0.5000) -- (1,0.5000);
\draw[thin, black!50] (0.5000,0) -- (0.5000,1);
\draw[thin, black!50] (0,0.6667) -- (1,0.6667);
\draw[thin, black!50] (0.6667,0) -- (0.6667,1);
\draw[thin, black!50] (0,0.8333) -- (1,0.8333);
\draw[thin, black!50] (0.8333,0) -- (0.8333,1);
\draw[thin, black!50] (0,1.0000) -- (1,1.0000);
\draw[thin, black!50] (1.0000,0) -- (1.0000,1);
\draw[thick, black] (0.0000,0.0000) -- (0.1667,0.0000);
\draw[thick, black] (0.0000,0.5000) -- (0.1667,0.5000);
\draw[thick, black] (0.0000,1.0000) -- (0.1667,1.0000);
\draw[thick, black] (0.1667,0.0000) -- (0.3333,0.0000);
\draw[thick, black] (0.1667,0.5000) -- (0.3333,0.5000);
\draw[thick, black] (0.1667,1.0000) -- (0.3333,1.0000);
\draw[thick, black] (0.3333,0.0000) -- (0.5000,0.0000);
\draw[very thick, blue!70!black] (0.3333,0.3333) -- (0.5000,0.3333);
\draw[very thick, blue!70!black] (0.3333,0.6667) -- (0.5000,0.6667);
\draw[thick, black] (0.3333,1.0000) -- (0.5000,1.0000);
\draw[thick, black] (0.5000,0.0000) -- (0.6667,0.0000);
\draw[very thick, blue!70!black] (0.5000,0.3333) -- (0.6667,0.3333);
\draw[very thick, blue!70!black] (0.5000,0.6667) -- (0.6667,0.6667);
\draw[thick, black] (0.5000,1.0000) -- (0.6667,1.0000);
\draw[thick, black] (0.6667,0.0000) -- (0.8333,0.0000);
\draw[thick, black] (0.6667,0.5000) -- (0.8333,0.5000);
\draw[thick, black] (0.6667,1.0000) -- (0.8333,1.0000);
\draw[thick, black] (0.8333,0.0000) -- (1.0000,0.0000);
\draw[thick, black] (0.8333,0.5000) -- (1.0000,0.5000);
\draw[thick, black] (0.8333,1.0000) -- (1.0000,1.0000);
\draw[thick, black] (0.0000,0.0000) -- (0.0000,0.1667);
\draw[thick, black] (0.0000,0.1667) -- (0.0000,0.3333);
\draw[thick, black] (0.0000,0.3333) -- (0.0000,0.5000);
\draw[thick, black] (0.0000,0.5000) -- (0.0000,0.6667);
\draw[thick, black] (0.0000,0.6667) -- (0.0000,0.8333);
\draw[thick, black] (0.0000,0.8333) -- (0.0000,1.0000);
\draw[very thick, blue!70!black] (0.3333,0.3333) -- (0.3333,0.5000);
\draw[very thick, blue!70!black] (0.3333,0.5000) -- (0.3333,0.6667);
\draw[thick, black] (0.5000,0.0000) -- (0.5000,0.1667);
\draw[thick, black] (0.5000,0.1667) -- (0.5000,0.3333);
\draw[thick, black] (0.5000,0.6667) -- (0.5000,0.8333);
\draw[thick, black] (0.5000,0.8333) -- (0.5000,1.0000);
\draw[very thick, blue!70!black] (0.6667,0.3333) -- (0.6667,0.5000);
\draw[very thick, blue!70!black] (0.6667,0.5000) -- (0.6667,0.6667);
\draw[thick, black] (1.0000,0.0000) -- (1.0000,0.1667);
\draw[thick, black] (1.0000,0.1667) -- (1.0000,0.3333);
\draw[thick, black] (1.0000,0.3333) -- (1.0000,0.5000);
\draw[thick, black] (1.0000,0.5000) -- (1.0000,0.6667);
\draw[thick, black] (1.0000,0.6667) -- (1.0000,0.8333);
\draw[thick, black] (1.0000,0.8333) -- (1.0000,1.0000);
\draw[very thick, red!70!black] (0.3333,0.5000) -- (0.5000,0.5000);
\draw[very thick, red!70!black] (0.5000,0.5000) -- (0.6667,0.5000);
\draw[very thick, red!70!black] (0.5000,0.3333) -- (0.5000,0.5000);
\draw[very thick, red!70!black] (0.5000,0.5000) -- (0.5000,0.6667);
\fill[black] (0.0000,0.5000) circle (0.9pt);
\fill[black] (0.3333,0.5000) circle (0.9pt);
\fill[black] (0.5000,0.0000) circle (0.9pt);
\fill[black] (0.5000,0.3333) circle (0.9pt);
\fill[black] (0.5000,0.6667) circle (0.9pt);
\fill[black] (0.5000,1.0000) circle (0.9pt);
\fill[black] (0.6667,0.5000) circle (0.9pt);
\fill[black] (1.0000,0.5000) circle (0.9pt);
\fill[red!70!black] (0.5000,0.5000) circle (0.9pt);
\end{tikzpicture}\\[2pt]
    {\small Patch type A}
  \end{minipage}\hfill
  \begin{minipage}[b]{0.24\textwidth}
    \centering
    \begin{tikzpicture}[scale=3.2]
\fill[red!18] (0.6667,0.6667) rectangle (0.8333,0.8333);
\fill[red!18] (0.3333,0.6667) rectangle (0.5000,0.8333);
\fill[red!18] (0.8333,0.8333) rectangle (1.0000,1.0000);
\fill[red!18] (0.0000,0.6667) rectangle (0.1667,0.8333);
\fill[red!18] (0.5000,0.6667) rectangle (0.6667,0.8333);
\fill[red!18] (0.1667,0.8333) rectangle (0.3333,1.0000);
\fill[red!18] (0.6667,0.5000) rectangle (0.8333,0.6667);
\fill[red!18] (0.8333,0.6667) rectangle (1.0000,0.8333);
\fill[red!18] (0.0000,0.5000) rectangle (0.1667,0.6667);
\fill[red!18] (0.1667,0.6667) rectangle (0.3333,0.8333);
\fill[red!18] (0.3333,0.5000) rectangle (0.5000,0.6667);
\fill[red!18] (0.6667,0.8333) rectangle (0.8333,1.0000);
\fill[red!18] (0.5000,0.5000) rectangle (0.6667,0.6667);
\fill[red!18] (0.0000,0.8333) rectangle (0.1667,1.0000);
\fill[red!18] (0.8333,0.5000) rectangle (1.0000,0.6667);
\fill[red!18] (0.3333,0.8333) rectangle (0.5000,1.0000);
\fill[red!18] (0.1667,0.5000) rectangle (0.3333,0.6667);
\fill[red!18] (0.5000,0.8333) rectangle (0.6667,1.0000);
\draw[thin, black!50] (0,0.0000) -- (1,0.0000);
\draw[thin, black!50] (0.0000,0) -- (0.0000,1);
\draw[thin, black!50] (0,0.1667) -- (1,0.1667);
\draw[thin, black!50] (0.1667,0) -- (0.1667,1);
\draw[thin, black!50] (0,0.3333) -- (1,0.3333);
\draw[thin, black!50] (0.3333,0) -- (0.3333,1);
\draw[thin, black!50] (0,0.5000) -- (1,0.5000);
\draw[thin, black!50] (0.5000,0) -- (0.5000,1);
\draw[thin, black!50] (0,0.6667) -- (1,0.6667);
\draw[thin, black!50] (0.6667,0) -- (0.6667,1);
\draw[thin, black!50] (0,0.8333) -- (1,0.8333);
\draw[thin, black!50] (0.8333,0) -- (0.8333,1);
\draw[thin, black!50] (0,1.0000) -- (1,1.0000);
\draw[thin, black!50] (1.0000,0) -- (1.0000,1);
\draw[thick, black] (0.0000,0.0000) -- (0.1667,0.0000);
\draw[very thick, blue!70!black] (0.0000,0.5000) -- (0.1667,0.5000);
\draw[very thick, blue!70!black] (0.0000,1.0000) -- (0.1667,1.0000);
\draw[thick, black] (0.1667,0.0000) -- (0.3333,0.0000);
\draw[very thick, blue!70!black] (0.1667,0.5000) -- (0.3333,0.5000);
\draw[very thick, blue!70!black] (0.1667,1.0000) -- (0.3333,1.0000);
\draw[thick, black] (0.3333,0.0000) -- (0.5000,0.0000);
\draw[thick, black] (0.3333,0.3333) -- (0.5000,0.3333);
\draw[very thick, blue!70!black] (0.3333,0.5000) -- (0.5000,0.5000);
\draw[very thick, blue!70!black] (0.3333,1.0000) -- (0.5000,1.0000);
\draw[thick, black] (0.5000,0.0000) -- (0.6667,0.0000);
\draw[thick, black] (0.5000,0.3333) -- (0.6667,0.3333);
\draw[very thick, blue!70!black] (0.5000,0.5000) -- (0.6667,0.5000);
\draw[very thick, blue!70!black] (0.5000,1.0000) -- (0.6667,1.0000);
\draw[thick, black] (0.6667,0.0000) -- (0.8333,0.0000);
\draw[very thick, blue!70!black] (0.6667,0.5000) -- (0.8333,0.5000);
\draw[very thick, blue!70!black] (0.6667,1.0000) -- (0.8333,1.0000);
\draw[thick, black] (0.8333,0.0000) -- (1.0000,0.0000);
\draw[very thick, blue!70!black] (0.8333,0.5000) -- (1.0000,0.5000);
\draw[very thick, blue!70!black] (0.8333,1.0000) -- (1.0000,1.0000);
\draw[thick, black] (0.0000,0.0000) -- (0.0000,0.1667);
\draw[thick, black] (0.0000,0.1667) -- (0.0000,0.3333);
\draw[thick, black] (0.0000,0.3333) -- (0.0000,0.5000);
\draw[very thick, blue!70!black] (0.0000,0.5000) -- (0.0000,0.6667);
\draw[very thick, blue!70!black] (0.0000,0.6667) -- (0.0000,0.8333);
\draw[very thick, blue!70!black] (0.0000,0.8333) -- (0.0000,1.0000);
\draw[thick, black] (0.3333,0.3333) -- (0.3333,0.5000);
\draw[thick, black] (0.5000,0.0000) -- (0.5000,0.1667);
\draw[thick, black] (0.5000,0.1667) -- (0.5000,0.3333);
\draw[thick, black] (0.5000,0.3333) -- (0.5000,0.5000);
\draw[thick, black] (0.6667,0.3333) -- (0.6667,0.5000);
\draw[thick, black] (1.0000,0.0000) -- (1.0000,0.1667);
\draw[thick, black] (1.0000,0.1667) -- (1.0000,0.3333);
\draw[thick, black] (1.0000,0.3333) -- (1.0000,0.5000);
\draw[very thick, blue!70!black] (1.0000,0.5000) -- (1.0000,0.6667);
\draw[very thick, blue!70!black] (1.0000,0.6667) -- (1.0000,0.8333);
\draw[very thick, blue!70!black] (1.0000,0.8333) -- (1.0000,1.0000);
\draw[very thick, red!70!black] (0.3333,0.6667) -- (0.5000,0.6667);
\draw[very thick, red!70!black] (0.5000,0.6667) -- (0.6667,0.6667);
\draw[very thick, red!70!black] (0.3333,0.5000) -- (0.3333,0.6667);
\draw[very thick, red!70!black] (0.5000,0.5000) -- (0.5000,0.6667);
\draw[very thick, red!70!black] (0.5000,0.6667) -- (0.5000,0.8333);
\draw[very thick, red!70!black] (0.5000,0.8333) -- (0.5000,1.0000);
\draw[very thick, red!70!black] (0.6667,0.5000) -- (0.6667,0.6667);
\fill[black] (0.0000,0.5000) circle (0.9pt);
\fill[black] (0.3333,0.5000) circle (0.9pt);
\fill[black] (0.5000,0.0000) circle (0.9pt);
\fill[black] (0.5000,0.3333) circle (0.9pt);
\fill[black] (0.5000,0.5000) circle (0.9pt);
\fill[black] (0.5000,1.0000) circle (0.9pt);
\fill[black] (0.6667,0.5000) circle (0.9pt);
\fill[black] (1.0000,0.5000) circle (0.9pt);
\fill[red!70!black] (0.5000,0.6667) circle (0.9pt);
\end{tikzpicture}\\[2pt]
    {\small Patch type B}
  \end{minipage}\hfill
  \begin{minipage}[b]{0.24\textwidth}
    \centering
    \begin{tikzpicture}[scale=3.2]
\fill[red!18] (0.6667,0.6667) rectangle (0.8333,0.8333);
\fill[red!18] (0.3333,0.6667) rectangle (0.5000,0.8333);
\fill[red!18] (0.8333,0.8333) rectangle (1.0000,1.0000);
\fill[red!18] (0.0000,0.6667) rectangle (0.1667,0.8333);
\fill[red!18] (0.5000,0.6667) rectangle (0.6667,0.8333);
\fill[red!18] (0.1667,0.8333) rectangle (0.3333,1.0000);
\fill[red!18] (0.6667,0.5000) rectangle (0.8333,0.6667);
\fill[red!18] (0.8333,0.6667) rectangle (1.0000,0.8333);
\fill[red!18] (0.0000,0.5000) rectangle (0.1667,0.6667);
\fill[red!18] (0.1667,0.6667) rectangle (0.3333,0.8333);
\fill[red!18] (0.6667,0.8333) rectangle (0.8333,1.0000);
\fill[red!18] (0.0000,0.8333) rectangle (0.1667,1.0000);
\fill[red!18] (0.8333,0.5000) rectangle (1.0000,0.6667);
\fill[red!18] (0.3333,0.8333) rectangle (0.5000,1.0000);
\fill[red!18] (0.1667,0.5000) rectangle (0.3333,0.6667);
\fill[red!18] (0.5000,0.8333) rectangle (0.6667,1.0000);
\draw[thin, black!50] (0,0.0000) -- (1,0.0000);
\draw[thin, black!50] (0.0000,0) -- (0.0000,1);
\draw[thin, black!50] (0,0.1667) -- (1,0.1667);
\draw[thin, black!50] (0.1667,0) -- (0.1667,1);
\draw[thin, black!50] (0,0.3333) -- (1,0.3333);
\draw[thin, black!50] (0.3333,0) -- (0.3333,1);
\draw[thin, black!50] (0,0.5000) -- (1,0.5000);
\draw[thin, black!50] (0.5000,0) -- (0.5000,1);
\draw[thin, black!50] (0,0.6667) -- (1,0.6667);
\draw[thin, black!50] (0.6667,0) -- (0.6667,1);
\draw[thin, black!50] (0,0.8333) -- (1,0.8333);
\draw[thin, black!50] (0.8333,0) -- (0.8333,1);
\draw[thin, black!50] (0,1.0000) -- (1,1.0000);
\draw[thin, black!50] (1.0000,0) -- (1.0000,1);
\draw[thick, black] (0.0000,0.0000) -- (0.1667,0.0000);
\draw[very thick, blue!70!black] (0.0000,0.5000) -- (0.1667,0.5000);
\draw[thick, black] (0.1667,0.0000) -- (0.3333,0.0000);
\draw[very thick, blue!70!black] (0.1667,0.5000) -- (0.3333,0.5000);
\draw[thick, black] (0.3333,0.0000) -- (0.5000,0.0000);
\draw[thick, black] (0.3333,0.3333) -- (0.5000,0.3333);
\draw[thick, black] (0.3333,0.5000) -- (0.5000,0.5000);
\draw[very thick, blue!70!black] (0.3333,0.6667) -- (0.5000,0.6667);
\draw[thick, black] (0.5000,0.0000) -- (0.6667,0.0000);
\draw[thick, black] (0.5000,0.3333) -- (0.6667,0.3333);
\draw[thick, black] (0.5000,0.5000) -- (0.6667,0.5000);
\draw[very thick, blue!70!black] (0.5000,0.6667) -- (0.6667,0.6667);
\draw[thick, black] (0.6667,0.0000) -- (0.8333,0.0000);
\draw[very thick, blue!70!black] (0.6667,0.5000) -- (0.8333,0.5000);
\draw[thick, black] (0.8333,0.0000) -- (1.0000,0.0000);
\draw[very thick, blue!70!black] (0.8333,0.5000) -- (1.0000,0.5000);
\draw[thick, black] (0.0000,0.0000) -- (0.0000,0.1667);
\draw[thick, black] (0.0000,0.1667) -- (0.0000,0.3333);
\draw[thick, black] (0.0000,0.3333) -- (0.0000,0.5000);
\draw[thick, black] (0.3333,0.3333) -- (0.3333,0.5000);
\draw[very thick, blue!70!black] (0.3333,0.5000) -- (0.3333,0.6667);
\draw[thick, black] (0.5000,0.0000) -- (0.5000,0.1667);
\draw[thick, black] (0.5000,0.1667) -- (0.5000,0.3333);
\draw[thick, black] (0.5000,0.3333) -- (0.5000,0.5000);
\draw[thick, black] (0.5000,0.5000) -- (0.5000,0.6667);
\draw[thick, black] (0.6667,0.3333) -- (0.6667,0.5000);
\draw[very thick, blue!70!black] (0.6667,0.5000) -- (0.6667,0.6667);
\draw[thick, black] (1.0000,0.0000) -- (1.0000,0.1667);
\draw[thick, black] (1.0000,0.1667) -- (1.0000,0.3333);
\draw[thick, black] (1.0000,0.3333) -- (1.0000,0.5000);
\draw[very thick, red!70!black] (0.0000,1.0000) -- (0.1667,1.0000);
\draw[very thick, red!70!black] (0.1667,1.0000) -- (0.3333,1.0000);
\draw[very thick, red!70!black] (0.3333,1.0000) -- (0.5000,1.0000);
\draw[very thick, red!70!black] (0.5000,1.0000) -- (0.6667,1.0000);
\draw[very thick, red!70!black] (0.6667,1.0000) -- (0.8333,1.0000);
\draw[very thick, red!70!black] (0.8333,1.0000) -- (1.0000,1.0000);
\draw[very thick, red!70!black] (0.0000,0.5000) -- (0.0000,0.6667);
\draw[very thick, red!70!black] (0.0000,0.6667) -- (0.0000,0.8333);
\draw[very thick, red!70!black] (0.0000,0.8333) -- (0.0000,1.0000);
\draw[very thick, red!70!black] (0.5000,0.6667) -- (0.5000,0.8333);
\draw[very thick, red!70!black] (0.5000,0.8333) -- (0.5000,1.0000);
\draw[very thick, red!70!black] (1.0000,0.5000) -- (1.0000,0.6667);
\draw[very thick, red!70!black] (1.0000,0.6667) -- (1.0000,0.8333);
\draw[very thick, red!70!black] (1.0000,0.8333) -- (1.0000,1.0000);
\fill[black] (0.0000,0.5000) circle (0.9pt);
\fill[black] (0.3333,0.5000) circle (0.9pt);
\fill[black] (0.5000,0.0000) circle (0.9pt);
\fill[black] (0.5000,0.3333) circle (0.9pt);
\fill[black] (0.5000,0.5000) circle (0.9pt);
\fill[black] (0.5000,0.6667) circle (0.9pt);
\fill[black] (0.6667,0.5000) circle (0.9pt);
\fill[black] (1.0000,0.5000) circle (0.9pt);
\fill[red!70!black] (0.5000,1.0000) circle (0.9pt);
\end{tikzpicture}\\[2pt]
    {\small Patch type C}
  \end{minipage}
  \caption{Example of a polytopal mesh $\mathcal{P}$ and induced star patches around selected vertices of $\mathcal{N}(\mathcal{P})$. Type A: a vertex shared by singleton polytopes only, whose patch is the classical vertex star. Type B: a vertex on the boundary of a non-trivial polytope, whose patch contains the whole polytope together with the neighbouring singleton cells. Type C: a vertex on the boundary of the domain, where the boundary faces incident to the vertex (in red) carry no Dirichlet condition.}
  \label{fig:star-types}
\end{figure}

In \fig{fig:star-types}, we show an example of this construction for a polytopal mesh $\mathcal{P}$ obtained by agglomerating the cells of a Cartesian mesh $\mathcal{T}$, chosen so that the colouring contains both trivial and non-trivial agglomerates (similar to $\Taggall$). The leftmost image shows the polytopal mesh $\mathcal{P}$, with faces and vertices highlighted. The remaining images show different types of star patches around selected vertices of $\mathcal{N}(\mathcal{P})$. The generating vertex and its incident faces are painted in red, and the homogeneous Dirichlet boundary in blue.

We now turn to the definition of the smoother induced by the star patches. Let $A : V \to V'$ be a linear operator representing the global system. For each star patch $\Star{\nu}$, we define the local operator $A_\nu : V_\nu \to (V_\nu)'$ as the restriction of $A$ to the local space $V_\nu$, that is, $A_\nu = I_\nu^T A I_\nu$, where $I_\nu : V_\nu \to V$ is the injection operator and $I_\nu^T : V' \to (V_\nu)'$ is its adjoint. We can then define the associated additive Schwarz preconditioner as
\[
  D(V,A,\mathcal{X}(\mathcal{P})) \coloneqq \sum_{\nu \in \mathcal{X}(\mathcal{P})} I_\nu A_\nu^{-1} I_\nu^T.
\]
In the notation of \eqref{eq:effective-smoother}, $D(V,A,\mathcal{X}(\mathcal{P}))$ plays the role of the smoothing operator $\Dl$.

We now introduce three distinct smoothers based on the above construction of star patches. Although patches around faces (and edges in three dimensions) could also be considered, we will limit ourselves to vertex-star patches. The main motivation is that they work for all spaces in the de Rham complex, and have been shown to be $p$-robust for body-fitted \ac{fe} methods~\cite{Pavarino1994,Schoberl2008,Miraci2021,BrubeckFarrell2024}.
Let $\Tl$, $\Tlaggall$ and $\Tlagg$ be the meshes defined on level $\ell$, with the associated space $\Vl$ and operator $\Al$. Then, define
\begin{align}
  D_\ell^{\rm act} & \coloneqq D(\Vl,\Al,\mathcal{N}(\Tl)) \label{def:smoother-classical} \\
  D_{\ell}^{\rm agg, all} &\coloneqq D(\Vl,\Al,\mathcal{N}(\Tlaggall)) \label{def:smoother-agg-all} \\
  D_{\ell}^{\rm agg} &\coloneqq D(\Vl,\Al,\tilde{\mathcal{N}}(\Tlaggall)), \label{def:smoother-agg}
\end{align}
where
\[
  \tilde{\mathcal{N}}(\Tlaggall) \coloneqq \{ z \in \mathcal{N}(\Tlaggall) : z \in \ol{P} \text{ for some } P \in \Tlagg \},
\]
that is, the set of vertices of $\Tlaggall$ that belong to some non-trivial aggregate in $\Tlagg$.
$D_\ell^{\rm act}$ is the classical vertex-star smoother on the active mesh $\Tl$. All its patches are type A as illustrated in \fig{fig:star-types}, except on the active mesh boundary where the patches are even smaller. $D_{\ell}^{\rm agg, all}$ is the polytopal vertex-star smoother on the aggregated mesh $\Tlaggall$, containing all patch types illustrated in \fig{fig:star-types}. Finally, $D_{\ell}^{\rm agg}$ is the polytopal vertex-star smoother restricted to the vertices of non-trivial aggregates, containing only patches of types B and C as illustrated in \fig{fig:star-types}. The idea is that the patches cover $\Tlagg$ as well as a single layer of interior active cells surrounding the aggregates.  The smoother $D_{\ell}^{\rm agg}$ is only positive semidefinite, since its patches do not cover the \acp{dof} away from the aggregates. It is used only within the auxiliary-space preconditioner of \sect{sec:auxspace}, where the active problem takes care of those.

\begin{remark}[Bounded overlap] \label{rem:bounded-overlap}
  A cell $T \in \mathcal{T}$ belongs to the star patch $\Star{\nu}$ if and only if $\nu$ lies on the closure of the polytope $P \in \mathcal{P}$ containing $T$. Thus, the maximum overlap $N_{\rm ov}$ is the maximum number of entities of $\mathcal{X}(\mathcal{P})$ lying on the closure of a single polytope. For the aggregated meshes, the aggregates contain a bounded number of cells (see \sect{subsec:aggregates}), so  $N_{\rm ov}$ is bounded independently of the mesh size and of the cut configuration. The overlap also governs the admissible damping of the smoother, since by the standard colouring argument $\lambda_{\max}(\Dl \Al) \le N_{\rm ov}$, so that the Richardson iteration \eqref{eq:effective-smoother} is stable for $\theta_\ell < 2/N_{\rm ov}$. In the experiments of \sect{sec:results}, $N_{\rm ov}$ is of the order of ten in two dimensions and of a few tens in three dimensions, which motivates the conservative damping parameters used there. This condition makes the Richardson iteration a contraction, which the V-cycle requires as a fixed-point iteration, but it is not required when the V-cycle is used as a preconditioner for a Krylov method, as in \sect{sec:results}, where only the conditioning of the preconditioned operator matters. Indeed, the tuned runs of \sect{subsec:results-gmg} violate $\theta_\ell < 2/N_{\rm ov}$ by a large margin.

  Furthermore, since the stabilisation only couples \acp{dof} within the same aggregate, the local spaces $V_\nu$ and $V_{\nu'}$ are orthogonal with respect to the bilinear forms of \sect{sec:formulations} unless $\Star{\nu}$ and $\Star{\nu'}$ overlap, and $N_{\rm ov}$ therefore also bounds the overlap in the sense of the strengthened Cauchy--Schwarz inequality of the additive Schwarz theory \cite[Sect.~2.3]{ToselliWidlund2005}.
\end{remark}

\begin{remark}[Computational cost]
  By enlarging the support of the star patches to include entire aggregates, we are increasing the number of \acp{dof} involved in each local solve, and thus the computational cost of the smoother. However, only patches associated with vertices of non-trivial aggregates are enlarged. These patches are therefore concentrated in a $(d-1)$-dimensional subset of the domain, and the overall impact on the computational cost remains manageable. The proportion of enlarged patches in $D_{\ell}^{\rm agg, all}$ tends to zero as the mesh is refined.
\end{remark}

\begin{remark}[High-order high-dimensional strategies]
  In high-order and high-dimensional settings, the number of \acp{dof} per patch can grow quickly, even for the classical vertex-star smoother $D_\ell^{\rm act}$. A dense factorisation of each local matrix then costs $\mathcal{O}(k^{3d})$ operations and $\mathcal{O}(k^{2d})$ storage per patch, which quickly becomes the dominant cost of the smoother. The problem is well known for body-fitted methods, and several remedies have been proposed: static condensation of the cell-interior \acp{dof}, which leaves patch problems posed on the interfaces only \cite{BrubeckFarrell2024,Dobrev2019}; sparsity-friendly bases and tensor-product fast diagonalisation on hexahedral cells \cite{BrubeckFarrell2022,Witte2021}; finer space decompositions with smaller patches \cite{Schoberl2008,BrubeckFarrell2024}; and inexact local solvers \cite{Wichrowski2025}.
  Of these, static condensation is the only one that is purely algebraic and does not rely on the tensor-product structure of the cells, which the cut cells and the aggregates lack. Static condensation is therefore the option implemented in our code for the high-order patch solves.
\end{remark}

\section{Homogeneous geometric multigrid}\label{sec:embedded-gmg}

We now use the abstract multigrid framework in \sect{sec:gmg} to construct a \emph{homogeneous} \ac{gmg} preconditioner for the unfitted operator $\Ah$, in the sense that the unfitted discretisation is used on every level of the hierarchy~\cite{lu2022}. With the spaces and transfer operators of \sect{subsec:spaces-transfer} fixed, the two remaining ingredients are the level operators and the smoother.

For each level $\ell$, we repeat the construction in \sect{sec:formulations} on the aggregated mesh $\Tlagg$ and define on $\Vl$ the unfitted bilinear form
\begin{equation}\label{eq:unfitted-level-form}
  \Ll(u,v) \coloneqq \Lx(u,v; \dom) + \tau^2 S^{\rm x}_\ell(u,v),
\end{equation}
where $S^{\rm x}_\ell$ is the stabilisation \eqref{eq:stabilisation-biforms} built from the aggregate-wise projections on $\Tlagg$. We denote by $\Al$ the operator induced by $\Ll$ on $\Vl$, so that $A_L^{\rm x} = \Ah$. By the inherited boundary, the consistency terms are integrated over the same physical domain $\dom$ on every level. The levels differ only through the \ac{fe} space and the stabilisation.

Finally, we take $\Dl = D_{\ell}^{\rm agg, all}$ of \eqref{def:smoother-agg-all} as the local smoothing operator for the multigrid method. As discussed in \sect{subsec:smoothers}, this choice enlarges the support of the subdomains in the cut region to correctly capture the kernel of $S^{\rm x}_\ell$, while falling back to the classical vertex-star smoother in regions away from the aggregates.

The homogeneous multigrid preconditioner for $\Ah$ is then
\begin{equation}\tag{P1}
  \Bgmg \coloneqq \VCycle\big(\{\Al\}_{\ell}, \{\Dl\}_{\ell}\big).
\end{equation}
The goal is a preconditioner such that the condition number of $\Bgmg \Ah$ is bounded independently of the mesh size, the number of levels and the cut configuration.

\section{Auxiliary-space geometric multigrid}\label{sec:auxspace}

 Recall from \sect{sec:formulations} that $\Ah, \Aact : \Vh \to \Vhd$ denote the operators induced by the unfitted and active bilinear forms, respectively. The objective of this section is to design a multigrid solver for the unfitted operator $\Ah$ which does not require rediscretising the unfitted problem on each level. Instead, we aim to use the active operator $\Aact$ and its associated hierarchy to construct an efficient multigrid method.

The first step is to recognise that the active operator $\Aact$ already provides a robust preconditioner for the unfitted operator $\Ah$. Indeed, the stabilisation in unfitted discretisations is designed such that property \eqref{eq:stab-property} holds uniformly with respect to the cut configuration and the mesh size.

Thus, a first valid approach is to use the active Riesz map as a preconditioner for the unfitted operator, that is
\begin{equation}\tag{P2}
  \Bx{ac} \coloneqq (\Aact)^{-1}.
\end{equation}
With this choice, the preconditioned operator $\Bx{ac} \Ah$ has its spectrum bounded independently of the cut configuration and the mesh size by the constants appearing in property \eqref{eq:stab-property}. Unfortunately, these constants can be large, particularly for $\Div$- and $\Curl$-conforming discretisations, leading to slow convergence of the preconditioned solver.

\begin{figure}
  \centering
  \includegraphics[width=0.23\textwidth]{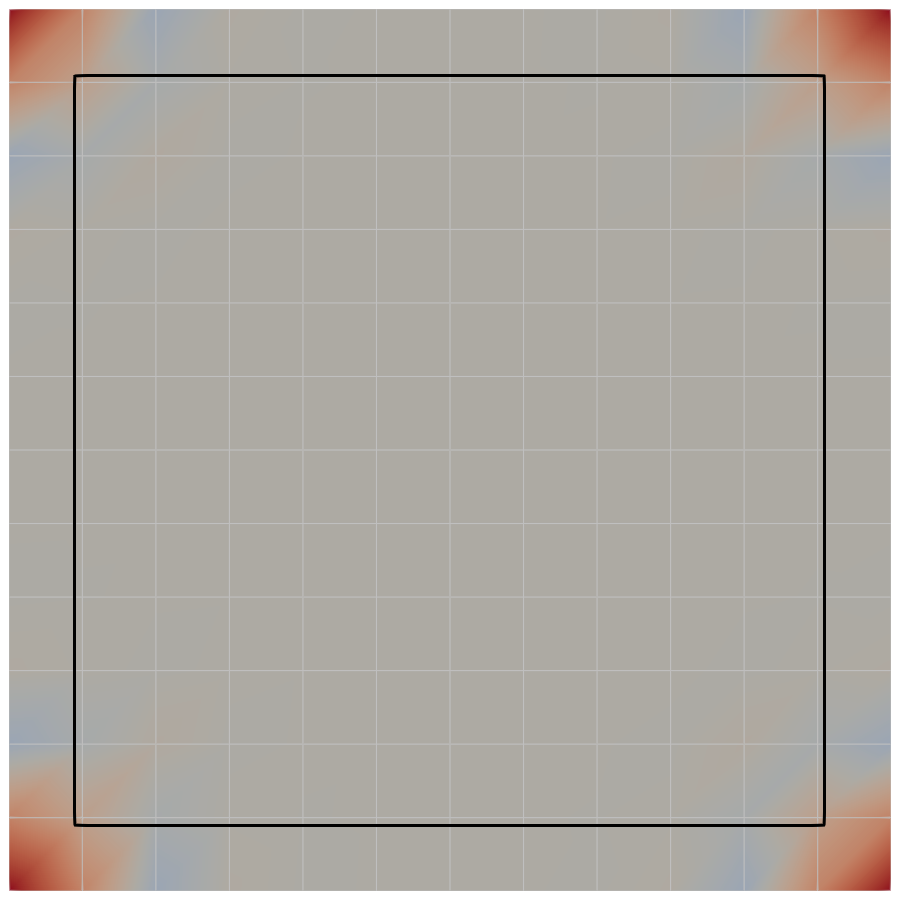}\hfill
  \includegraphics[width=0.23\textwidth]{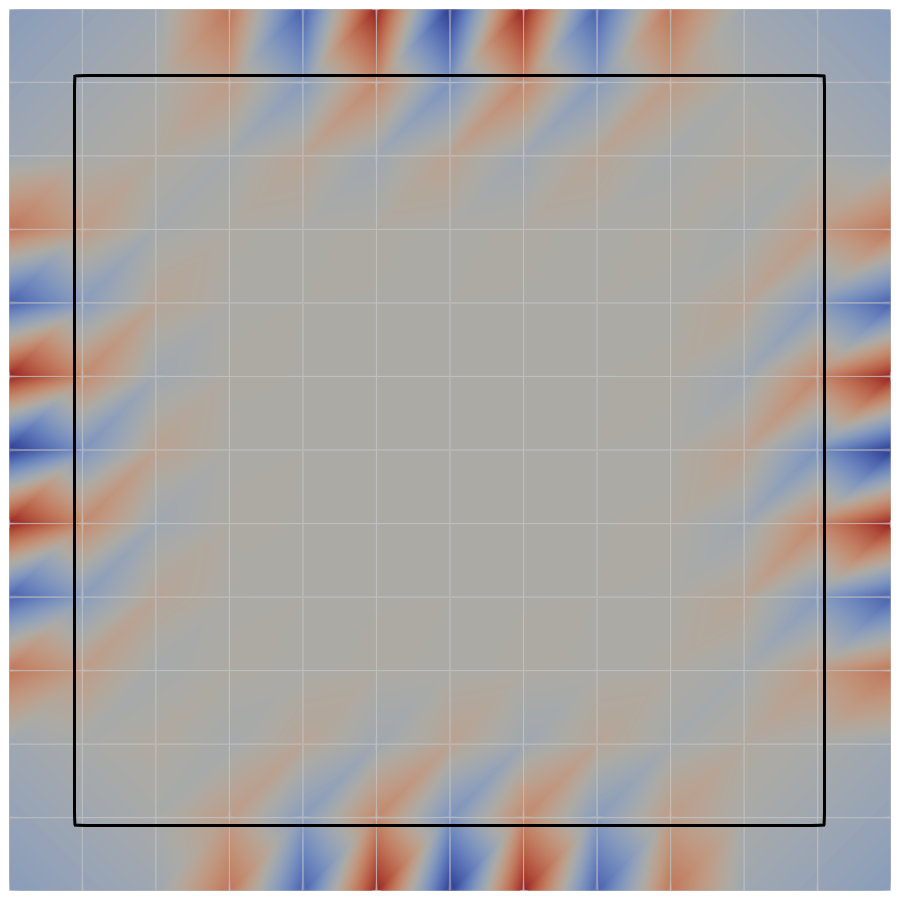}\hfill
  \includegraphics[width=0.23\textwidth]{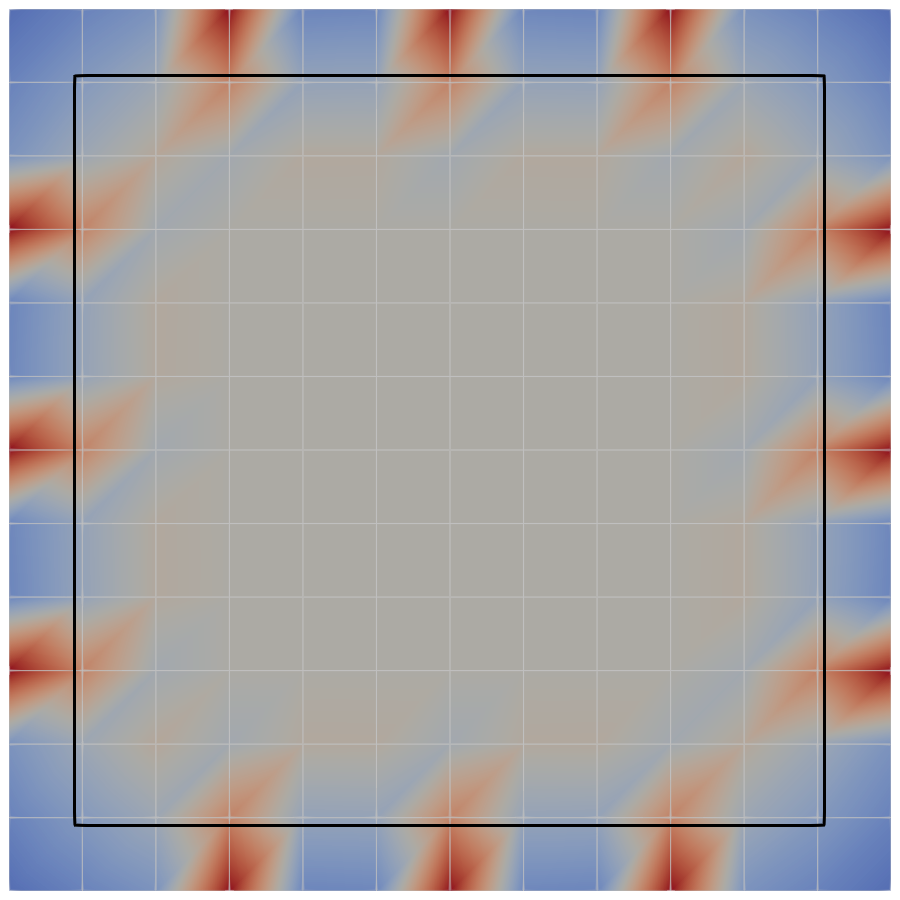}\hfill
  \includegraphics[width=0.23\textwidth]{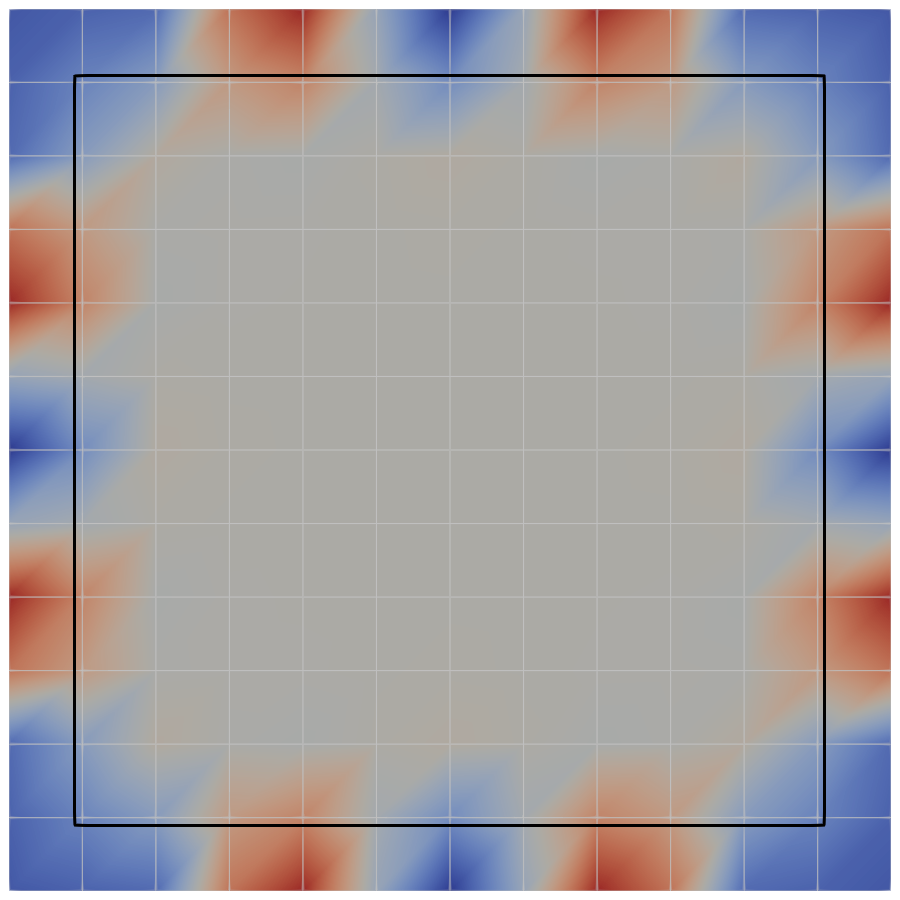}\\[4pt]
  \includegraphics[width=0.23\textwidth]{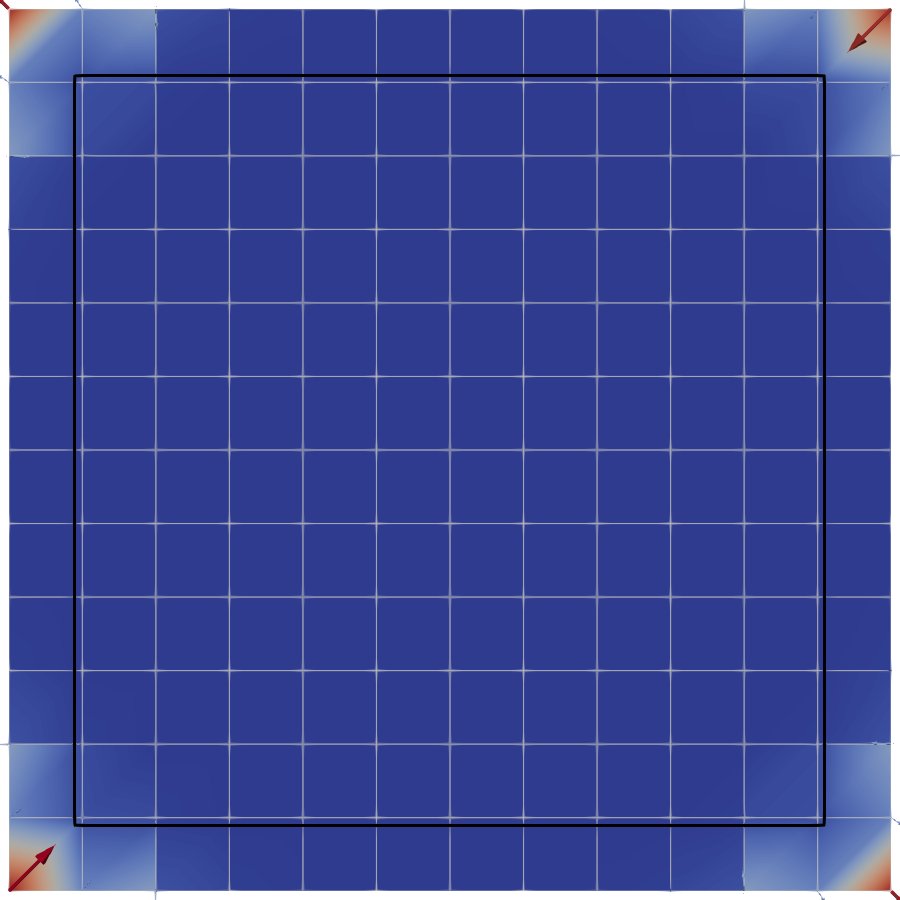}\hfill
  \includegraphics[width=0.23\textwidth]{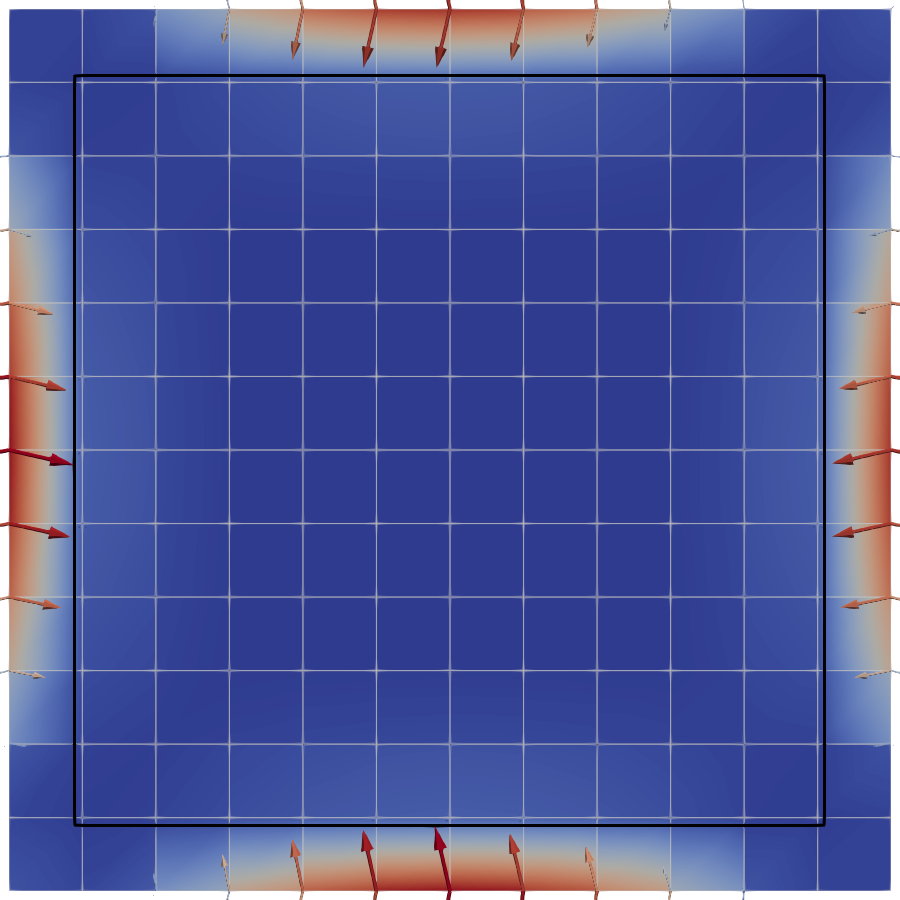}\hfill
  \includegraphics[width=0.23\textwidth]{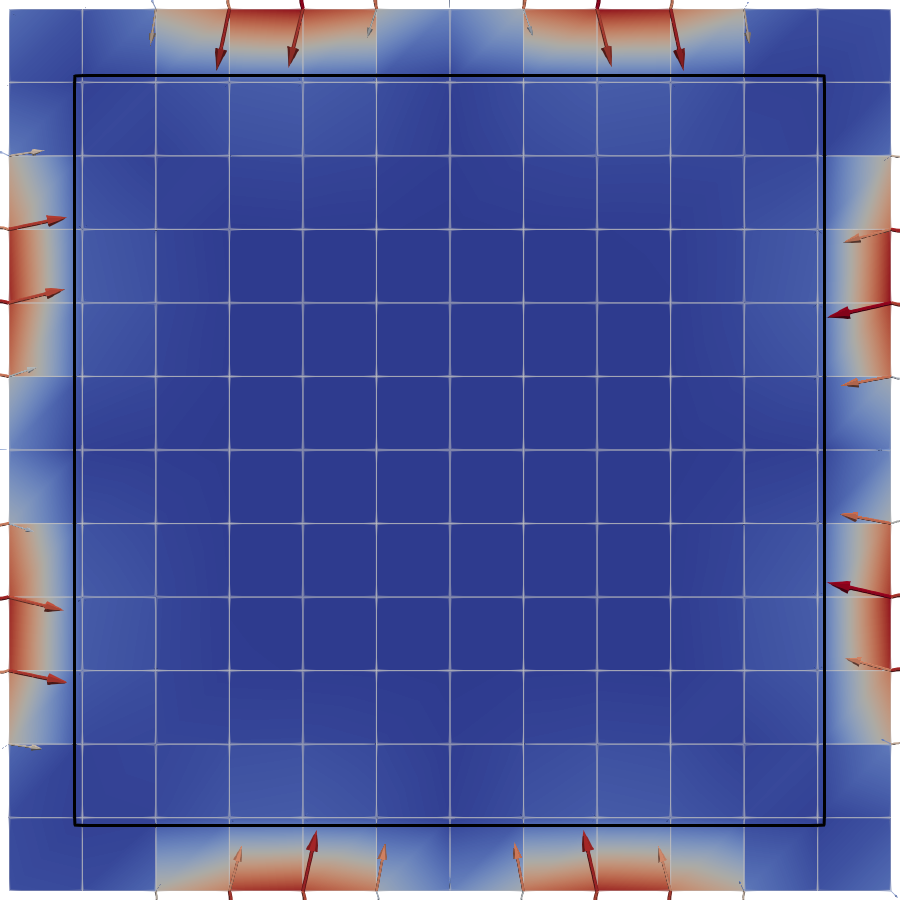}\hfill
  \includegraphics[width=0.23\textwidth]{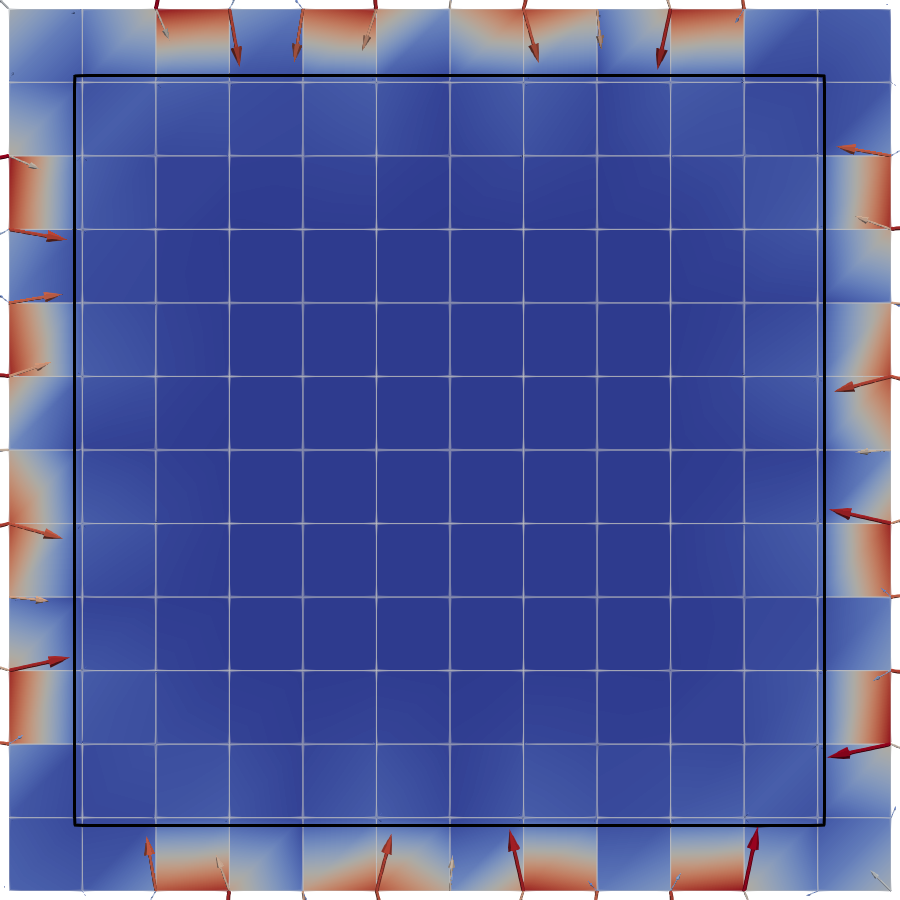}
  \caption{Worst-controlled modes of $(\Aact)^{-1}\Ah$ on the cut square. Top: $H(\Grad)$
  (scalar value, diverging colour map). Bottom: $H(\Div)$ (flux magnitude and arrows).
  All modes live in the fictitious cut layer between $\bnd$ and the active mesh boundary,
  concentrated at the corners or oscillating along the boundary layer.}
  \label{fig:modes}
\end{figure}

 We examine the eigenvectors of the generalised eigenvalue problem $\Ah x = \lambda \Aact x$ with the smallest eigenvalues, which correspond to the modes that are poorly controlled by $\Bx{ac}$. \fig{fig:modes} shows a selection of these modes on a cut-square domain.
 These modes are \emph{localised on the cut cell layer}. This is the region controlled by the aggregation stabilisation, and where $\Aact$ and $\Ah$ disagree. Thanks to their local support, a cheap solver acting only on the cut-cell aggregate patches can resolve them, without touching the global modes that $\Aact$ already handles well.

To recover robustness for the modes $\Aact$ misses, we propose an auxiliary-space preconditioner \cite{Xu1996} that combines the active problem with a smoother $\Dh$ localised around the cut cell region. This is a conceptually similar approach to \cite{gross2023cutfemprec}, but instead of relying on a splitting of the \acp{dof} we keep the original problems and rely on the smoother to handle the problematic modes.

We now address the choice of $\Dh$. The rationale follows the same principles as for the homogeneous multigrid solver of \sect{sec:embedded-gmg}. Thus, we choose $\Dh \coloneqq D_{L}^{\rm agg}$ of \eqref{def:smoother-agg} as our smoother, where $L$ denotes the finest level, so that $V_L = \Vh$. Choosing $D_{L}^{\rm agg, all}$ would also lead to a robust method, but the patches that do not intersect an aggregate act on modes that $(\Aact)^{-1}$ already resolves, and can therefore be omitted.

Let $\Rh$ be the effective smoother derived from $\Dh$ as in \eqref{eq:effective-smoother}, that is, $n_L$ steps of the damped Richardson iteration with $\Dh$ applied to $\Ah$.
Since the auxiliary space is the active space itself, the intergrid transfers are the identity and, in additive form, the auxiliary-space preconditioner reads
\begin{equation}\label{eq:auxspace-add-precond}\tag{P3}
  \Bx{add} \coloneqq \Rh + (\Aact)^{-1}.
\end{equation}
Since $\Rh$ is symmetric (see \sect{subsec:vcycle}), so is $\Bx{add}$.

We also define a multiplicative version of the preconditioner as a single V-cycle iteration, namely the two-level instance of \eqref{eq:error-propagation} with identity transfers in which the coarse solve is the active problem,
\begin{equation}\label{eq:auxspace-mult-precond}\tag{P4}
  \Id - \Bx{mult} \Ah \coloneqq (\Id - \Rh^T \Ah)(\Id - (\Aact)^{-1} \Ah)(\Id - \Rh \Ah).
\end{equation}
$\Bx{mult}$ is symmetric and, for $\theta_L$ small enough, positive definite, and is therefore a valid preconditioner. As a stationary iteration, however, it need not be convergent, since $\Id - (\Aact)^{-1}\Ah$ is not a contraction when $\lambda_{\max}((\Aact)^{-1}\Ah) > 2$, which may well happen for the values of $\tau$ used in practice.

As we will see in the numerical experiments, both the additive and multiplicative auxiliary-space preconditioners exhibit robust performance with respect to the cut configuration and the mesh size, and significantly improve the convergence of iterative solvers for the unfitted problem. The remaining challenge is the efficient solution of the active problem, which is still a large system.

To address this final hurdle, we replace the exact solve of the active problem in the auxiliary-space preconditioner with a \ac{gmg} V-cycle for the active problem. The key difficulty lies in the fact that the active triangulation $\Th$ has fine-level features and thus cannot be coarsened in a straightforward manner without losing accuracy near the fine-scale structures. In \cite{gross2023cutfemprec} the analogous interior problem is left to any spectrally equivalent solver, and algebraic multigrid is used in practice. However, algebraic multigrid performs poorly on $H(\Div)$- and $H(\Curl)$-conforming discretisations unless combined with auxiliary-space constructions of their own \cite{HiptmairXu2007}. We prefer a geometric construction that treats the three spaces alike.

Instead, we use a \ac{gmg} on the hierarchy of active meshes $\{\Tl\}_{\ell}$ of \sect{subsec:hierarchy}, that is, the cells of each background mesh $\Tlbg$ that meet $\dom$. The active meshes $\Tl$, and hence the spaces $\Vl$, are not nested, but the transfer operators of \sect{subsec:spaces-transfer} account for this, and we only require the background mesh to accommodate a hierarchy, which, as discussed, is easy to accomplish.
Furthermore, the construction exploits the fact that $\Aact$ is a \emph{body-fitted} Riesz map on the active domain: unlike the unfitted operator $\Ah$, it carries no cut-cell stabilisation and is therefore amenable to a textbook multigrid such as the ones in \cite{Arnold2000}.

On the level spaces $\Vl$ of \sect{subsec:spaces-transfer}, we define the active bilinear form
\[
  \Llact(u,v) \coloneqq \Lambda^{\rm x}(u,v; \doml).
\]
We denote by $\Alact$ the operator induced by $\Llact$ on $\Vl$. The prolongation $\Il$ and restriction $\Il^T$ are the common ones of \eqref{eq:prolongation}.

For the choice of smoother, we can draw from the extensive literature on multigrid methods for standard Riesz maps. In this work we take the classical vertex-star smoother of \sect{subsec:smoothers} on the active mesh, that is, $\Dlact \coloneqq D(\Vl,\Alact,\mathcal{N}(\Tl))$, the vertex-star smoother \eqref{def:smoother-classical} with the active operator in place of the unfitted one. The motivation for this choice is that the vertex-star smoother is known to be robust for all the considered discretisations~\cite{Arnold2000,farrell2021}.

With these ingredients, the active multigrid preconditioner is
\begin{equation}\label{eq:active-gmg-precond}
  \Bagmg \coloneqq \VCycle\big(\{\Alact\}_{\ell}, \{\Dlact\}_{\ell}\big),
\end{equation}
and the corresponding preconditioner for $\Ah$ is obtained by substituting the exact solve in \eqref{eq:auxspace-add-precond} and \eqref{eq:auxspace-mult-precond} with $\Bagmg$. In what follows we only use the multiplicative variant, and we denote by $\Bmg$ the resulting \emph{multilevel auxiliary-space preconditioner}, defined by
\begin{equation}\tag{P5}
  \Id - \Bmg \Ah \coloneqq (\Id - \Rh^T \Ah)(\Id - \Bagmg \Ah)(\Id - \Rh \Ah).
\end{equation}
It is the second multigrid preconditioner of this paper.

\begin{remark}[Comparison with \cite{gross2023cutfemprec}]
  In \cite{gross2023cutfemprec}, the authors present a non-overlapping domain decomposition preconditioner for Poisson problems. The main idea is to split the original problem into a solve of the interior problem on $\indom$ with homogeneous Dirichlet boundary conditions and a solve on the span of the cut basis functions, whose matrix is uniformly well conditioned after diagonal scaling, effectively decoupling the two regions. In the spirit of \eqref{eq:auxspace-add-precond}, one could extend the construction to other discretisations by considering
  \[
    \Bx{split} \coloneqq D_{\rm cut} + E_{\rm in}\, A_{\rm in, 0}^{-1} E_{\rm in}^T,
  \]
  where $D_{\rm cut}$ is a smoother on the cut domain, $A_{\rm in, 0}^{-1}$ is the solve (in practice, by any optimal solver) of the interior problem on $\indom$ with homogeneous Dirichlet boundary conditions and $E_{\rm in}$ is the extension by zero from the space with homogeneous Dirichlet conditions on $\partial\indom$ to $\Vh$.

  From our exploratory experiments, we observe that for $\Hgrad{\dom}$ problems, this preconditioner is only $h$-robust when considering Dirichlet boundary conditions, that is, for $\widehat{\Lgradh}$ when $\bnd = \bnddir$. For $\Hgrad{\dom}$ with Neumann boundary conditions, and for $\Hdiv{\dom}$ and $\Hcurl{\dom}$ problems with any boundary condition, we found that the preconditioner is not $h$-robust.
\end{remark}

\begin{remark}[Convergence of $\Bx{add}$]
  The convergence of $\Bx{add}$ is covered by the standard theory of additive Schwarz methods \cite{ToselliWidlund2005,Xu1992}, for the space decomposition $\Vh = \Vh + \sum_\nu V_\nu$. Its condition number is bounded in terms of two constants: the maximum overlap of the subspaces, bounded by Remark~\ref{rem:bounded-overlap}, and a stable decomposition $u = w + \sum_\nu w_\nu$ of every $u \in \Vh$ into $w \in \Vh$ and $w_\nu \in V_\nu$. The trivial decomposition $w = u$, $w_\nu = 0$ then satisfies the stability requirement with constants that can be directly related to the constants in \eqref{eq:stab-property}.

  This bound is not sharp, since, as the experiments of \sect{sec:results} show, $\Bx{add}$ performs much better than $\Bx{ac}$ in practice. Proving sharper bounds would require similar tools as the analysis of the multigrid methods, and we leave it as future work.
\end{remark}

\section{Numerical results}\label{sec:results}

All the methods in this paper have been implemented in the Julia programming language, on top of the \gridap{} ecosystem. The \gridap{} library \cite{Badia2020gridap,Verdugo2022gridap} provides the core \ac{fe} machinery: meshes, finite element spaces of all conformities, and assembly. \gridapembedded{} extends it with the tools required by unfitted methods: level-set cutting, classification of cells and facets, cut-cell quadratures and cell aggregation. Finally, \gridapsolvers{} \cite{Manyer2024gridapsolvers} provides the linear solver infrastructure, including Krylov solvers, the \ac{gmg} framework and patch-based smoothers.

Every geometry considered below is unfitted in the unit box $[0,1]^d$, $d \in \{2,3\}$, which is discretised by a Cartesian background mesh $\Tbg$ made of $n^d$ cells. The domain is described by a level-set function. Cells are cut by the zero level-set using a marching-cubes-type algorithm, so that the discrete boundary is piecewise linear, and quadratures on the cut cells are moment-based and computed from a simplicial sub-triangulation of the cut regions~\cite{Badia-JComputPhys-2022k}.

Besides robustness with respect to the mesh size, unfitted discretisations require robustness with respect to the cut configuration. We quantify it in terms of the \emph{cut ratio} $\eps \in (0,1)$, defined as the interior volume fraction $|T \cap \dom| / |T|$ of a cell $T$ cut by a single flat piece of the boundary $\bnd$ parallel to a grid line of $\Tbg$, as in \cite{badia2026darcy}. Equivalently, the flat piece lies at a distance $\eps h$ from the grid line on its interior side, with corner cells of the domain possibly having fractions of the order of $\eps^d$.

\begin{figure}[ht]
  \centering
  \includegraphics[width=0.19\textwidth]{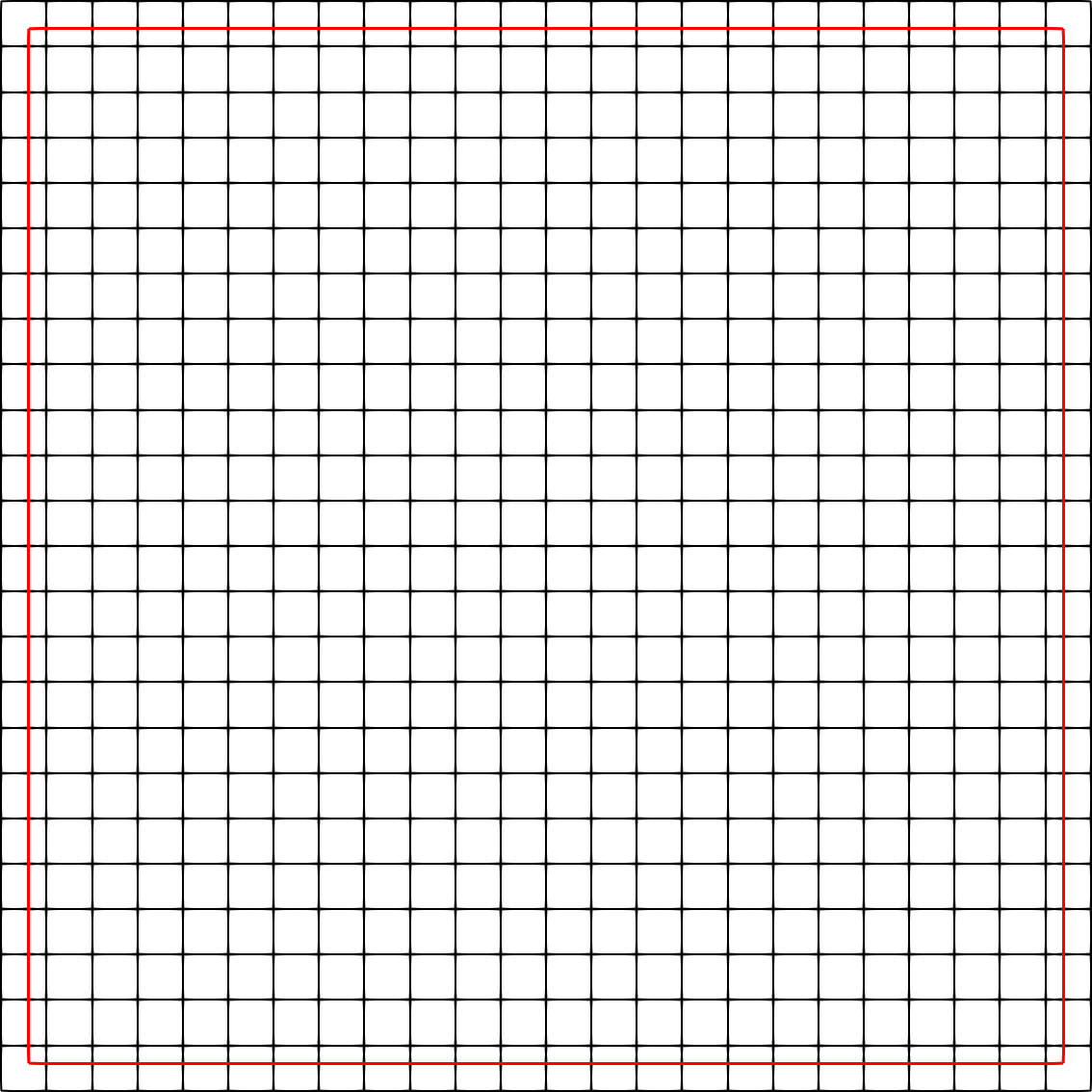}\hfill
  \includegraphics[width=0.19\textwidth]{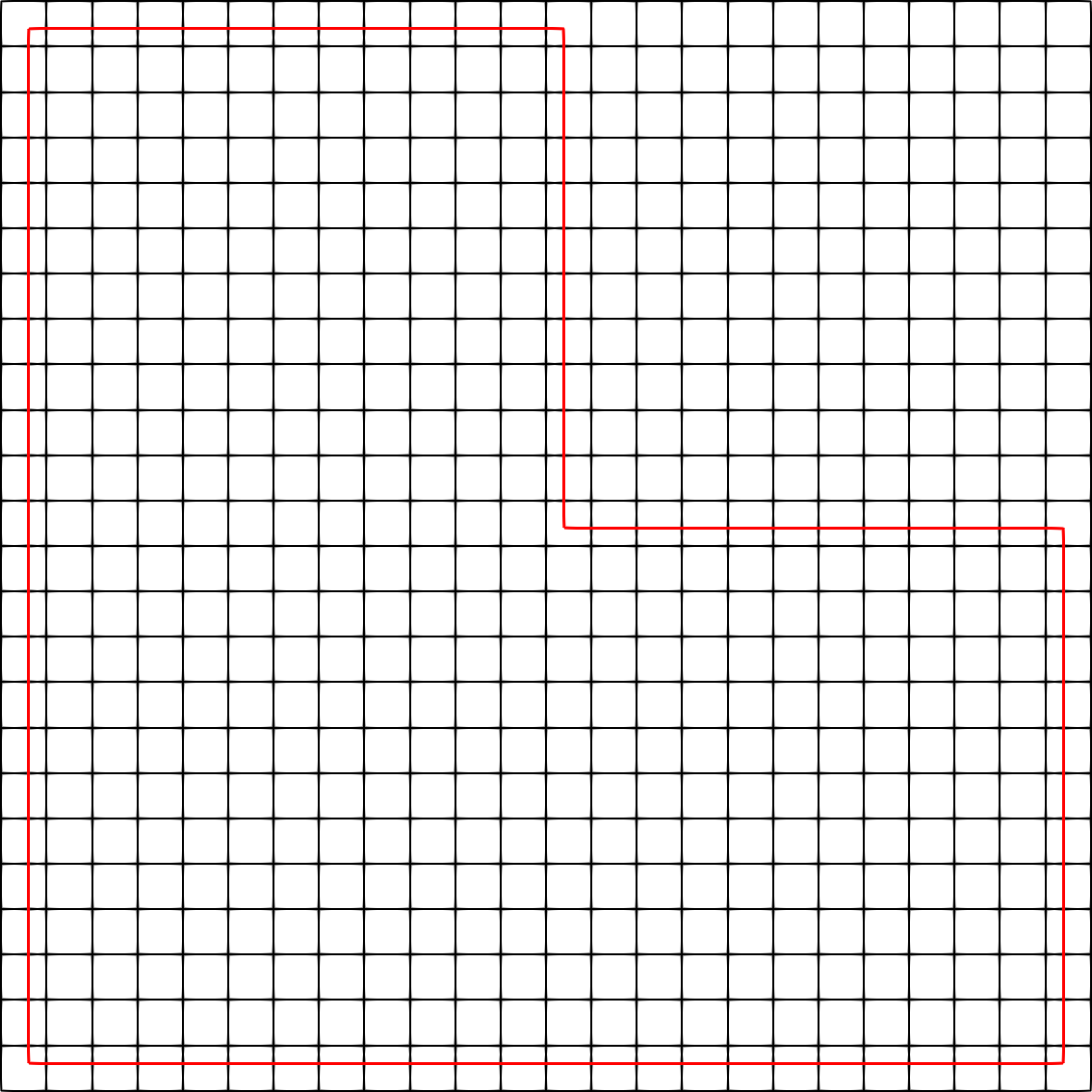}\hfill
  \includegraphics[width=0.19\textwidth]{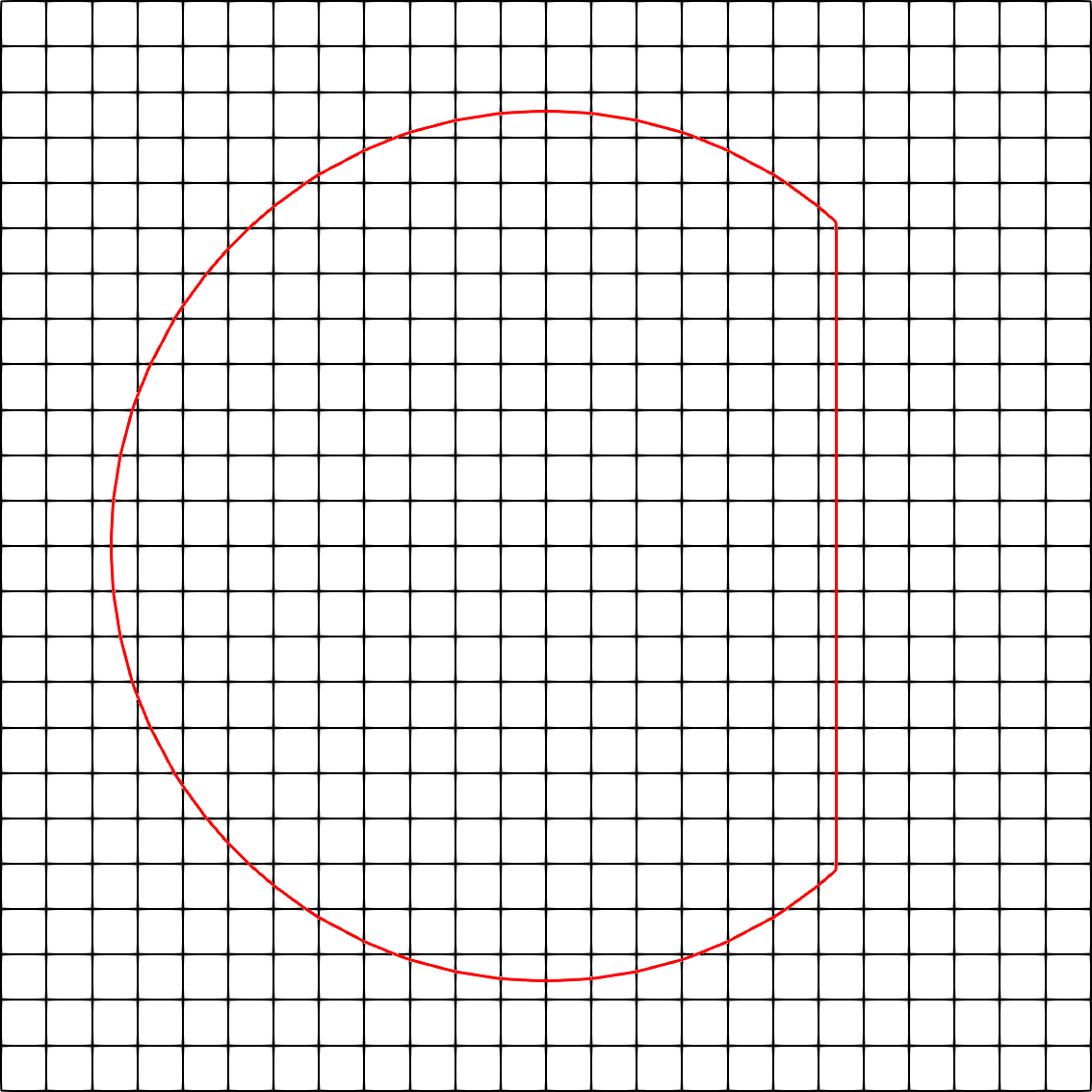}\hfill
  \includegraphics[width=0.19\textwidth]{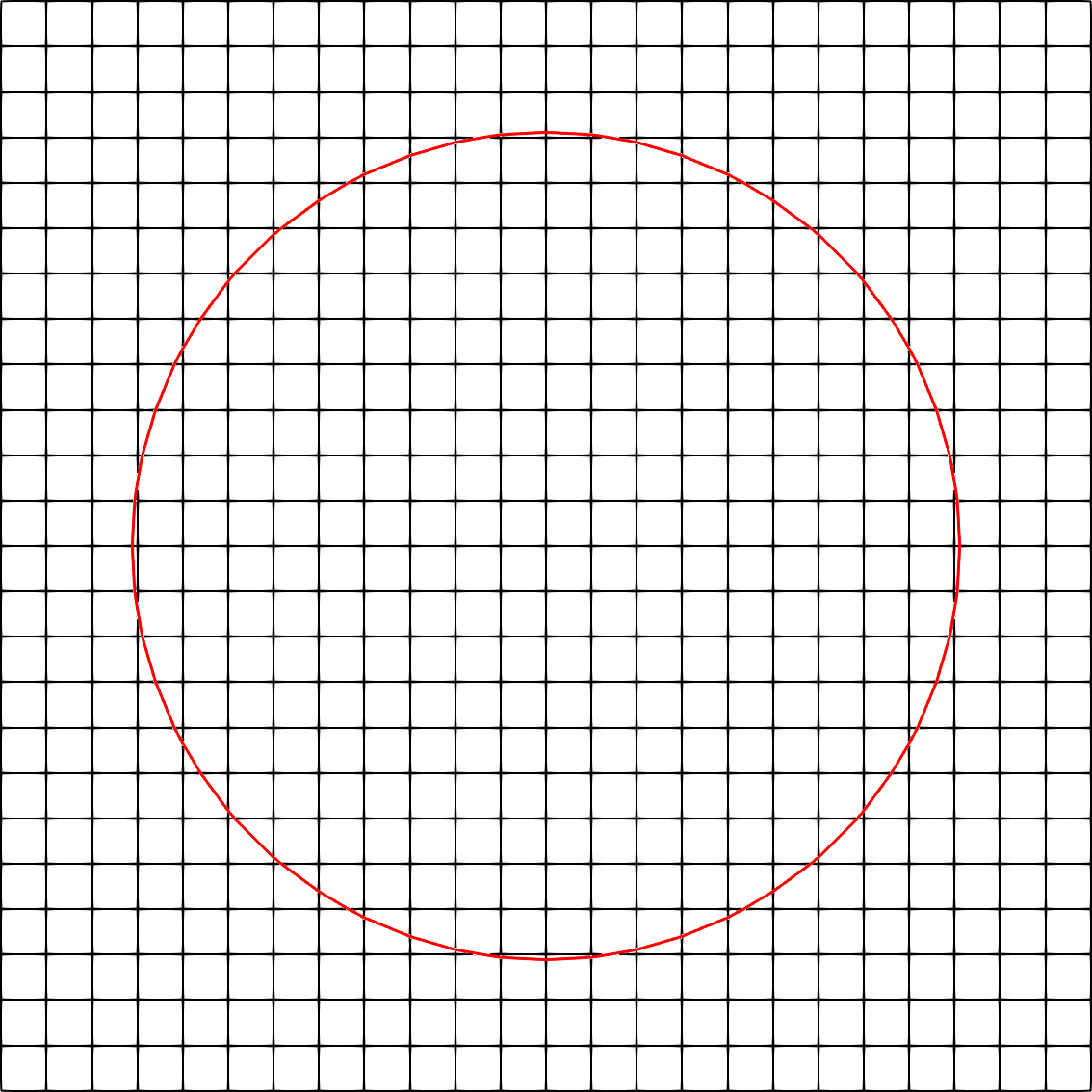}\hfill
  \includegraphics[width=0.19\textwidth]{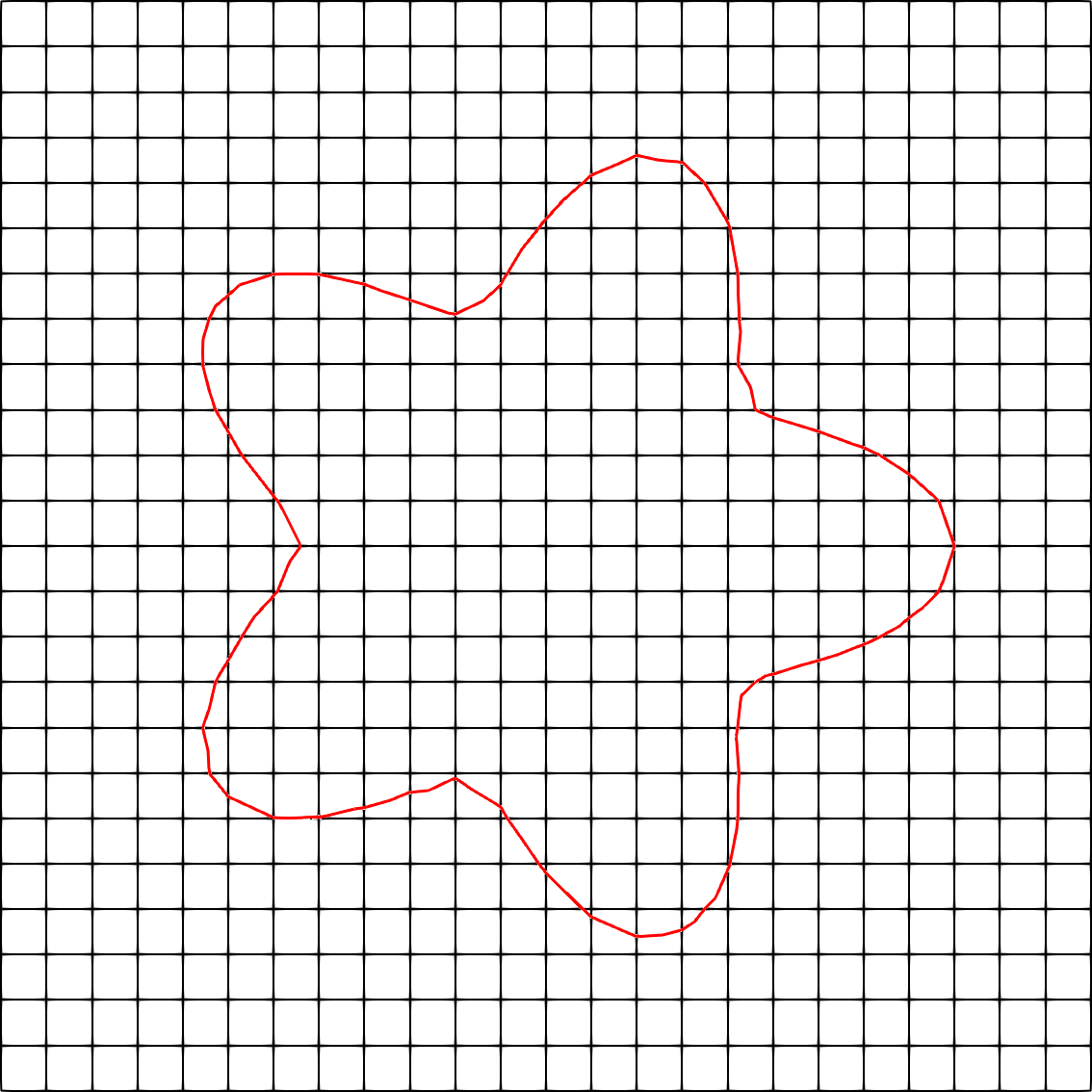}\\[6pt]
  \includegraphics[width=0.19\textwidth]{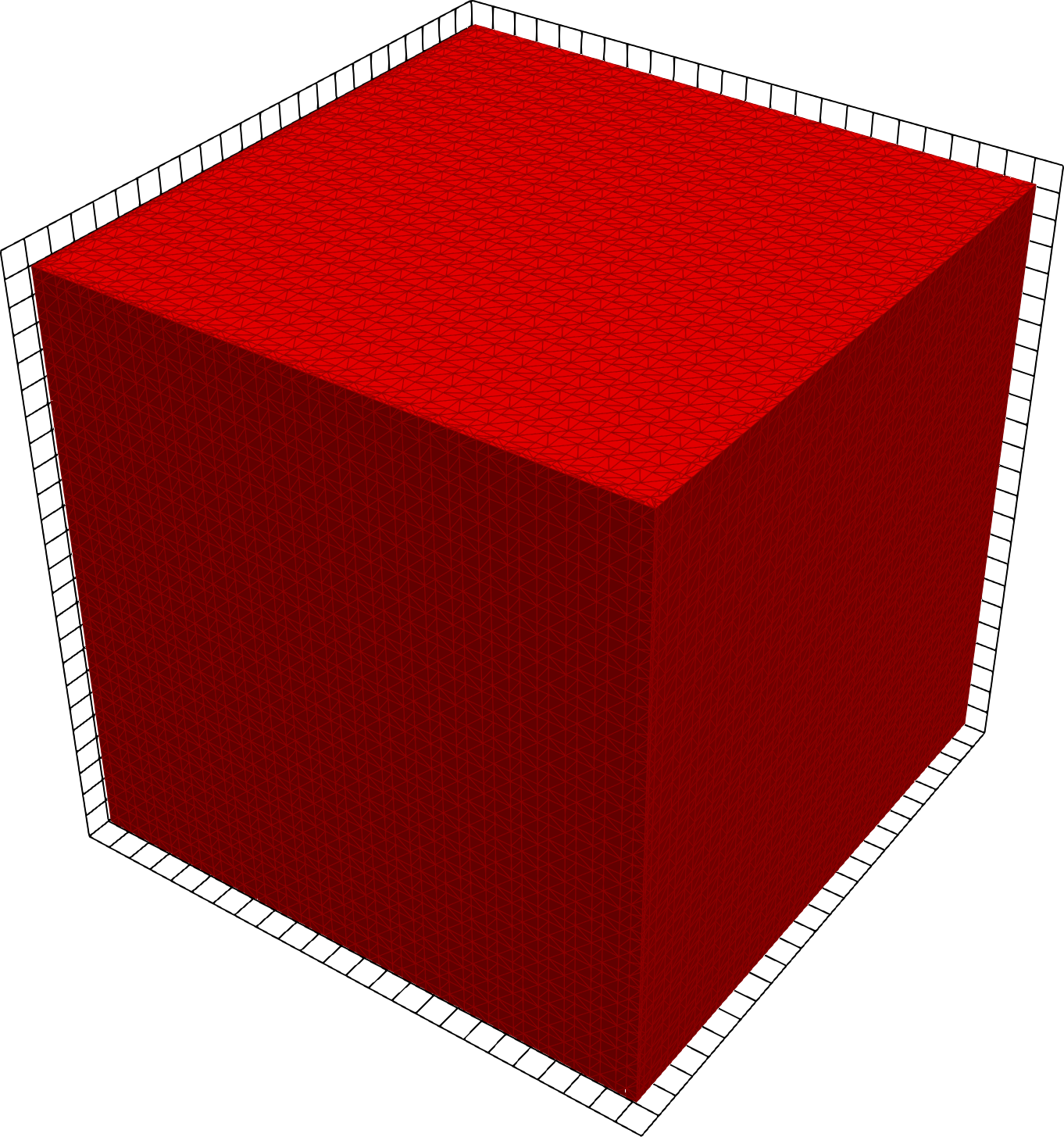}\hfill
  \includegraphics[width=0.19\textwidth]{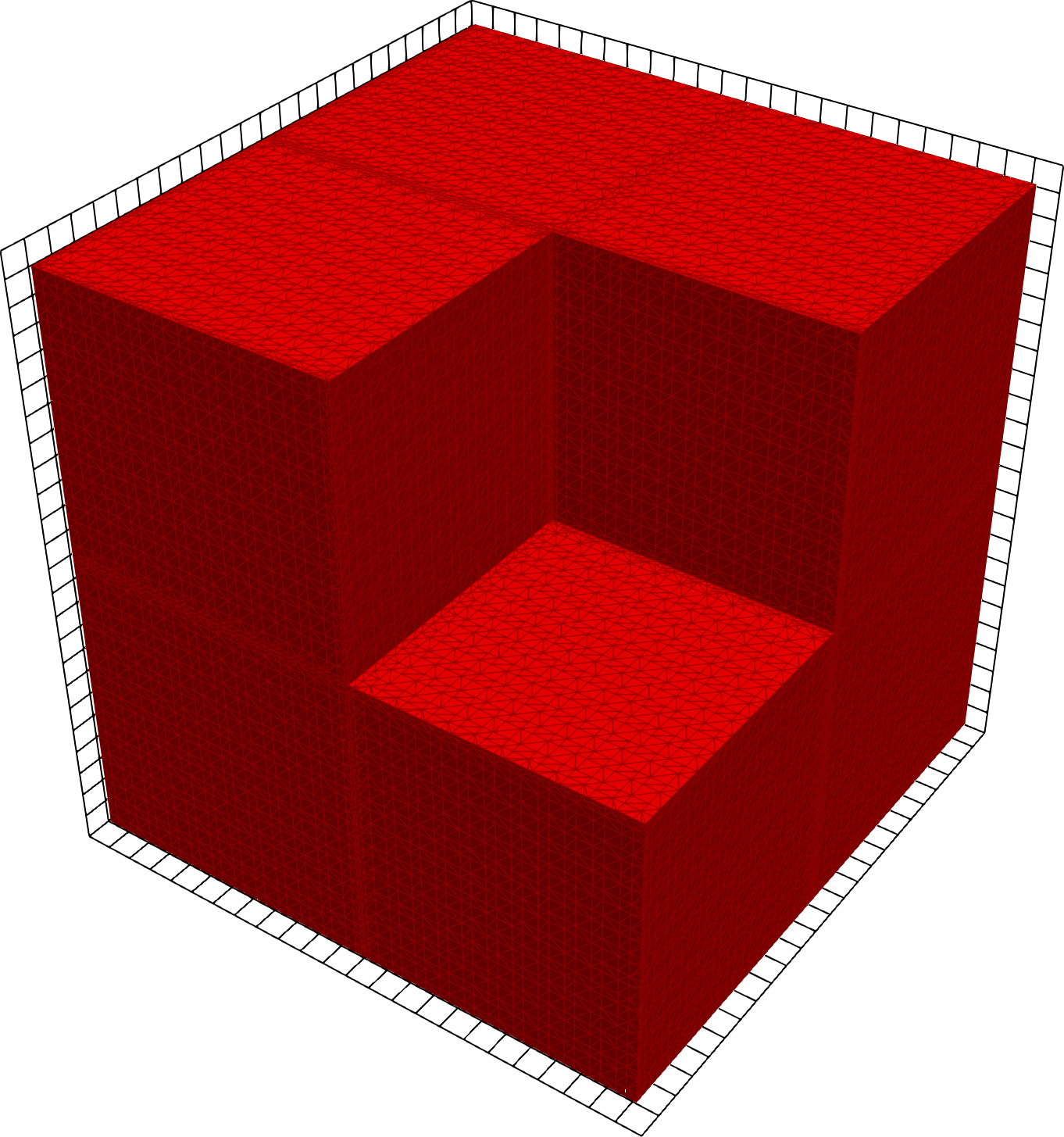}\hfill
  \includegraphics[width=0.19\textwidth]{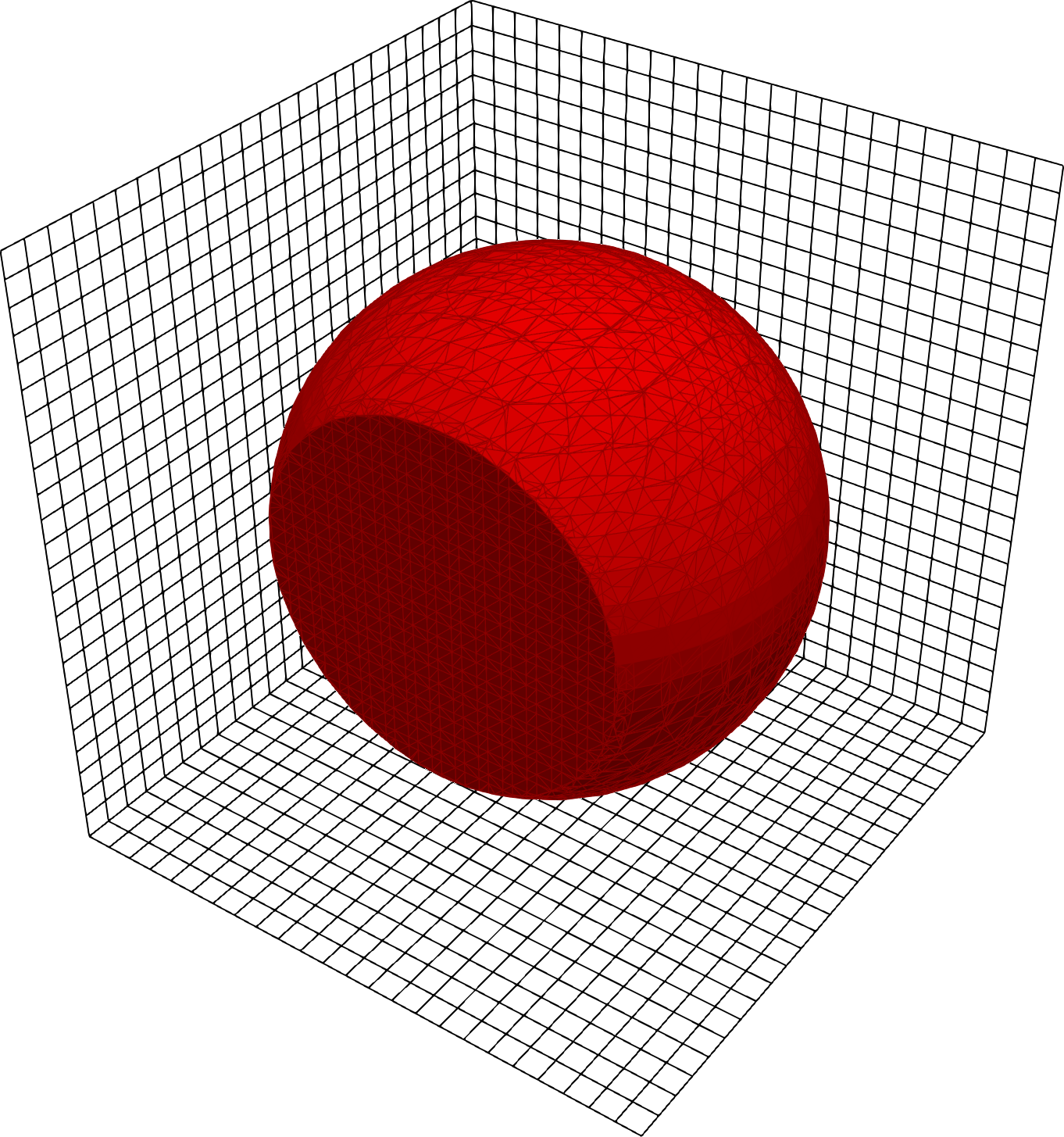}\hfill
  \includegraphics[width=0.19\textwidth]{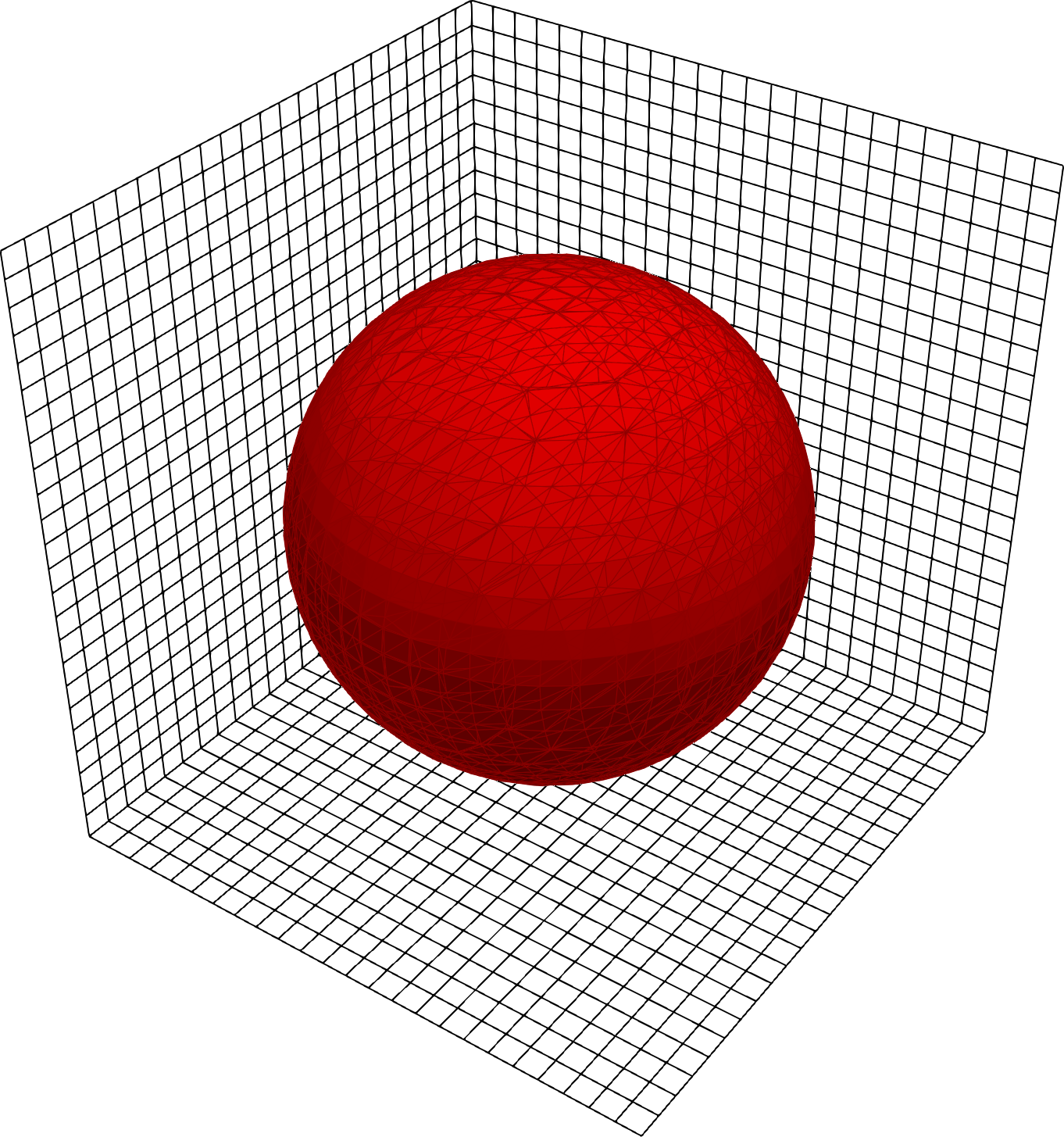}\hfill
  \includegraphics[width=0.19\textwidth]{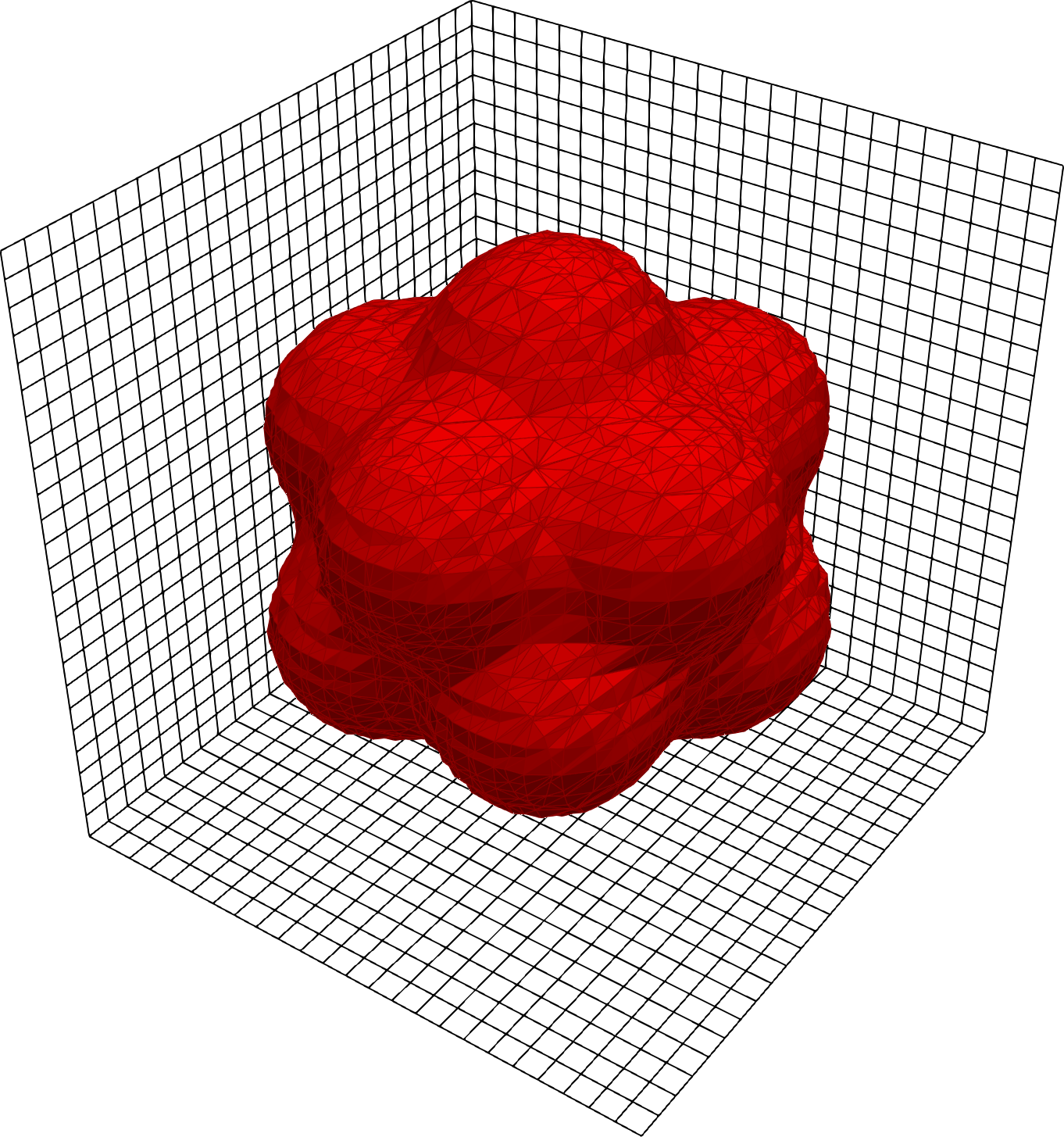}
  \caption{Benchmark geometries on a background mesh with $n = 24$. Top, from left to right:
    Square, L-shape, C-Disk, Disk and Flower. Bottom, from left to right: Cube, Tetris,
    C-Sphere, Sphere and Popcorn. The $\eps$-parametrised geometries are shown for
    $\eps = 0.4$.}
  \label{fig:meshes}
\end{figure}

We consider two families of geometries, depicted in \fig{fig:meshes}. The first family is parametrised by $\eps$: for every $n$, the domain is placed relative to the background mesh so that the cut ratio of every flat piece of its boundary is pinned at $\eps$, which allows us to test robustness with respect to small cuts. The domain therefore changes slightly with $n$. In 2D, this family comprises the Square, the L-shape and the C-Disk, a disk truncated by a vertical line; in 3D, their analogues: the Cube, the Tetris, a cube with an octant removed, and the C-Sphere, a sphere truncated by a plane. The second family consists of fixed shapes, independent of both $n$ and $\eps$, whose curved boundaries produce arbitrary cut configurations: the Disk and the five-petal Flower in 2D, and the Sphere and the Popcorn in 3D.
\tab{tab:meshes} reports the size of the resulting discretisations. The counts of the $\eps$-parametrised geometries do not depend on $\eps$, since varying $\eps$ only moves the boundary within the same layer of cut cells and hence leaves the active mesh unchanged.

\begin{table}[htbp]
  \centering
  \scriptsize
  \setlength{\tabcolsep}{1.5pt}
  \begin{tabular}{lr*{11}{r}}
\toprule
 & & \multicolumn{2}{c}{\# Cells} & \multicolumn{3}{c}{\# DOFs $H(\Grad)$} & \multicolumn{3}{c}{\# DOFs $\bm{H}(\Div)$} & \multicolumn{3}{c}{\# DOFs $\bm{H}(\Curl)$} \\
\cmidrule(lr){3-4} \cmidrule(lr){5-7} \cmidrule(lr){8-10} \cmidrule(lr){11-13}
 & $n$ & IN & CUT & $k = 1$ & $k = 2$ & $k = 3$ & $k = 0$ & $k = 1$ & $k = 2$ & $k = 0$ & $k = 1$ & $k = 2$ \\
\midrule
\multirow{6}{*}{Square} & 24 & 484 & 92 & 625 & 2\,401 & 5\,329 & 1\,200 & 4\,704 & 10\,512 & 1\,200 & 4\,704 & 10\,512 \\
 & 48 & 2\,116 & 188 & 2\,401 & 9\,409 & 21\,025 & 4\,704 & 18\,624 & 41\,760 & 4\,704 & 18\,624 & 41\,760 \\
 & 96 & 8\,836 & 380 & 9\,409 & 37\,249 & 83\,521 & 18\,624 & 74\,112 & 166\,464 & 18\,624 & 74\,112 & 166\,464 \\
 & 192 & 36\,100 & 764 & 37\,249 & 148\,225 & 332\,929 & 74\,112 & 295\,680 & 664\,704 & 74\,112 & 295\,680 & 664\,704 \\
 & 384 & 145\,924 & 1\,532 & 148\,225 & 591\,361 & 1\,329\,409 & 295\,680 & 1\,181\,184 & 2\,656\,512 & 295\,680 & 1\,181\,184 & 2\,656\,512 \\
 & 768 & 586\,756 & 3\,068 & 591\,361 & 2\,362\,369 & 5\,313\,025 & 1\,181\,184 & 4\,721\,664 & 10\,621\,440 & 1\,181\,184 & 4\,721\,664 & 10\,621\,440 \\
\midrule
\multirow{6}{*}{L-shape} & 24 & 363 & 92 & 504 & 1\,917 & 4\,240 & 958 & 3\,736 & 8\,334 & 958 & 3\,736 & 8\,334 \\
 & 48 & 1\,587 & 188 & 1\,872 & 7\,293 & 16\,264 & 3\,646 & 14\,392 & 32\,238 & 3\,646 & 14\,392 & 32\,238 \\
 & 96 & 6\,627 & 380 & 7\,200 & 28\,413 & 63\,640 & 14\,206 & 56\,440 & 126\,702 & 14\,206 & 56\,440 & 126\,702 \\
 & 192 & 27\,075 & 764 & 28\,224 & 112\,125 & 251\,704 & 56\,062 & 223\,480 & 502\,254 & 56\,062 & 223\,480 & 502\,254 \\
 & 384 & 109\,443 & 1\,532 & 111\,744 & 445\,437 & 1\,001\,080 & 222\,718 & 889\,336 & 1\,999\,854 & 222\,718 & 889\,336 & 1\,999\,854 \\
 & 768 & 440\,067 & 3\,068 & 444\,672 & 1\,775\,613 & 3\,992\,824 & 887\,806 & 3\,548\,152 & 7\,981\,038 & 887\,806 & 3\,548\,152 & 7\,981\,038 \\
\midrule
\multirow{6}{*}{C-Disk} & 24 & 228 & 70 & 336 & 1\,267 & 2\,794 & 633 & 2\,458 & 5\,475 & 633 & 2\,458 & 5\,475 \\
 & 48 & 948 & 142 & 1\,164 & 4\,507 & 10\,030 & 2\,253 & 8\,866 & 19\,839 & 2\,253 & 8\,866 & 19\,839 \\
 & 96 & 3\,924 & 280 & 4\,347 & 17\,101 & 38\,263 & 8\,550 & 33\,916 & 76\,098 & 8\,550 & 33\,916 & 76\,098 \\
 & 192 & 15\,894 & 556 & 16\,731 & 66\,361 & 148\,891 & 33\,180 & 132\,160 & 296\,940 & 33\,180 & 132\,160 & 296\,940 \\
 & 384 & 64\,076 & 1\,114 & 65\,750 & 261\,879 & 588\,388 & 130\,939 & 522\,638 & 1\,175\,097 & 130\,939 & 522\,638 & 1\,175\,097 \\
 & 768 & 257\,138 & 2\,230 & 260\,486 & 1\,039\,707 & 2\,337\,664 & 519\,853 & 2\,077\,178 & 4\,671\,975 & 519\,853 & 2\,077\,178 & 4\,671\,975 \\
\midrule
\multirow{6}{*}{Disk} & 24 & 224 & 76 & 341 & 1\,281 & 2\,821 & 640 & 2\,480 & 5\,520 & 640 & 2\,480 & 5\,520 \\
 & 48 & 968 & 148 & 1\,193 & 4\,617 & 10\,273 & 2\,308 & 9\,080 & 20\,316 & 2\,308 & 9\,080 & 20\,316 \\
 & 96 & 4\,044 & 292 & 4\,485 & 17\,641 & 39\,469 & 8\,820 & 34\,984 & 78\,492 & 8\,820 & 34\,984 & 78\,492 \\
 & 192 & 16\,420 & 580 & 17\,293 & 68\,585 & 153\,877 & 34\,292 & 136\,584 & 306\,876 & 34\,292 & 136\,584 & 306\,876 \\
 & 384 & 66\,324 & 1\,164 & 68\,073 & 271\,121 & 609\,145 & 135\,560 & 541\,072 & 1\,216\,536 & 135\,560 & 541\,072 & 1\,216\,536 \\
 & 768 & 266\,400 & 2\,332 & 269\,901 & 1\,077\,265 & 2\,422\,093 & 538\,632 & 2\,152\,192 & 4\,840\,680 & 538\,632 & 2\,152\,192 & 4\,840\,680 \\
\midrule
\multirow{6}{*}{Flower} & 24 & 134 & 72 & 245 & 901 & 1\,969 & 450 & 1\,724 & 3\,822 & 450 & 1\,724 & 3\,822 \\
 & 48 & 598 & 156 & 834 & 3\,175 & 7\,024 & 1\,587 & 6\,190 & 13\,809 & 1\,587 & 6\,190 & 13\,809 \\
 & 96 & 2\,528 & 306 & 2\,990 & 11\,647 & 25\,972 & 5\,823 & 22\,982 & 51\,477 & 5\,823 & 22\,982 & 51\,477 \\
 & 192 & 10\,432 & 610 & 11\,350 & 44\,783 & 100\,300 & 22\,391 & 88\,950 & 199\,677 & 22\,391 & 88\,950 & 199\,677 \\
 & 384 & 42\,398 & 1\,220 & 44\,231 & 175\,697 & 394\,399 & 87\,848 & 350\,168 & 786\,960 & 87\,848 & 350\,168 & 786\,960 \\
 & 768 & 170\,774 & 2\,432 & 174\,425 & 695\,261 & 1\,562\,509 & 347\,630 & 1\,388\,084 & 3\,121\,362 & 347\,630 & 1\,388\,084 & 3\,121\,362 \\
\midrule[\heavyrulewidth]
\midrule
\multirow{4}{*}{Cube} & 12 & 1\,000 & 728 & 2\,197 & 15\,625 & 50\,653 & 5\,616 & 43\,200 & 143\,856 & 6\,084 & 45\,000 & 147\,852 \\
 & 24 & 10\,648 & 3\,176 & 15\,625 & 117\,649 & 389\,017 & 43\,200 & 338\,688 & 1\,135\,296 & 45\,000 & 345\,744 & 1\,151\,064 \\
 & 48 & 97\,336 & 13\,256 & 117\,649 & 912\,673 & 3\,048\,625 & 338\,688 & 2\,681\,856 & 9\,020\,160 & 345\,744 & 2\,709\,792 & 9\,082\,800 \\
 & 96 & 830\,584 & 54\,152 & 912\,673 & 7\,189\,057 & 24\,137\,569 & 2\,681\,856 & 21\,344\,256 & 71\,912\,448 & 2\,709\,792 & 21\,455\,424 & 72\,162\,144 \\
\midrule
\multirow{4}{*}{Tetris} & 12 & 875 & 728 & 2\,072 & 14\,625 & 47\,278 & 5\,241 & 40\,200 & 133\,731 & 5\,709 & 42\,000 & 137\,727 \\
 & 24 & 9\,317 & 3\,176 & 14\,294 & 107\,001 & 353\,080 & 39\,207 & 306\,744 & 1\,027\,485 & 41\,007 & 313\,800 & 1\,043\,253 \\
 & 48 & 85\,169 & 13\,256 & 105\,482 & 815\,337 & 2\,720\,116 & 302\,187 & 2\,389\,848 & 8\,034\,633 & 309\,243 & 2\,417\,784 & 8\,097\,273 \\
 & 96 & 726\,761 & 54\,152 & 808\,850 & 6\,358\,473 & 21\,334\,348 & 2\,370\,387 & 18\,852\,504 & 63\,502\,785 & 2\,398\,323 & 18\,963\,672 & 63\,752\,481 \\
\midrule
\multirow{4}{*}{C-Sphere} & 12 & 268 & 396 & 946 & 6\,379 & 20\,284 & 2\,244 & 16\,944 & 56\,052 & 2\,525 & 18\,010 & 58\,407 \\
 & 24 & 2\,716 & 1\,636 & 5\,338 & 38\,643 & 126\,028 & 13\,984 & 108\,160 & 360\,864 & 14\,969 & 111\,986 & 369\,387 \\
 & 48 & 23\,992 & 6\,380 & 33\,898 & 256\,851 & 851\,092 & 94\,528 & 742\,576 & 2\,490\,840 & 98\,053 & 756\,450 & 2\,521\,887 \\
 & 96 & 204\,092 & 25\,528 & 243\,041 & 1\,890\,201 & 6\,319\,201 & 702\,060 & 5\,563\,680 & 18\,718\,020 & 715\,480 & 5\,616\,920 & 18\,837\,480 \\
\midrule
\multirow{4}{*}{Sphere} & 12 & 208 & 416 & 919 & 6\,109 & 19\,315 & 2\,136 & 16\,032 & 52\,920 & 2\,430 & 17\,148 & 55\,386 \\
 & 24 & 2\,464 & 1\,568 & 4\,993 & 35\,977 & 117\,145 & 12\,996 & 100\,368 & 334\,692 & 13\,956 & 104\,088 & 342\,972 \\
 & 48 & 22\,400 & 6\,248 & 32\,111 & 242\,805 & 803\,971 & 89\,292 & 700\,944 & 2\,350\,620 & 92\,754 & 714\,564 & 2\,381\,094 \\
 & 96 & 191\,208 & 25\,136 & 229\,575 & 1\,783\,229 & 5\,959\,027 & 662\,040 & 5\,244\,288 & 17\,640\,936 & 675\,270 & 5\,296\,764 & 17\,758\,674 \\
\midrule
\multirow{4}{*}{Popcorn} & 12 & 216 & 448 & 975 & 6\,489 & 20\,527 & 2\,270 & 17\,048 & 56\,286 & 2\,580 & 18\,224 & 58\,884 \\
 & 24 & 2\,524 & 1\,796 & 5\,393 & 38\,725 & 125\,917 & 13\,970 & 107\,720 & 359\,010 & 15\,042 & 111\,884 & 368\,286 \\
 & 48 & 23\,344 & 7\,216 & 34\,551 & 260\,177 & 860\,239 & 95\,538 & 748\,872 & 2\,510\,082 & 99\,528 & 764\,568 & 2\,545\,200 \\
 & 96 & 200\,304 & 28\,772 & 244\,217 & 1\,892\,657 & 6\,319\,777 & 702\,112 & 5\,557\,360 & 18\,689\,112 & 717\,252 & 5\,617\,408 & 18\,823\,836 \\
\bottomrule
\end{tabular}

  \caption{Number of interior (IN) and cut (CUT) cells, and number of \acp{dof} of the $H(\Grad)$-, $\bm{H}(\Div)$- and $\bm{H}(\Curl)$-conforming spaces of polynomial degree $k$ on the active mesh, for each geometry and background mesh of $n^d$ cells.}
  \label{tab:meshes}
\end{table}

In all experiments, we solve the linear systems of the stabilised problems $\Lgradh$, $\Ldivh$ and $\Lcurlh$ of \sect{sec:formulations}, with manufactured solutions. Unless stated otherwise, the stabilisation weight is $\tau^2 = 10$. The $H(\Curl)$ problems are only considered in three dimensions, since in two dimensions the N\'ed\'elec space is the rotation of the Raviart--Thomas space and the $H(\Curl)$ Riesz map is the rotated $H(\Div)$ one, so that the corresponding iteration counts, also for the $H(\Curl)$ block of the Hodge Laplacian of \sect{subsec:results-multiphysics}, coincide with those of $H(\Div)$.
The right-hand sides are computed from the exact solutions $u = \prod_{i=1}^d \sin(2\pi x_i)$ for $\Lgradh$, and $\bm{u} = \bm{s} + \Grad \varphi$ for $\Ldivh$ and $\Lcurlh$, with $\varphi = \pi^{-1} \prod_{i=1}^d \sin(\pi x_i)$ and the divergence-free field $\bm{s} = (\sin(\pi y),\, \sin(\pi x))$ in two dimensions and $\bm{s} = (\sin(\pi y),\, \sin(\pi z),\, \sin(\pi x))$ in three dimensions.

Every preconditioner is used within an FGMRES solver, stopped once the residual has been reduced by a factor $10^{-10}$. All smoothers perform $n_\ell = 10$ damped Richardson iterations \eqref{eq:effective-smoother}. The aggregate-based star smoothers $D_\ell^{\rm agg, all}$ and $D_\ell^{\rm agg}$ use the damping $\theta_\ell = 0.1$ in 2D and $\theta_\ell = 0.05$ in 3D, whereas the smoother $ \Dlact$ on the active problem, whose patches are regular vertex stars, uses $\theta_\ell = 0.2$ in both dimensions. The Nitsche penalty is $\gamma = 10$ for the lowest order of each space and $\gamma = 40$ for the higher orders, and every cut cell is aggregated. The coarsest-level problems, as well as the active problem in \sect{subsec:results-auxspace}, are solved by a sparse direct solver.

\subsection{Convergence of the discretisations}\label{subsec:results-convergence}

\begin{figure}[htbp]
  \centering
  \includegraphics[width=\textwidth]{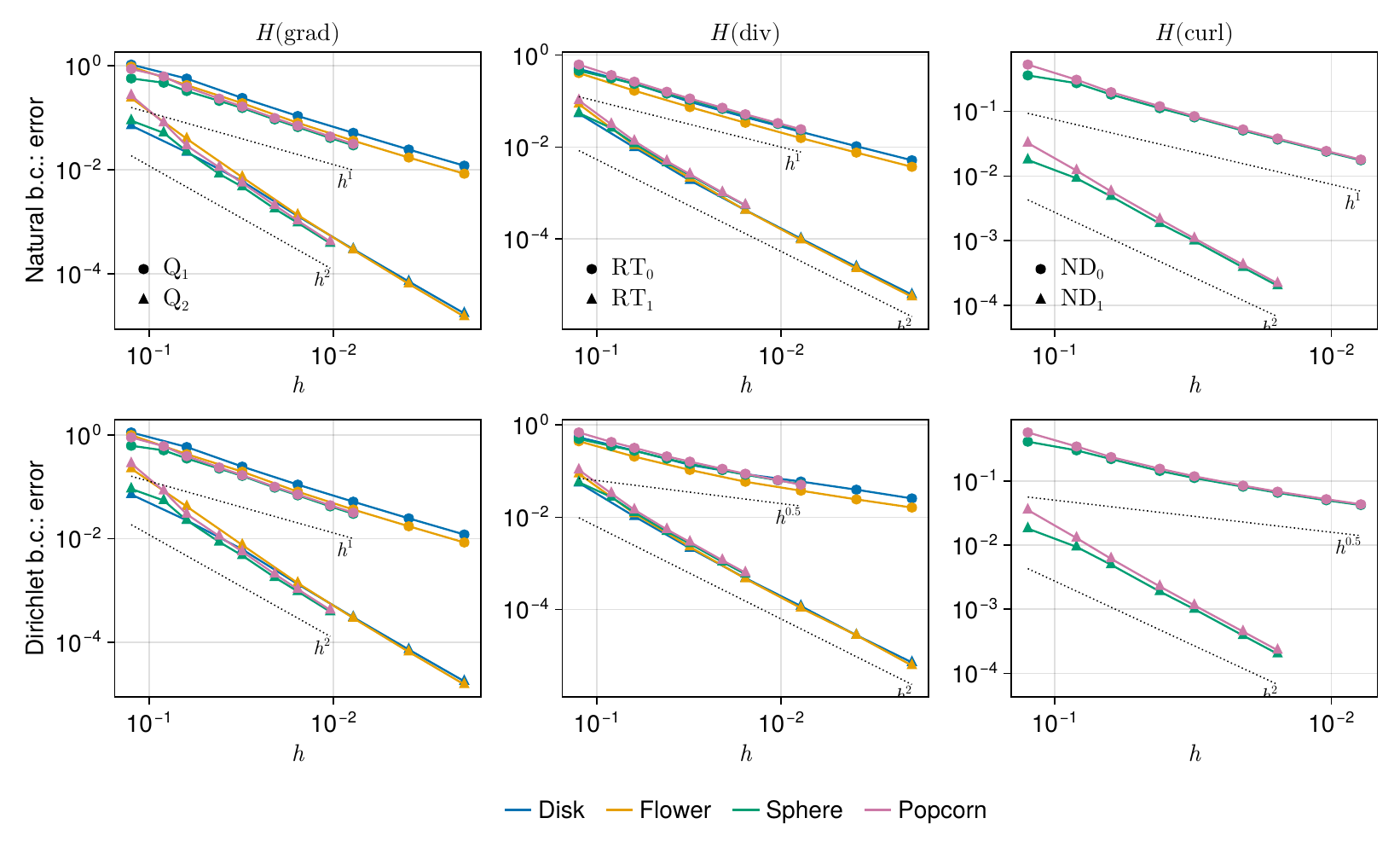}
  \caption{Convergence of the discretisations (\sect{subsec:results-convergence}). Discretisation error of $\Lgradh$, $\Ldivh$ and $\Lcurlh$ on the fixed geometries for the two lowest orders of each space, with natural (top) and Dirichlet (bottom) boundary conditions. Dotted lines indicate the reference rates.}
  \label{fig:convergence}
\end{figure}

Before turning to the preconditioners, we verify the convergence of the stabilised discretisations on the fixed geometries, where refining the mesh resolves the same domain. \fig{fig:convergence} reports the discretisation error of $\Lgradh$, $\Ldivh$ and $\Lcurlh$ for the two lowest orders of each space, solved by a direct solver, with natural boundary conditions and with Dirichlet boundary conditions on the whole of $\bnd$ imposed through \eqref{eq:nitsche-weakforms}. The error is measured in the norm $\norm{\rm x}{v}^2 \coloneqq \norm{L^2(\dom)}{v}^2 + \norm{L^2(\dom)}{\mathrm{D}^{\rm x} v}^2$, to which the term $h^{-1} \norm{L^2(\bnd)}{{\rm tr}^{\rm x}\, v}^2$ is added in the Dirichlet case, where ${\rm tr}^{\rm x} v \in \{v,\, \bm{v}\times\normal,\, \bm{v}\cdot\normal\}$ is the trace penalised by $b^{\rm x}$. All discretisations converge at the optimal rates, namely $k$ for $\Qk{k}$ and $k+1$ for $\mathbb{ND}_k$ and $\RT{k}$. With Dirichlet boundary conditions, the boundary term limits the rate guaranteed by Nitsche\'s method to $k+1/2$ for the vector-valued spaces. This rate is observed for $\mathbb{ND}_0$ and $\RT{0}$, whereas the higher-order discretisations converge at $k+1$, beyond the theoretical bound.

\subsection{Auxiliary space preconditioners}\label{subsec:results-auxspace}

We first compare the auxiliary-space preconditioners $\Bx{add}$ and $\Bx{mult}$ of \sect{sec:auxspace}, with smoother $\Dh = D_L^{\rm agg}$, against the active preconditioner $\Bx{ac}$. All three apply $(\Aact)^{-1}$ exactly, by means of a sparse direct solver, so that these experiments isolate the effect of the smoother.

We report iteration numbers of the FGMRES solver for each preconditioner in Tables~\ref{tab:auxspace-eps} and \ref{tab:auxspace-fixed}, for the lowest-order discretisations and all geometries. On the one hand, \tab{tab:auxspace-eps} contains iteration counts for the $\eps$-parametrised geometries. Moving along a row refines the mesh at a fixed cut ratio $\eps$, while moving down the rows of a preconditioner reduces $\eps$ on a fixed mesh. On the other hand, \tab{tab:auxspace-fixed} contains iteration counts for the fixed geometries. In this case, the geometries are fixed and each row refines the mesh, so that the cut configurations change arbitrarily from one mesh to the next.

We observe that $\Bx{ac}$ is robust with respect to both the cut ratio $\eps$ and the mesh size. This is expected from the stabilisation assumption \eqref{eq:stab-property}. However, as discussed in \sect{sec:auxspace}, the constants in \eqref{eq:stab-property} can be very large depending on the regularity of the domain and the discretisation itself. We can see the iteration counts reflect this behaviour, being very high on some geometries (over 200 iterations on the Tetris) and for $\Hcurl{\dom}$.
As expected, the auxiliary-space preconditioners $\Bx{add}$ and $\Bx{mult}$ solve this issue, and generally perform much better across all tested geometries and discretisations, while maintaining robustness with respect to both the cut ratio $\eps$ and the mesh size.

\begin{table}[htbp]
  \centering
  \scriptsize
  \setlength{\tabcolsep}{1.5pt}
  \begin{tabular}{lll|w{c}{1.5em}w{c}{1.5em}w{c}{1.5em}w{c}{1.5em}w{c}{1.5em}|w{c}{1.5em}w{c}{1.5em}w{c}{1.5em}w{c}{1.5em}w{c}{1.5em}|w{c}{1.5em}w{c}{1.5em}w{c}{1.5em}w{c}{1.5em}w{c}{1.5em}|w{c}{1.5em}w{c}{1.5em}w{c}{1.5em}|w{c}{1.5em}w{c}{1.5em}w{c}{1.5em}|w{c}{1.5em}w{c}{1.5em}w{c}{1.5em}}
\toprule
 & & & \multicolumn{5}{c|}{Square} & \multicolumn{5}{c|}{L-shape} & \multicolumn{5}{c|}{C-Disk} & \multicolumn{3}{c|}{Cube} & \multicolumn{3}{c|}{Tetris} & \multicolumn{3}{c}{C-Sphere} \\
\cmidrule(lr){4-8} \cmidrule(lr){9-13} \cmidrule(lr){14-18} \cmidrule(lr){19-21} \cmidrule(lr){22-24} \cmidrule(lr){25-27}
 & & $\eps$ \textbackslash{} $n$ & 24 & 48 & 96 & 192 & 384 & 24 & 48 & 96 & 192 & 384 & 24 & 48 & 96 & 192 & 384 & 12 & 24 & 48 & 12 & 24 & 48 & 12 & 24 & 48 \\
\midrule\midrule
\multirow{12}{*}{$\Lgrad$} & \multirow{4}{*}{$B_{\rm ac}^{\rm g}$} & 0.5 & 26 & 31 & 30 & 29 & 27 & 36 & 34 & 34 & 32 & 31 & 55 & 65 & 70 & 72 & 78 & 42 & 47 & 44 & 78 & 71 & 66 & 76 & 86 & 86 \\
 &  & 0.1 & 39 & 46 & 47 & 46 & 44 & 55 & 54 & 52 & 49 & 48 & 62 & 74 & 76 & 79 & 81 & 65 & 102 & 98 & 160 & 160 & 147 & 86 & 95 & 93 \\
 &  & 0.01 & 40 & 49 & 52 & 51 & 49 & 62 & 60 & 57 & 55 & 53 & 69 & 81 & 83 & 85 & 88 & 69 & 119 & 120 & 191 & 198 & 180 & 86 & 105 & 102 \\
 &  & 0.001 & 39 & 49 & 53 & 51 & 50 & 62 & 60 & 58 & 55 & 53 & 69 & 81 & 83 & 86 & 89 & 66 & 119 & 123 & 193 & 201 & 184 & 86 & 106 & 103 \\
\cmidrule(lr){2-27}
 & \multirow{4}{*}{$B_{\rm add}^{\rm g}$} & 0.5 & 24 & 24 & 24 & 23 & 22 & 28 & 28 & 26 & 26 & 25 & 42 & 44 & 47 & 48 & 50 & 31 & 32 & 31 & 48 & 47 & 44 & 63 & 73 & 72 \\
 &  & 0.1 & 26 & 28 & 27 & 26 & 25 & 32 & 32 & 30 & 30 & 28 & 41 & 45 & 47 & 48 & 50 & 38 & 41 & 39 & 76 & 68 & 64 & 64 & 74 & 73 \\
 &  & 0.01 & 26 & 29 & 28 & 26 & 26 & 33 & 32 & 32 & 30 & 29 & 41 & 45 & 47 & 48 & 50 & 39 & 43 & 41 & 83 & 74 & 70 & 66 & 77 & 76 \\
 &  & 0.001 & 26 & 29 & 28 & 26 & 26 & 34 & 32 & 32 & 30 & 29 & 41 & 45 & 47 & 48 & 50 & 39 & 43 & 41 & 84 & 75 & 71 & 67 & 77 & 77 \\
\cmidrule(lr){2-27}
 & \multirow{4}{*}{$B_{\rm mult}^{\rm g}$} & 0.5 & 10 & 10 & 9 & 9 & 9 & 10 & 10 & 10 & 9 & 9 & 13 & 13 & 13 & 14 & 14 & 23 & 24 & 23 & 28 & 27 & 25 & 34 & 39 & 38 \\
 &  & 0.1 & 10 & 10 & 10 & 9 & 9 & 11 & 11 & 11 & 10 & 10 & 13 & 13 & 13 & 14 & 14 & 26 & 28 & 27 & 43 & 37 & 35 & 35 & 39 & 39 \\
 &  & 0.01 & 10 & 10 & 10 & 9 & 9 & 12 & 11 & 11 & 10 & 10 & 13 & 13 & 13 & 14 & 14 & 26 & 28 & 28 & 47 & 41 & 38 & 35 & 41 & 41 \\
 &  & 0.001 & 10 & 10 & 10 & 9 & 9 & 12 & 11 & 11 & 10 & 10 & 13 & 13 & 13 & 14 & 14 & 25 & 28 & 28 & 47 & 41 & 39 & 36 & 41 & 41 \\
\midrule\midrule
\multirow{12}{*}{$\Ldiv$} & \multirow{4}{*}{$B_{\rm ac}^{\rm d}$} & 0.5 & 29 & 27 & 25 & 23 & 22 & 31 & 30 & 29 & 27 & 25 & 55 & 62 & 62 & 62 & 63 & 38 & 37 & 35 & 53 & 50 & 47 & 87 & 83 & 82 \\
 &  & 0.1 & 46 & 44 & 40 & 37 & 35 & 48 & 46 & 44 & 41 & 38 & 64 & 64 & 65 & 62 & 63 & 72 & 68 & 65 & 105 & 95 & 89 & 96 & 93 & 90 \\
 &  & 0.01 & 50 & 48 & 45 & 41 & 37 & 53 & 51 & 48 & 44 & 42 & 68 & 68 & 68 & 66 & 65 & 85 & 80 & 76 & 119 & 109 & 104 & 99 & 98 & 94 \\
 &  & 0.001 & 51 & 48 & 46 & 41 & 39 & 53 & 51 & 49 & 45 & 42 & 68 & 68 & 68 & 67 & 66 & 87 & 80 & 77 & 120 & 111 & 105 & 100 & 98 & 95 \\
\cmidrule(lr){2-27}
 & \multirow{4}{*}{$B_{\rm add}^{\rm d}$} & 0.5 & 25 & 23 & 23 & 21 & 20 & 27 & 27 & 27 & 25 & 24 & 40 & 41 & 40 & 40 & 39 & 26 & 26 & 25 & 36 & 35 & 34 & 48 & 49 & 49 \\
 &  & 0.1 & 27 & 27 & 24 & 23 & 21 & 33 & 31 & 30 & 28 & 28 & 40 & 41 & 40 & 40 & 39 & 30 & 29 & 28 & 43 & 42 & 40 & 48 & 49 & 49 \\
 &  & 0.01 & 28 & 27 & 25 & 24 & 21 & 33 & 32 & 31 & 29 & 28 & 39 & 41 & 40 & 40 & 39 & 31 & 30 & 28 & 46 & 45 & 43 & 48 & 49 & 49 \\
 &  & 0.001 & 28 & 27 & 25 & 24 & 21 & 33 & 32 & 30 & 29 & 28 & 39 & 41 & 40 & 40 & 39 & 31 & 30 & 29 & 46 & 45 & 43 & 48 & 49 & 49 \\
\cmidrule(lr){2-27}
 & \multirow{4}{*}{$B_{\rm mult}^{\rm d}$} & 0.5 & 9 & 8 & 8 & 8 & 7 & 9 & 9 & 8 & 8 & 7 & 11 & 11 & 10 & 10 & 10 & 13 & 12 & 12 & 14 & 13 & 13 & 14 & 16 & 16 \\
 &  & 0.1 & 9 & 9 & 9 & 8 & 7 & 10 & 10 & 9 & 9 & 8 & 11 & 11 & 10 & 10 & 10 & 13 & 13 & 12 & 16 & 16 & 15 & 14 & 16 & 16 \\
 &  & 0.01 & 9 & 9 & 9 & 8 & 7 & 11 & 10 & 9 & 9 & 8 & 11 & 11 & 11 & 10 & 10 & 13 & 13 & 12 & 17 & 17 & 15 & 14 & 16 & 15 \\
 &  & 0.001 & 10 & 9 & 9 & 8 & 8 & 11 & 10 & 9 & 9 & 8 & 11 & 11 & 11 & 10 & 10 & 13 & 13 & 12 & 17 & 17 & 15 & 14 & 16 & 15 \\
\midrule\midrule
\multirow{12}{*}{$\Lcurl$} & \multirow{4}{*}{$B_{\rm ac}^{\rm c}$} & 0.5 &  &  &  &  &  &  &  &  &  &  &  &  &  &  &  & 54 & 50 & 49 & 79 & 75 & 71 & 100 & 96 & 96 \\
 &  & 0.1 &  &  &  &  &  &  &  &  &  &  &  &  &  &  &  & 121 & 113 & 105 & 186 & 168 & 157 & 112 & 108 & 107 \\
 &  & 0.01 &  &  &  &  &  &  &  &  &  &  &  &  &  &  &  & 147 & 138 & 126 & 223 & 207 & 193 & 117 & 116 & 116 \\
 &  & 0.001 &  &  &  &  &  &  &  &  &  &  &  &  &  &  &  & 151 & 140 & 129 & 227 & 212 & 196 & 117 & 117 & 116 \\
\cmidrule(lr){2-27}
 & \multirow{4}{*}{$B_{\rm add}^{\rm c}$} & 0.5 &  &  &  &  &  &  &  &  &  &  &  &  &  &  &  & 34 & 32 & 31 & 49 & 45 & 43 & 66 & 67 & 67 \\
 &  & 0.1 &  &  &  &  &  &  &  &  &  &  &  &  &  &  &  & 44 & 40 & 37 & 76 & 66 & 62 & 69 & 70 & 68 \\
 &  & 0.01 &  &  &  &  &  &  &  &  &  &  &  &  &  &  &  & 45 & 42 & 39 & 82 & 72 & 68 & 70 & 71 & 69 \\
 &  & 0.001 &  &  &  &  &  &  &  &  &  &  &  &  &  &  &  & 46 & 42 & 39 & 83 & 73 & 69 & 70 & 71 & 69 \\
\cmidrule(lr){2-27}
 & \multirow{4}{*}{$B_{\rm mult}^{\rm c}$} & 0.5 &  &  &  &  &  &  &  &  &  &  &  &  &  &  &  & 23 & 21 & 20 & 26 & 24 & 23 & 33 & 34 & 33 \\
 &  & 0.1 &  &  &  &  &  &  &  &  &  &  &  &  &  &  &  & 26 & 24 & 22 & 39 & 34 & 31 & 33 & 35 & 34 \\
 &  & 0.01 &  &  &  &  &  &  &  &  &  &  &  &  &  &  &  & 26 & 24 & 23 & 43 & 36 & 34 & 34 & 35 & 34 \\
 &  & 0.001 &  &  &  &  &  &  &  &  &  &  &  &  &  &  &  & 26 & 24 & 23 & 43 & 36 & 34 & 34 & 35 & 35 \\
\bottomrule
\end{tabular}

  \caption{Auxiliary-space preconditioners (\sect{subsec:results-auxspace}). FGMRES iterations of $\Bx{ac}$, $\Bx{add}$ and $\Bx{mult}$ for the lowest-order Riesz maps $\Lxh$ on the $\eps$-parametrised geometries.}
  \label{tab:auxspace-eps}
\end{table}

\begin{table}[htbp]
  \centering
  \scriptsize
  \setlength{\tabcolsep}{1.5pt}
  \begin{tabular}{ll|w{c}{1.5em}w{c}{1.5em}w{c}{1.5em}w{c}{1.5em}w{c}{1.5em}|w{c}{1.5em}w{c}{1.5em}w{c}{1.5em}w{c}{1.5em}w{c}{1.5em}|w{c}{1.5em}w{c}{1.5em}w{c}{1.5em}|w{c}{1.5em}w{c}{1.5em}w{c}{1.5em}}
\toprule
 & & \multicolumn{5}{c|}{Disk} & \multicolumn{5}{c|}{Flower} & \multicolumn{3}{c|}{Sphere} & \multicolumn{3}{c}{Popcorn} \\
\cmidrule(lr){3-7} \cmidrule(lr){8-12} \cmidrule(lr){13-15} \cmidrule(lr){16-18}
 & $n$ & 24 & 48 & 96 & 192 & 384 & 24 & 48 & 96 & 192 & 384 & 12 & 24 & 48 & 12 & 24 & 48 \\
\midrule\midrule
\multirow{3}{*}{$\Lgrad$} & $B_{\rm ac}^{\rm g}$ & 44 & 57 & 67 & 70 & 76 & 71 & 74 & 77 & 79 & 83 & 53 & 83 & 87 & 97 & 120 & 117 \\
 & $B_{\rm add}^{\rm g}$ & 40 & 45 & 47 & 50 & 51 & 45 & 49 & 52 & 50 & 52 & 46 & 69 & 70 & 68 & 77 & 76 \\
 & $B_{\rm mult}^{\rm g}$ & 12 & 13 & 13 & 14 & 14 & 13 & 14 & 15 & 14 & 15 & 25 & 36 & 38 & 34 & 38 & 39 \\
\midrule\midrule
\multirow{3}{*}{$\Ldiv$} & $B_{\rm ac}^{\rm d}$ & 57 & 60 & 60 & 59 & 62 & 62 & 65 & 67 & 67 & 64 & 79 & 79 & 80 & 85 & 90 & 85 \\
 & $B_{\rm add}^{\rm d}$ & 42 & 40 & 40 & 41 & 40 & 39 & 42 & 44 & 41 & 40 & 48 & 51 & 49 & 44 & 48 & 49 \\
 & $B_{\rm mult}^{\rm d}$ & 11 & 10 & 10 & 11 & 10 & 11 & 11 & 11 & 11 & 11 & 13 & 14 & 15 & 12 & 15 & 15 \\
\midrule\midrule
\multirow{3}{*}{$\Lcurl$} & $B_{\rm ac}^{\rm c}$ &  &  &  &  &  &  &  &  &  &  & 99 & 96 & 97 & 101 & 123 & 113 \\
 & $B_{\rm add}^{\rm c}$ &  &  &  &  &  &  &  &  &  &  & 71 & 68 & 64 & 63 & 72 & 69 \\
 & $B_{\rm mult}^{\rm c}$ &  &  &  &  &  &  &  &  &  &  & 33 & 34 & 32 & 30 & 33 & 33 \\
\bottomrule
\end{tabular}

  \caption{Auxiliary-space preconditioners (\sect{subsec:results-auxspace}). FGMRES iterations of $\Bx{ac}$, $\Bx{add}$ and $\Bx{mult}$ for the lowest-order Riesz maps $\Lxh$ on the fixed geometries.}
  \label{tab:auxspace-fixed}
\end{table}

\subsection{Multigrid preconditioners}\label{subsec:results-gmg}

We now replace the exact solves by multigrid. We compare the homogeneous multigrid preconditioner $\Bgmg$ of \sect{sec:embedded-gmg}, with smoothers $D_\ell^{\rm agg, all}$, and the multilevel auxiliary-space preconditioner $\Bmg$, obtained from $\Bx{mult}$ by replacing $(\Aact)^{-1}$ with the active multigrid preconditioner $\Bagmg$ of \eqref{eq:active-gmg-precond}, with smoothers $ \Dlact$. We use both preconditioners within an FGMRES solver. Each preconditioner is applied as a single V-cycle, and the coarsest level is solved exactly. The hierarchies are built by uniform coarsening of the finest background mesh as described in \sect{sec:gmg}.

Results are reported in Tables \ref{tab:gmg-eps} and \ref{tab:gmg-fixed}, for $\eps$-parametrised and fixed geometries, respectively.
The two tables explore the hierarchy in different ways. In \tab{tab:gmg-eps}, the coarsest mesh is the same for every run ($24^2$ cells in 2D and $12^3$ cells in 3D). The finest mesh size $n$ therefore determines the number of levels $L$, such that a finer mesh corresponds to a deeper hierarchy with constant coarsest mesh. Each row of the table then tests robustness with respect to the mesh size and the number of levels simultaneously, for a given $\eps$.
In \tab{tab:gmg-fixed}, the mesh size and the number of levels are instead varied independently. Every finest mesh is solved with hierarchies that go from a two-level method ($L = 1$) to the deepest one, whose coarsest mesh is the one of \tab{tab:gmg-eps}. Along a row, the number of levels is fixed and the finest mesh is refined, so that the coarsest mesh is refined with it; this tests robustness with respect to the mesh size. Along a column, the finest mesh is fixed and levels are added by coarsening further; this tests robustness with respect to the number of levels. The last entry of each column is the setting of \tab{tab:gmg-eps}.

The main exception is the Tetris, where the iteration counts of $\Bgmg$ grow with the number of levels and as $\eps$ decreases. This is a hard geometry, with a wall of cut cells along the re-entrant edge on every level of the hierarchy. $\Bmg$ is unaffected. The growth is not intrinsic to the method but a consequence of the uniform choice of parameters, in particular the small damping $\theta_\ell = 0.05$. \tab{tab:gmg-tetris} repeats the Tetris runs of \tab{tab:gmg-eps} with $\tau^2 = 1$ and $\theta_\ell = 0.4, 0.16, 0.14$ for $\Lgrad$, $\Lcurl$ and $\Ldiv$, respectively. Iteration counts drop significantly, and the dependence on $\eps$ essentially disappears. Such tuning could go further, e.g.\ with level-dependent parameters. This illustrates how the results shown here could be significantly improved with careful parameter selection. We have instead opted for a uniform choice of parameters across all geometries and levels, prioritising comparability between runs.

The above discussion motivates a few remarks on the parameter tuning of these multigrid preconditioners. As a general rule, we have observed that lower values of the stabilisation weight $\tau^2$ lead to faster convergence as long as they do not compromise the well-posedness of the level problems.
Conversely, higher values of $\theta_\ell$ lead to faster convergence, but they must be chosen below the stability limit of the patch smoothers to avoid divergence. This limit is, as discussed in \sect{subsec:smoothers}, closely tied to the maximum overlap of the patches. For body-fitted problems, this overlap is moderate and predictable, allowing for relatively straightforward parameter selection. In contrast, our aggregation-based approach can lead to irregular patch overlaps, especially on challenging geometries like the Popcorn.
This is a practical advantage of $\Bmg$ over $\Bgmg$, since only its finest level carries the stabilisation terms and the aggregation-based irregular patch overlaps, so only its finest-level parameters need tuning. In contrast, $\Bgmg$ potentially requires tuning of parameters across the entire hierarchy.

\begin{table}[htbp]
  \centering
  \scriptsize
  \setlength{\tabcolsep}{1.5pt}
  \begin{tabular}{lll|w{c}{1.5em}w{c}{1.5em}w{c}{1.5em}w{c}{1.5em}w{c}{1.5em}|w{c}{1.5em}w{c}{1.5em}w{c}{1.5em}w{c}{1.5em}w{c}{1.5em}|w{c}{1.5em}w{c}{1.5em}w{c}{1.5em}w{c}{1.5em}w{c}{1.5em}|w{c}{1.5em}w{c}{1.5em}w{c}{1.5em}|w{c}{1.5em}w{c}{1.5em}w{c}{1.5em}|w{c}{1.5em}w{c}{1.5em}w{c}{1.5em}}
\toprule
 & & & \multicolumn{5}{c|}{Square} & \multicolumn{5}{c|}{L-shape} & \multicolumn{5}{c|}{C-Disk} & \multicolumn{3}{c|}{Cube} & \multicolumn{3}{c|}{Tetris} & \multicolumn{3}{c}{C-Sphere} \\
\cmidrule(lr){4-8} \cmidrule(lr){9-13} \cmidrule(lr){14-18} \cmidrule(lr){19-21} \cmidrule(lr){22-24} \cmidrule(lr){25-27}
 & & $\eps$ \textbackslash{} $n$ & 48 & 96 & 192 & 384 & 768 & 48 & 96 & 192 & 384 & 768 & 48 & 96 & 192 & 384 & 768 & 24 & 48 & 96 & 24 & 48 & 96 & 24 & 48 & 96 \\
\midrule\midrule
\multirow{8}{*}{$\Lgrad$} & \multirow{4}{*}{$B_{\rm hom}^{\rm g}$} & 0.5 & 10 & 10 & 11 & 11 & 11 & 14 & 18 & 21 & 23 & 25 & 18 & 20 & 22 & 23 & 24 & 15 & 15 & 15 & 27 & 35 & 43 & 45 & 45 & 49 \\
 &  & 0.1 & 9 & 10 & 10 & 11 & 11 & 16 & 20 & 22 & 25 & 27 & 18 & 20 & 22 & 23 & 24 & 16 & 16 & 15 & 42 & 57 & 69 & 45 & 45 & 49 \\
 &  & 0.01 & 9 & 10 & 10 & 11 & 11 & 17 & 20 & 24 & 26 & 28 & 18 & 20 & 22 & 23 & 24 & 16 & 16 & 15 & 55 & 75 & 91 & 45 & 46 & 49 \\
 &  & 0.001 & 9 & 10 & 10 & 11 & 11 & 17 & 20 & 24 & 26 & 28 & 18 & 20 & 22 & 23 & 24 & 16 & 16 & 15 & 59 & 81 & 98 & 45 & 46 & 49 \\
\cmidrule(lr){2-27}
 & \multirow{4}{*}{$B_{\rm aux}^{\rm g}$} & 0.5 & 10 & 10 & 10 & 9 & 9 & 11 & 12 & 12 & 12 & 12 & 13 & 14 & 15 & 16 & 16 & 23 & 23 & 22 & 26 & 25 & 25 & 38 & 38 & 41 \\
 &  & 0.1 & 10 & 10 & 10 & 10 & 10 & 12 & 13 & 13 & 13 & 12 & 14 & 14 & 15 & 16 & 16 & 27 & 26 & 25 & 36 & 34 & 33 & 39 & 39 & 41 \\
 &  & 0.01 & 10 & 10 & 10 & 10 & 10 & 12 & 13 & 13 & 13 & 13 & 14 & 14 & 15 & 16 & 16 & 28 & 27 & 25 & 39 & 37 & 35 & 40 & 41 & 42 \\
 &  & 0.001 & 10 & 10 & 10 & 10 & 10 & 13 & 13 & 13 & 13 & 13 & 14 & 14 & 15 & 16 & 16 & 28 & 27 & 26 & 40 & 38 & 36 & 40 & 41 & 43 \\
\midrule\midrule
\multirow{8}{*}{$\Ldiv$} & \multirow{4}{*}{$B_{\rm hom}^{\rm d}$} & 0.5 & 9 & 9 & 9 & 9 & 10 & 12 & 15 & 17 & 18 & 19 & 15 & 15 & 16 & 16 & 16 & 9 & 10 & 10 & 14 & 17 & 20 & 16 & 16 & 16 \\
 &  & 0.1 & 8 & 9 & 9 & 9 & 9 & 13 & 16 & 18 & 20 & 21 & 15 & 15 & 16 & 16 & 16 & 9 & 9 & 10 & 17 & 22 & 25 & 16 & 16 & 16 \\
 &  & 0.01 & 8 & 9 & 9 & 9 & 9 & 14 & 17 & 19 & 20 & 21 & 15 & 15 & 16 & 16 & 16 & 9 & 9 & 9 & 19 & 25 & 28 & 16 & 15 & 16 \\
 &  & 0.001 & 8 & 9 & 9 & 9 & 9 & 14 & 17 & 19 & 20 & 21 & 15 & 15 & 16 & 16 & 16 & 9 & 9 & 9 & 20 & 25 & 28 & 16 & 15 & 16 \\
\cmidrule(lr){2-27}
 & \multirow{4}{*}{$B_{\rm aux}^{\rm d}$} & 0.5 & 8 & 8 & 8 & 7 & 7 & 10 & 10 & 10 & 11 & 11 & 11 & 12 & 12 & 12 & 12 & 12 & 12 & 11 & 15 & 16 & 16 & 16 & 16 & 16 \\
 &  & 0.1 & 9 & 9 & 8 & 8 & 7 & 11 & 11 & 11 & 11 & 11 & 11 & 12 & 12 & 12 & 12 & 13 & 12 & 12 & 17 & 17 & 16 & 16 & 16 & 16 \\
 &  & 0.01 & 9 & 9 & 8 & 8 & 7 & 11 & 11 & 11 & 11 & 11 & 11 & 12 & 12 & 12 & 12 & 13 & 12 & 12 & 17 & 17 & 17 & 16 & 16 & 16 \\
 &  & 0.001 & 9 & 9 & 8 & 8 & 7 & 11 & 11 & 11 & 11 & 11 & 11 & 12 & 12 & 12 & 12 & 13 & 12 & 12 & 18 & 17 & 17 & 16 & 16 & 16 \\
\midrule\midrule
\multirow{8}{*}{$\Lcurl$} & \multirow{4}{*}{$B_{\rm hom}^{\rm c}$} & 0.5 &  &  &  &  &  &  &  &  &  &  &  &  &  &  &  & 13 & 13 & 14 & 24 & 31 & 39 & 38 & 38 & 38 \\
 &  & 0.1 &  &  &  &  &  &  &  &  &  &  &  &  &  &  &  & 13 & 13 & 13 & 38 & 52 & 62 & 38 & 38 & 37 \\
 &  & 0.01 &  &  &  &  &  &  &  &  &  &  &  &  &  &  &  & 13 & 13 & 13 & 51 & 69 & 82 & 38 & 39 & 37 \\
 &  & 0.001 &  &  &  &  &  &  &  &  &  &  &  &  &  &  &  & 13 & 13 & 13 & 54 & 75 & 87 & 38 & 39 & 37 \\
\cmidrule(lr){2-27}
 & \multirow{4}{*}{$B_{\rm aux}^{\rm c}$} & 0.5 &  &  &  &  &  &  &  &  &  &  &  &  &  &  &  & 21 & 19 & 19 & 25 & 25 & 25 & 34 & 33 & 34 \\
 &  & 0.1 &  &  &  &  &  &  &  &  &  &  &  &  &  &  &  & 23 & 22 & 19 & 32 & 31 & 31 & 34 & 34 & 34 \\
 &  & 0.01 &  &  &  &  &  &  &  &  &  &  &  &  &  &  &  & 23 & 22 & 20 & 35 & 33 & 32 & 35 & 34 & 34 \\
 &  & 0.001 &  &  &  &  &  &  &  &  &  &  &  &  &  &  &  & 24 & 22 & 20 & 35 & 33 & 33 & 35 & 34 & 35 \\
\bottomrule
\end{tabular}

  \caption{Multigrid preconditioners (\sect{subsec:results-gmg}). FGMRES iterations of $\Bgmg$ and $\Bmg$ for the lowest-order Riesz maps $\Lxh$ on the $\eps$-parametrised geometries, with a fixed coarsest mesh.}
  \label{tab:gmg-eps}
\end{table}

\begin{table}[htbp]
  \centering
  \scriptsize
  \setlength{\tabcolsep}{1.5pt}
  \begin{tabular}{llc|w{c}{1.5em}w{c}{1.5em}w{c}{1.5em}w{c}{1.5em}w{c}{1.5em}|w{c}{1.5em}w{c}{1.5em}w{c}{1.5em}w{c}{1.5em}w{c}{1.5em}|w{c}{1.5em}w{c}{1.5em}w{c}{1.5em}|w{c}{1.5em}w{c}{1.5em}w{c}{1.5em}}
\toprule
 & & & \multicolumn{5}{c|}{Disk} & \multicolumn{5}{c|}{Flower} & \multicolumn{3}{c|}{Sphere} & \multicolumn{3}{c}{Popcorn} \\
\cmidrule(lr){4-8} \cmidrule(lr){9-13} \cmidrule(lr){14-16} \cmidrule(lr){17-19}
 & & $L$ \textbackslash{} $n$ & 48 & 96 & 192 & 384 & 768 & 48 & 96 & 192 & 384 & 768 & 24 & 48 & 96 & 24 & 48 & 96 \\
\midrule\midrule
\multirow{10}{*}{$\Lgrad$} & \multirow{5}{*}{$B_{\rm hom}^{\rm g}$} & 1 & 17 & 18 & 21 & 22 & 22 & 18 & 21 & 22 & 22 & 23 & 39 & 45 & 52 & 44 & 55 & 54 \\
 &  & 2 &  & 20 & 22 & 23 & 23 &  & 22 & 23 & 23 & 23 &  & 42 & 48 &  & 52 & 51 \\
 &  & 3 &  &  & 22 & 23 & 24 &  &  & 23 & 24 & 24 &  &  & 49 &  &  & 54 \\
 &  & 4 &  &  &  & 24 & 25 &  &  &  & 24 & 24 &  &  &  &  &  &  \\
 &  & 5 &  &  &  &  & 26 &  &  &  &  & 24 &  &  &  &  &  &  \\
\cmidrule(lr){2-19}
 & \multirow{5}{*}{$B_{\rm aux}^{\rm g}$} & 1 & 13 & 14 & 15 & 15 & 15 & 15 & 15 & 15 & 16 & 15 & 36 & 37 & 40 & 37 & 38 & 41 \\
 &  & 2 &  & 14 & 15 & 15 & 16 &  & 16 & 15 & 16 & 16 &  & 37 & 40 &  & 39 & 42 \\
 &  & 3 &  &  & 15 & 15 & 16 &  &  & 16 & 16 & 16 &  &  & 40 &  &  & 42 \\
 &  & 4 &  &  &  & 15 & 16 &  &  &  & 16 & 16 &  &  &  &  &  &  \\
 &  & 5 &  &  &  &  & 16 &  &  &  &  & 17 &  &  &  &  &  &  \\
\midrule\midrule
\multirow{10}{*}{$\Ldiv$} & \multirow{5}{*}{$B_{\rm hom}^{\rm d}$} & 1 & 14 & 15 & 15 & 15 & 14 & 15 & 16 & 15 & 15 & 14 & 15 & 16 & 16 & 17 & 16 & 16 \\
 &  & 2 &  & 15 & 16 & 16 & 15 &  & 17 & 17 & 16 & 15 &  & 16 & 16 &  & 17 & 16 \\
 &  & 3 &  &  & 17 & 16 & 16 &  &  & 17 & 17 & 16 &  &  & 17 &  &  & 17 \\
 &  & 4 &  &  &  & 17 & 16 &  &  &  & 17 & 16 &  &  &  &  &  &  \\
 &  & 5 &  &  &  &  & 17 &  &  &  &  & 17 &  &  &  &  &  &  \\
\cmidrule(lr){2-19}
 & \multirow{5}{*}{$B_{\rm aux}^{\rm d}$} & 1 & 11 & 11 & 11 & 11 & 10 & 12 & 12 & 11 & 11 & 11 & 14 & 15 & 16 & 15 & 15 & 16 \\
 &  & 2 &  & 11 & 11 & 11 & 11 &  & 12 & 12 & 11 & 11 &  & 15 & 16 &  & 15 & 16 \\
 &  & 3 &  &  & 12 & 11 & 11 &  &  & 12 & 12 & 11 &  &  & 16 &  &  & 16 \\
 &  & 4 &  &  &  & 11 & 11 &  &  &  & 12 & 12 &  &  &  &  &  &  \\
 &  & 5 &  &  &  &  & 11 &  &  &  &  & 12 &  &  &  &  &  &  \\
\midrule\midrule
\multirow{6}{*}{$\Lcurl$} & \multirow{3}{*}{$B_{\rm hom}^{\rm c}$} & 1 &  &  &  &  &  &  &  &  &  &  & 34 & 37 & 40 & 37 & 43 & 39 \\
 &  & 2 &  &  &  &  &  &  &  &  &  &  &  & 35 & 37 &  & 42 & 37 \\
 &  & 3 &  &  &  &  &  &  &  &  &  &  &  &  & 38 &  &  & 40 \\
\cmidrule(lr){2-19}
 & \multirow{3}{*}{$B_{\rm aux}^{\rm c}$} & 1 &  &  &  &  &  &  &  &  &  &  & 33 & 32 & 33 & 33 & 32 & 34 \\
 &  & 2 &  &  &  &  &  &  &  &  &  &  &  & 32 & 34 &  & 33 & 34 \\
 &  & 3 &  &  &  &  &  &  &  &  &  &  &  &  & 34 &  &  & 35 \\
\bottomrule
\end{tabular}

  \caption{Multigrid preconditioners (\sect{subsec:results-gmg}). FGMRES iterations of $\Bgmg$ and $\Bmg$ for the lowest-order Riesz maps $\Lxh$ on the fixed geometries, with $L+1$ levels.}
  \label{tab:gmg-fixed}
\end{table}

\begin{table}[htbp]
  \centering
  \scriptsize
  \setlength{\tabcolsep}{1.5pt}
  \begin{tabular}{ll|w{c}{1.5em}w{c}{1.5em}w{c}{1.5em}|w{c}{1.5em}w{c}{1.5em}w{c}{1.5em}|w{c}{1.5em}w{c}{1.5em}w{c}{1.5em}|w{c}{1.5em}w{c}{1.5em}w{c}{1.5em}}
\toprule
 & & \multicolumn{3}{c|}{$\eps = 0.5$} & \multicolumn{3}{c|}{$\eps = 0.1$} & \multicolumn{3}{c|}{$\eps = 0.01$} & \multicolumn{3}{c}{$\eps = 0.001$} \\
\cmidrule(lr){3-5} \cmidrule(lr){6-8} \cmidrule(lr){9-11} \cmidrule(lr){12-14}
 & $\theta_\ell$ \textbackslash{} $n$ & 24 & 48 & 96 & 24 & 48 & 96 & 24 & 48 & 96 & 24 & 48 & 96 \\
\midrule\midrule
$B_{\rm hom}^{\rm g}$ & 0.40 & 9 & 11 & 13 & 10 & 12 & 14 & 11 & 13 & 15 & 11 & 14 & 15 \\
$B_{\rm hom}^{\rm d}$ & 0.14 & 8 & 10 & 12 & 8 & 10 & 12 & 9 & 11 & 12 & 9 & 11 & 12 \\
$B_{\rm hom}^{\rm c}$ & 0.16 & 10 & 14 & 17 & 13 & 17 & 19 & 15 & 19 & 22 & 15 & 20 & 23 \\
\bottomrule
\end{tabular}

  \caption{Multigrid preconditioners (\sect{subsec:results-gmg}). FGMRES iterations of $\Bgmg$ for the lowest-order Riesz maps $\Lxh$ on the Tetris, in the setting of \tab{tab:gmg-eps}, with $\tau^2 = 1$ and the damping $\theta_\ell$ of the second column.}
  \label{tab:gmg-tetris}
\end{table}

\subsection{Robustness with respect to the parameter \texorpdfstring{$\kappa$}{kappa}} \label{subsec:results-kappa}

Another important property of the patch smoothers developed in \cite{Arnold2000} is their robustness with respect to the problem parameters, which we now verify for the unfitted multigrid preconditioners.
In the previous experiments, we kept $\kappa = 1$. We now explore how the performance of the multigrid preconditioners changes as $\kappa$ increases.

Tables \ref{tab:gmg-kappa-emb} and \ref{tab:gmg-kappa-aux} repeat the experiments of \tab{tab:gmg-eps} on the fixed geometries (up to $n = 384$ in two dimensions) for $\kappa \in \{1, 10, 100\}$, for $\Bgmg$ and $\Bmg$ respectively. As expected, we observe that both preconditioners remain robust with respect to the parameter $\kappa$. In particular, iteration counts for $\Lgradh$ do not depend on $\kappa$, and those for $\Ldivh$ and $\Lcurlh$ decrease as $\kappa$ grows, for both preconditioners.

\begin{table}[htbp]
  \centering
  \scriptsize
  \setlength{\tabcolsep}{1.5pt}
  \begin{tabular}{llc|w{c}{1.5em}w{c}{1.5em}w{c}{1.5em}w{c}{1.5em}|w{c}{1.5em}w{c}{1.5em}w{c}{1.5em}w{c}{1.5em}|w{c}{1.5em}w{c}{1.5em}w{c}{1.5em}|w{c}{1.5em}w{c}{1.5em}w{c}{1.5em}}
\toprule
 & & & \multicolumn{4}{c|}{Disk} & \multicolumn{4}{c|}{Flower} & \multicolumn{3}{c|}{Sphere} & \multicolumn{3}{c}{Popcorn} \\
\cmidrule(lr){4-7} \cmidrule(lr){8-11} \cmidrule(lr){12-14} \cmidrule(lr){15-17}
 & & $\kappa$ \textbackslash{} $n$ & 48 & 96 & 192 & 384 & 48 & 96 & 192 & 384 & 24 & 48 & 96 & 24 & 48 & 96 \\
\midrule\midrule
\multirow{3}{*}{$\Lgrad$} & \multirow{3}{*}{$B_{\rm hom}^{\rm g}$} & 1 & 17 & 20 & 22 & 24 & 18 & 22 & 23 & 24 & 39 & 42 & 49 & 44 & 52 & 54 \\
 &  & 10 & 17 & 20 & 22 & 24 & 18 & 22 & 23 & 24 & 39 & 42 & 49 & 44 & 52 & 54 \\
 &  & 100 & 17 & 20 & 22 & 24 & 18 & 22 & 23 & 24 & 39 & 42 & 49 & 44 & 52 & 54 \\
\midrule\midrule
\multirow{3}{*}{$\Ldiv$} & \multirow{3}{*}{$B_{\rm hom}^{\rm d}$} & 1 & 14 & 15 & 17 & 17 & 15 & 17 & 17 & 17 & 15 & 16 & 17 & 17 & 17 & 17 \\
 &  & 10 & 11 & 12 & 13 & 12 & 12 & 13 & 13 & 13 & 12 & 13 & 13 & 14 & 13 & 13 \\
 &  & 100 & 8 & 9 & 9 & 9 & 8 & 9 & 9 & 9 & 9 & 10 & 9 & 10 & 10 & 10 \\
\midrule\midrule
\multirow{3}{*}{$\Lcurl$} & \multirow{3}{*}{$B_{\rm hom}^{\rm c}$} & 1 &  &  &  &  &  &  &  &  & 34 & 35 & 38 & 37 & 42 & 40 \\
 &  & 10 &  &  &  &  &  &  &  &  & 28 & 29 & 31 & 32 & 34 & 32 \\
 &  & 100 &  &  &  &  &  &  &  &  & 25 & 25 & 26 & 28 & 29 & 26 \\
\bottomrule
\end{tabular}

  \caption{Parameter robustness (\sect{subsec:results-kappa}). FGMRES iterations of $\Bgmg$ for the lowest-order Riesz maps $\Lxh$ with weight $\kappa$ on the fixed geometries, with a fixed coarsest mesh.}
  \label{tab:gmg-kappa-emb}
\end{table}

\begin{table}[htbp]
  \centering
  \scriptsize
  \setlength{\tabcolsep}{1.5pt}
  \begin{tabular}{llc|w{c}{1.5em}w{c}{1.5em}w{c}{1.5em}w{c}{1.5em}|w{c}{1.5em}w{c}{1.5em}w{c}{1.5em}w{c}{1.5em}|w{c}{1.5em}w{c}{1.5em}w{c}{1.5em}|w{c}{1.5em}w{c}{1.5em}w{c}{1.5em}}
\toprule
 & & & \multicolumn{4}{c|}{Disk} & \multicolumn{4}{c|}{Flower} & \multicolumn{3}{c|}{Sphere} & \multicolumn{3}{c}{Popcorn} \\
\cmidrule(lr){4-7} \cmidrule(lr){8-11} \cmidrule(lr){12-14} \cmidrule(lr){15-17}
 & & $\kappa$ \textbackslash{} $n$ & 48 & 96 & 192 & 384 & 48 & 96 & 192 & 384 & 24 & 48 & 96 & 24 & 48 & 96 \\
\midrule\midrule
\multirow{3}{*}{$\Lgrad$} & \multirow{3}{*}{$B_{\rm aux}^{\rm g}$} & 1 & 13 & 14 & 15 & 15 & 15 & 16 & 16 & 16 & 36 & 37 & 40 & 37 & 39 & 42 \\
 &  & 10 & 13 & 14 & 15 & 15 & 15 & 16 & 16 & 16 & 36 & 37 & 40 & 37 & 39 & 42 \\
 &  & 100 & 13 & 14 & 15 & 15 & 15 & 16 & 16 & 16 & 36 & 37 & 40 & 37 & 39 & 42 \\
\midrule\midrule
\multirow{3}{*}{$\Ldiv$} & \multirow{3}{*}{$B_{\rm aux}^{\rm d}$} & 1 & 11 & 11 & 12 & 11 & 12 & 12 & 12 & 12 & 14 & 15 & 16 & 15 & 15 & 16 \\
 &  & 10 & 9 & 9 & 9 & 9 & 9 & 10 & 9 & 9 & 11 & 11 & 11 & 11 & 12 & 12 \\
 &  & 100 & 6 & 7 & 7 & 7 & 7 & 7 & 7 & 7 & 8 & 7 & 7 & 8 & 8 & 7 \\
\midrule\midrule
\multirow{3}{*}{$\Lcurl$} & \multirow{3}{*}{$B_{\rm aux}^{\rm c}$} & 1 &  &  &  &  &  &  &  &  & 33 & 32 & 34 & 33 & 33 & 35 \\
 &  & 10 &  &  &  &  &  &  &  &  & 28 & 28 & 29 & 28 & 29 & 30 \\
 &  & 100 &  &  &  &  &  &  &  &  & 25 & 26 & 28 & 26 & 26 & 28 \\
\bottomrule
\end{tabular}

  \caption{Parameter robustness (\sect{subsec:results-kappa}). FGMRES iterations of $\Bmg$ for the lowest-order Riesz maps $\Lxh$ with weight $\kappa$ on the fixed geometries, with a fixed coarsest mesh.}
  \label{tab:gmg-kappa-aux}
\end{table}

\subsection{Higher-order discretisations}\label{subsec:results-highorder}

In general, \ac{gmg} preconditioners are not robust with respect to the polynomial order. However, in the context of body-fitted discretisations, vertex-star smoothers are known to be robust with respect to the polynomial order \cite{Miraci2020,Miraci2021,Innerberger2024,ChaumontFreletVohralik2026}. We now assess numerically whether this robustness also holds for our unfitted \ac{gmg} preconditioners.

\tab{tab:gmg-highorder} repeats the experiments in \tab{tab:gmg-fixed} for higher polynomial orders: $\Qk{2}$, $\Qk{3}$, $\RT{1}$ and $\RT{2}$ on the two-dimensional fixed geometries, and $\RT{1}$ on the three-dimensional ones. We can see that the iteration counts of $\Bmg$ remain stable as the polynomial order increases, and those of $\Bgmg$ remain bounded. Both are also lower than the lowest-order counts of \tab{tab:gmg-fixed}, indicating that the unfitted \ac{gmg} preconditioners retain robustness with respect to the polynomial order.

\begin{table}[htbp]
  \centering
  \scriptsize
  \setlength{\tabcolsep}{1pt}
  \begin{tabular}{llc|w{c}{1.5em}w{c}{1.5em}w{c}{1.5em}w{c}{1.5em}w{c}{1.5em}|w{c}{1.5em}w{c}{1.5em}w{c}{1.5em}w{c}{1.5em}w{c}{1.5em}|w{c}{1.5em}w{c}{1.5em}w{c}{1.5em}w{c}{1.5em}w{c}{1.5em}|w{c}{1.5em}w{c}{1.5em}w{c}{1.5em}w{c}{1.5em}w{c}{1.5em}|w{c}{1.5em}w{c}{1.5em}w{c}{1.5em}|w{c}{1.5em}w{c}{1.5em}w{c}{1.5em}}
\toprule
 & & & \multicolumn{10}{c|}{Disk} & \multicolumn{10}{c|}{Flower} & \multicolumn{3}{c|}{Sphere} & \multicolumn{3}{c}{Popcorn} \\
\cmidrule(lr){4-13} \cmidrule(lr){14-23} \cmidrule(lr){24-26} \cmidrule(lr){27-29}
 & & $L$ \textbackslash{} $n$ & 48 & 96 & 192 & 384 & 768 & 48 & 96 & 192 & 384 & 768 & 48 & 96 & 192 & 384 & 768 & 48 & 96 & 192 & 384 & 768 & 24 & 48 & 96 & 24 & 48 & 96 \\
\midrule\midrule
 & & & \multicolumn{5}{c|}{$k = 2$} & \multicolumn{5}{c|}{$k = 3$} & \multicolumn{5}{c|}{$k = 2$} & \multicolumn{5}{c|}{$k = 3$} & \multicolumn{3}{c|}{} & \multicolumn{3}{c}{} \\
\cmidrule(lr){1-29}
\multirow{10}{*}{$\Lgrad$} & \multirow{5}{*}{$B_{\rm hom}^{\rm g}$} & 1 & 10 & 10 & 11 & 12 & 12 & 12 & 14 & 14 & 20 & 20 & 13 & 11 & 13 & 12 & 12 & 18 & 16 & 20 & 21 & 21 &  &  &  &  &  &  \\
 &  & 2 &  & 10 & 11 & 12 & 12 &  & 14 & 14 & 19 & 19 &  & 11 & 13 & 12 & 12 &  & 15 & 20 & 20 & 21 &  &  &  &  &  &  \\
 &  & 3 &  &  & 11 & 12 & 12 &  &  & 14 & 19 & 19 &  &  & 13 & 12 & 12 &  &  & 20 & 20 & 21 &  &  &  &  &  &  \\
 &  & 4 &  &  &  & 12 & 12 &  &  &  & 19 & 19 &  &  &  & 12 & 12 &  &  &  & 20 & 21 &  &  &  &  &  &  \\
 &  & 5 &  &  &  &  & 12 &  &  &  &  & 19 &  &  &  &  & 12 &  &  &  &  & 21 &  &  &  &  &  &  \\
\cmidrule(lr){2-29}
 & \multirow{5}{*}{$B_{\rm aux}^{\rm g}$} & 1 & 11 & 12 & 12 & 12 & 12 & 11 & 11 & 12 & 12 & 12 & 12 & 12 & 12 & 12 & 12 & 13 & 12 & 12 & 12 & 12 &  &  &  &  &  &  \\
 &  & 2 &  & 12 & 12 & 12 & 13 &  & 12 & 12 & 12 & 12 &  & 13 & 13 & 12 & 12 &  & 12 & 12 & 12 & 12 &  &  &  &  &  &  \\
 &  & 3 &  &  & 13 & 13 & 13 &  &  & 13 & 12 & 12 &  &  & 13 & 13 & 13 &  &  & 13 & 13 & 12 &  &  &  &  &  &  \\
 &  & 4 &  &  &  & 13 & 13 &  &  &  & 13 & 13 &  &  &  & 14 & 14 &  &  &  & 13 & 13 &  &  &  &  &  &  \\
 &  & 5 &  &  &  &  & 13 &  &  &  &  & 13 &  &  &  &  & 14 &  &  &  &  & 13 &  &  &  &  &  &  \\
\midrule\midrule
 & & & \multicolumn{5}{c|}{$k = 1$} & \multicolumn{5}{c|}{$k = 2$} & \multicolumn{5}{c|}{$k = 1$} & \multicolumn{5}{c|}{$k = 2$} & \multicolumn{3}{c|}{$k = 1$} & \multicolumn{3}{c}{$k = 1$} \\
\cmidrule(lr){1-29}
\multirow{10}{*}{$\Ldiv$} & \multirow{5}{*}{$B_{\rm hom}^{\rm d}$} & 1 & 8 & 8 & 8 & 8 & 7 & 9 & 10 & 10 & 13 & 12 & 10 & 8 & 9 & 8 & 8 & 13 & 12 & 13 & 13 & 13 & 11 & 14 & 16 & 14 & 17 & 17 \\
 &  & 2 &  & 8 & 8 & 8 & 7 &  & 9 & 10 & 12 & 11 &  & 8 & 9 & 8 & 7 &  & 11 & 13 & 12 & 12 &  & 14 & 17 &  & 17 & 17 \\
 &  & 3 &  &  & 8 & 8 & 7 &  &  & 10 & 12 & 11 &  &  & 9 & 8 & 8 &  &  & 13 & 12 & 12 &  &  & 17 &  &  & 17 \\
 &  & 4 &  &  &  & 8 & 8 &  &  &  & 12 & 11 &  &  &  & 8 & 8 &  &  &  & 12 & 12 &  &  &  &  &  &  \\
 &  & 5 &  &  &  &  & 8 &  &  &  &  & 11 &  &  &  &  & 8 &  &  &  &  & 12 &  &  &  &  &  &  \\
\cmidrule(lr){2-29}
 & \multirow{5}{*}{$B_{\rm aux}^{\rm d}$} & 1 & 9 & 9 & 8 & 8 & 8 & 8 & 8 & 8 & 8 & 8 & 10 & 9 & 9 & 8 & 8 & 10 & 9 & 8 & 8 & 7 & 14 & 16 & 17 & 16 & 18 & 17 \\
 &  & 2 &  & 9 & 9 & 9 & 9 &  & 9 & 8 & 8 & 8 &  & 10 & 10 & 9 & 9 &  & 9 & 9 & 8 & 8 &  & 16 & 17 &  & 19 & 18 \\
 &  & 3 &  &  & 9 & 9 & 9 &  &  & 9 & 9 & 8 &  &  & 10 & 10 & 9 &  &  & 10 & 9 & 8 &  &  & 17 &  &  & 18 \\
 &  & 4 &  &  &  & 9 & 9 &  &  &  & 9 & 9 &  &  &  & 10 & 10 &  &  &  & 10 & 9 &  &  &  &  &  &  \\
 &  & 5 &  &  &  &  & 9 &  &  &  &  & 9 &  &  &  &  & 10 &  &  &  &  & 9 &  &  &  &  &  &  \\
\bottomrule
\end{tabular}

  \caption{Higher-order discretisations (\sect{subsec:results-highorder}). FGMRES iterations of $\Bgmg$ and $\Bmg$ for the Riesz maps $\Lxh$ of order $k$ on the fixed geometries, with $L+1$ levels.}
  \label{tab:gmg-highorder}
\end{table}


\subsection{Face-based stabilisation}\label{subsec:results-ghost-penalty}

The preconditioners of this work are not limited to the projection-based stabilisation in \eqref{eq:stabilisation}. The rationale behind the design of the smoothers and the multigrid hierarchy remains valid for stabilisations whose support is local to the aggregates.
Specifically, we consider a ghost-penalty stabilisation \cite{BurmanHansbo2012} defined in terms of facet jumps, which has recently been applied to the whole de Rham complex in \cite{DiPietroDroniouNilsson2025}. Following \cite{BadiaNeivaVerdugo2022}, the jumps are restricted to the facets interior to the aggregates, that is, to the set $\mathcal{F}_h^{\rm ag}$ of facets shared by two cells of the same aggregate in $\Tagg$.

Let $\mathrm{D}^{\rm x} \in \{\Grad, \Curl, \Div\}$ denote the differential operator of the space $V_h^{\rm x}$. Given $u \in V_h^{\rm x}$, let $p$ denote its polynomial degree in the direction normal to a facet, i.e.\ $p = k$ for $\Qk{k}$ and $p = k+1$ for $\RT{k}$ and $\mathbb{ND}_k$. The face-based counterpart of \eqref{eq:stabilisation-biforms} is then defined as
\begin{equation}\label{eq:facet-stabilisation-biforms}
  S^{\rm x}_{\mathcal{F}}(u,v) \coloneqq \sum_{F \in \mathcal{F}_h^{\rm ag}} \Big[ \sum_{j=0}^{p} \beta_{j} \inner{\jump{\partial_n^{j} u}}{\jump{\partial_n^{j} v}}{F} + \kappa^2 \sum_{j=1}^{p} h_F^{-2}\, \beta_{j} \inner{\jump{\partial_n^{j-1} \mathrm{D}^{\rm x} u}}{\jump{\partial_n^{j-1} \mathrm{D}^{\rm x} v}}{F} \Big],
\end{equation}
where $\jump{\cdot}$ denotes the jump across $F$, $\partial_n^{j}$ the $j$-th normal derivative, $h_F = h$ the facet size, and the weights $\beta_j$ are adapted from \cite{StickoKreiss2019}:
\begin{equation}
  \beta_{j} \coloneqq \frac{1}{j!\,\sqrt{2j+1}} \Big( \frac{h_F}{p} \Big)^{2j+1}.
\end{equation}

The corresponding bilinear forms $\Lxh$ are built as in \eqref{eq:global-biforms}. Both stabilisations control the same modes, namely the deviation of the discrete solution from a single polynomial on each aggregate, so that the smoothers of \sect{subsec:smoothers} require no modification. We repeat the experiments of \tab{tab:gmg-fixed} and \tab{tab:gmg-highorder} with this stabilisation and $\tau^2 = 1$, and report the results in \tab{tab:gmg-facestab} and \tab{tab:gmg-facestab-highorder}.
As for the projection-based stabilisation, we observe that the face-based stabilisation also leads to robust multigrid performance.

\begin{table}[htbp]
  \centering
  \scriptsize
  \setlength{\tabcolsep}{1.5pt}
  \begin{tabular}{llc|w{c}{1.5em}w{c}{1.5em}w{c}{1.5em}w{c}{1.5em}w{c}{1.5em}|w{c}{1.5em}w{c}{1.5em}w{c}{1.5em}w{c}{1.5em}w{c}{1.5em}|w{c}{1.5em}w{c}{1.5em}w{c}{1.5em}|w{c}{1.5em}w{c}{1.5em}w{c}{1.5em}}
\toprule
 & & & \multicolumn{5}{c|}{Disk} & \multicolumn{5}{c|}{Flower} & \multicolumn{3}{c|}{Sphere} & \multicolumn{3}{c}{Popcorn} \\
\cmidrule(lr){4-8} \cmidrule(lr){9-13} \cmidrule(lr){14-16} \cmidrule(lr){17-19}
 & & $L$ \textbackslash{} $n$ & 48 & 96 & 192 & 384 & 768 & 48 & 96 & 192 & 384 & 768 & 24 & 48 & 96 & 24 & 48 & 96 \\
\midrule\midrule
\multirow{10}{*}{$\Lgrad$} & \multirow{5}{*}{$B_{\rm hom}^{\rm g}$} & 1 & 10 & 10 & 11 & 11 & 11 & 10 & 11 & 11 & 11 & 11 & 17 & 18 & 20 & 18 & 20 & 20 \\
 &  & 2 &  & 10 & 11 & 11 & 11 &  & 11 & 11 & 11 & 11 &  & 18 & 20 &  & 20 & 20 \\
 &  & 3 &  &  & 11 & 11 & 11 &  &  & 12 & 12 & 11 &  &  & 20 &  &  & 20 \\
 &  & 4 &  &  &  & 11 & 11 &  &  &  & 12 & 12 &  &  &  &  &  &  \\
 &  & 5 &  &  &  &  & 12 &  &  &  &  & 12 &  &  &  &  &  &  \\
\cmidrule(lr){2-19}
 & \multirow{5}{*}{$B_{\rm aux}^{\rm g}$} & 1 & 9 & 10 & 10 & 10 & 10 & 10 & 10 & 10 & 10 & 10 & 18 & 18 & 19 & 18 & 19 & 19 \\
 &  & 2 &  & 10 & 10 & 11 & 11 &  & 11 & 11 & 10 & 11 &  & 18 & 19 &  & 19 & 19 \\
 &  & 3 &  &  & 11 & 11 & 11 &  &  & 11 & 11 & 11 &  &  & 19 &  &  & 19 \\
 &  & 4 &  &  &  & 11 & 11 &  &  &  & 12 & 12 &  &  &  &  &  &  \\
 &  & 5 &  &  &  &  & 12 &  &  &  &  & 12 &  &  &  &  &  &  \\
\midrule\midrule
\multirow{10}{*}{$\Ldiv$} & \multirow{5}{*}{$B_{\rm hom}^{\rm d}$} & 1 & 10 & 11 & 11 & 10 & 10 & 10 & 12 & 11 & 11 & 10 & 10 & 10 & 10 & 11 & 10 & 10 \\
 &  & 2 &  & 11 & 11 & 11 & 11 &  & 12 & 11 & 11 & 11 &  & 10 & 11 &  & 11 & 10 \\
 &  & 3 &  &  & 12 & 11 & 11 &  &  & 12 & 11 & 11 &  &  & 11 &  &  & 11 \\
 &  & 4 &  &  &  & 11 & 11 &  &  &  & 12 & 11 &  &  &  &  &  &  \\
 &  & 5 &  &  &  &  & 11 &  &  &  &  & 11 &  &  &  &  &  &  \\
\cmidrule(lr){2-19}
 & \multirow{5}{*}{$B_{\rm aux}^{\rm d}$} & 1 & 9 & 9 & 9 & 9 & 8 & 10 & 9 & 9 & 8 & 8 & 9 & 10 & 10 & 10 & 10 & 10 \\
 &  & 2 &  & 9 & 9 & 9 & 9 &  & 10 & 10 & 9 & 9 &  & 10 & 10 &  & 10 & 10 \\
 &  & 3 &  &  & 9 & 9 & 9 &  &  & 10 & 10 & 9 &  &  & 10 &  &  & 10 \\
 &  & 4 &  &  &  & 9 & 9 &  &  &  & 10 & 10 &  &  &  &  &  &  \\
 &  & 5 &  &  &  &  & 9 &  &  &  &  & 10 &  &  &  &  &  &  \\
\midrule\midrule
\multirow{6}{*}{$\Lcurl$} & \multirow{3}{*}{$B_{\rm hom}^{\rm c}$} & 1 &  &  &  &  &  &  &  &  &  &  & 24 & 27 & 27 & 25 & 29 & 27 \\
 &  & 2 &  &  &  &  &  &  &  &  &  &  &  & 26 & 26 &  & 29 & 26 \\
 &  & 3 &  &  &  &  &  &  &  &  &  &  &  &  & 27 &  &  & 27 \\
\cmidrule(lr){2-19}
 & \multirow{3}{*}{$B_{\rm aux}^{\rm c}$} & 1 &  &  &  &  &  &  &  &  &  &  & 25 & 25 & 25 & 25 & 25 & 25 \\
 &  & 2 &  &  &  &  &  &  &  &  &  &  &  & 25 & 25 &  & 26 & 25 \\
 &  & 3 &  &  &  &  &  &  &  &  &  &  &  &  & 25 &  &  & 25 \\
\bottomrule
\end{tabular}

  \caption{Face-based stabilisation (\sect{subsec:results-ghost-penalty}). FGMRES iterations of $\Bgmg$ and $\Bmg$ for the lowest-order Riesz maps $\Lxh$ with the stabilisation \eqref{eq:facet-stabilisation-biforms} on the fixed geometries, with $L+1$ levels.}
  \label{tab:gmg-facestab}
\end{table}

\begin{table}[htbp]
  \centering
  \scriptsize
  \setlength{\tabcolsep}{1.5pt}
  \begin{tabular}{llc|w{c}{1.5em}w{c}{1.5em}w{c}{1.5em}w{c}{1.5em}w{c}{1.5em}|w{c}{1.5em}w{c}{1.5em}w{c}{1.5em}w{c}{1.5em}w{c}{1.5em}|w{c}{1.5em}w{c}{1.5em}w{c}{1.5em}w{c}{1.5em}w{c}{1.5em}|w{c}{1.5em}w{c}{1.5em}w{c}{1.5em}w{c}{1.5em}w{c}{1.5em}}
\toprule
 & & & \multicolumn{10}{c|}{Disk} & \multicolumn{10}{c}{Flower} \\
\cmidrule(lr){4-13} \cmidrule(lr){14-23}
 & & $L$ \textbackslash{} $n$ & 48 & 96 & 192 & 384 & 768 & 48 & 96 & 192 & 384 & 768 & 48 & 96 & 192 & 384 & 768 & 48 & 96 & 192 & 384 & 768 \\
\midrule\midrule
 & & & \multicolumn{5}{c|}{$k = 2$} & \multicolumn{5}{c|}{$k = 3$} & \multicolumn{5}{c|}{$k = 2$} & \multicolumn{5}{c}{$k = 3$} \\
\cmidrule(lr){1-23}
\multirow{10}{*}{$\Lgrad$} & \multirow{5}{*}{$B_{\rm hom}^{\rm g}$} & 1 & 8 & 8 & 8 & 9 & 9 & 9 & 10 & 10 & 12 & 11 & 8 & 8 & 9 & 9 & 9 & 11 & 10 & 11 & 12 & 12 \\
 &  & 2 &  & 8 & 8 & 9 & 9 &  & 10 & 10 & 12 & 11 &  & 8 & 9 & 9 & 9 &  & 10 & 11 & 12 & 12 \\
 &  & 3 &  &  & 8 & 9 & 9 &  &  & 10 & 12 & 11 &  &  & 9 & 9 & 9 &  &  & 11 & 12 & 12 \\
 &  & 4 &  &  &  & 9 & 9 &  &  &  & 12 & 11 &  &  &  & 9 & 9 &  &  &  & 12 & 12 \\
 &  & 5 &  &  &  &  & 9 &  &  &  &  & 11 &  &  &  &  & 9 &  &  &  &  & 12 \\
\cmidrule(lr){2-23}
 & \multirow{5}{*}{$B_{\rm aux}^{\rm g}$} & 1 & 9 & 9 & 9 & 9 & 9 & 9 & 9 & 9 & 9 & 9 & 9 & 9 & 9 & 9 & 9 & 9 & 9 & 8 & 9 & 8 \\
 &  & 2 &  & 9 & 9 & 10 & 10 &  & 9 & 9 & 9 & 10 &  & 10 & 10 & 10 & 10 &  & 10 & 10 & 9 & 9 \\
 &  & 3 &  &  & 10 & 10 & 10 &  &  & 10 & 10 & 10 &  &  & 11 & 11 & 10 &  &  & 10 & 10 & 10 \\
 &  & 4 &  &  &  & 10 & 11 &  &  &  & 10 & 10 &  &  &  & 11 & 11 &  &  &  & 11 & 11 \\
 &  & 5 &  &  &  &  & 11 &  &  &  &  & 10 &  &  &  &  & 11 &  &  &  &  & 11 \\
\midrule\midrule
 & & & \multicolumn{5}{c|}{$k = 1$} & \multicolumn{5}{c|}{$k = 2$} & \multicolumn{5}{c|}{$k = 1$} & \multicolumn{5}{c}{$k = 2$} \\
\cmidrule(lr){1-23}
\multirow{10}{*}{$\Ldiv$} & \multirow{5}{*}{$B_{\rm hom}^{\rm d}$} & 1 & 7 & 7 & 7 & 7 & 7 & 9 & 10 & 9 & 11 & 11 & 8 & 7 & 7 & 7 & 7 & 11 & 10 & 11 & 12 & 11 \\
 &  & 2 &  & 7 & 7 & 7 & 7 &  & 10 & 10 & 11 & 11 &  & 7 & 7 & 7 & 7 &  & 10 & 12 & 12 & 11 \\
 &  & 3 &  &  & 7 & 7 & 7 &  &  & 10 & 11 & 11 &  &  & 7 & 7 & 7 &  &  & 11 & 12 & 11 \\
 &  & 4 &  &  &  & 7 & 7 &  &  &  & 11 & 11 &  &  &  & 7 & 7 &  &  &  & 12 & 11 \\
 &  & 5 &  &  &  &  & 7 &  &  &  &  & 11 &  &  &  &  & 7 &  &  &  &  & 11 \\
\cmidrule(lr){2-23}
 & \multirow{5}{*}{$B_{\rm aux}^{\rm d}$} & 1 & 8 & 8 & 8 & 8 & 7 & 8 & 8 & 8 & 8 & 7 & 9 & 9 & 8 & 8 & 7 & 9 & 9 & 8 & 8 & 7 \\
 &  & 2 &  & 9 & 8 & 8 & 8 &  & 9 & 8 & 8 & 8 &  & 9 & 9 & 8 & 8 &  & 10 & 9 & 8 & 8 \\
 &  & 3 &  &  & 9 & 9 & 8 &  &  & 9 & 9 & 8 &  &  & 10 & 9 & 8 &  &  & 10 & 9 & 8 \\
 &  & 4 &  &  &  & 9 & 9 &  &  &  & 9 & 9 &  &  &  & 10 & 9 &  &  &  & 10 & 9 \\
 &  & 5 &  &  &  &  & 9 &  &  &  &  & 9 &  &  &  &  & 9 &  &  &  &  & 10 \\
\bottomrule
\end{tabular}

  \caption{Face-based stabilisation (\sect{subsec:results-ghost-penalty}). FGMRES iterations of $\Bgmg$ and $\Bmg$ for the Riesz maps $\Lxh$ of order $k$ with the stabilisation \eqref{eq:facet-stabilisation-biforms} on the two-dimensional fixed geometries, with $L+1$ levels.}
  \label{tab:gmg-facestab-highorder}
\end{table}


\subsection{Dirichlet boundary conditions}\label{subsec:results-bcs}

So far, we have only considered problems with natural boundary conditions. We now briefly discuss the performance of the multigrid preconditioners when Dirichlet boundary conditions are imposed weakly through \eqref{eq:nitsche-weakforms}.

From an implementation point of view, we simply replace $\Lxh$ and $\Lxact$ with their Nitsche counterparts $\widehat{\Lxh}$ and $\widehat{\Lxact}$ in the construction of the preconditioners of \sect{sec:embedded-gmg} and \sect{sec:auxspace}. We consider the Poisson problem in primal form, which corresponds to the limit $\kappa \to \infty$ of $\widehat{\Lgradh}$ and is well posed in the presence of Dirichlet boundary conditions, as discussed in \sect{sec:formulations}, with the same exact solution as $\Lgradh$. We also consider the vector-valued Riesz maps $\widehat{\Ldivh}$ and $\widehat{\Lcurlh}$.

The Nitsche penalty on the coarse levels requires a choice. Rediscretising \eqref{eq:nitsche-weakforms} on every level with the same constant $\gamma$ gives a boundary term $\gamma/h_\ell$, the natural Nitsche scaling of each level. As an alternative, one can take $\gamma_\ell = \gamma\, h_\ell / h_L$, that is, carry the penalty $\gamma/h_L$ of the finest level to every level, at the price of over-penalising the coarse problems. With this choice, the rediscretised operator of level $\ell-1$ coincides with the Galerkin coarse operator $\Il^T \Al \Il$ except for the stabilisation term, whose aggregates are built on each level.
This question is not unique to our framework, but arises in multilevel settings for most penalty-based discretisations. For instance, we find similar discussions for interior-penalty \cite{BrennerZhao2005,AntoniettiSartiVerani2015,PanPersson2022} and immersed \cite{KothariKrause2021,SaberiMeschkeVogel2023} discretisations, with conclusions that depend on the setting. For the experiments below, we adopt the second option, that is, keeping the Nitsche penalty constant across levels. In our testing, this choice has proven to be more robust.

Results are reported in \tab{tab:gmg-dirichlet} for $\Bgmg$ and $\Bmg$, with Dirichlet boundary conditions on the whole boundary, in the setting and layout of \tab{tab:gmg-fixed}. For the Poisson problem, both preconditioners remain robust in the presence of Dirichlet boundary conditions. For the lowest-order vector-valued spaces $\RT{0}$ and $\mathbb{ND}_0$, however, the iteration counts grow with both the number of levels and the mesh size, for both preconditioners alike. This corresponds to the degraded convergence order observed for these spaces in \sect{subsec:results-convergence}: the lowest-order trace is too poor to represent the Dirichlet datum on a cut boundary. Also in line with the approximability results in \sect{subsec:results-convergence}, we find that enlarging the space is enough to recover robustness. For a higher polynomial order, tested here at $k = 1$, the iteration counts regain robustness for $\Bmg$ and $\Bgmg$.

\begin{table}[htbp]
  \centering
  \scriptsize
  \setlength{\tabcolsep}{1.5pt}
  \begin{tabular}{llc|w{c}{1.5em}w{c}{1.5em}w{c}{1.5em}w{c}{1.5em}w{c}{1.5em}|w{c}{1.5em}w{c}{1.5em}w{c}{1.5em}w{c}{1.5em}w{c}{1.5em}|w{c}{1.5em}w{c}{1.5em}w{c}{1.5em}|w{c}{1.5em}w{c}{1.5em}w{c}{1.5em}}
\toprule
 & & & \multicolumn{5}{c|}{Disk} & \multicolumn{5}{c|}{Flower} & \multicolumn{3}{c|}{Sphere} & \multicolumn{3}{c}{Popcorn} \\
\cmidrule(lr){4-8} \cmidrule(lr){9-13} \cmidrule(lr){14-16} \cmidrule(lr){17-19}
 & & $L$ \textbackslash{} $n$ & 48 & 96 & 192 & 384 & 768 & 48 & 96 & 192 & 384 & 768 & 24 & 48 & 96 & 24 & 48 & 96 \\
\midrule\midrule
\multirow{10}{*}{Poisson} & \multirow{5}{*}{$B_{\rm hom}^{\rm g}$} & 1 & 15 & 15 & 15 & 15 & 15 & 16 & 15 & 15 & 15 & 15 & 22 & 25 & 28 & 23 & 27 & 27 \\
 &  & 2 &  & 16 & 16 & 16 & 17 &  & 16 & 16 & 17 & 17 &  & 25 & 27 &  & 26 & 27 \\
 &  & 3 &  &  & 17 & 17 & 17 &  &  & 17 & 17 & 18 &  &  & 28 &  &  & 27 \\
 &  & 4 &  &  &  & 18 & 18 &  &  &  & 18 & 18 &  &  &  &  &  &  \\
 &  & 5 &  &  &  &  & 18 &  &  &  &  & 19 &  &  &  &  &  &  \\
\cmidrule(lr){2-19}
 & \multirow{5}{*}{$B_{\rm aux}^{\rm g}$} & 1 & 12 & 12 & 12 & 13 & 13 & 13 & 13 & 13 & 13 & 13 & 24 & 25 & 26 & 25 & 26 & 27 \\
 &  & 2 &  & 12 & 13 & 13 & 13 &  & 13 & 13 & 13 & 13 &  & 25 & 26 &  & 26 & 27 \\
 &  & 3 &  &  & 13 & 13 & 13 &  &  & 13 & 14 & 14 &  &  & 26 &  &  & 27 \\
 &  & 4 &  &  &  & 13 & 14 &  &  &  & 14 & 14 &  &  &  &  &  &  \\
 &  & 5 &  &  &  &  & 14 &  &  &  &  & 14 &  &  &  &  &  &  \\
\midrule\midrule
\multirow{10}{*}{$\Ldiv$, $k = 0$} & \multirow{5}{*}{$B_{\rm hom}^{\rm d}$} & 1 & 21 & 21 & 28 & 27 & 26 & 32 & 43 & 47 & 73 & 77 & 34 & 58 & 106 & 22 & 38 & 65 \\
 &  & 2 &  & 36 & 44 & 46 & 38 &  & 55 & 72 & 104 & 113 &  & 61 & 107 &  & 40 & 68 \\
 &  & 3 &  &  & 65 & 73 & 74 &  &  & 95 & 135 & 161 &  &  & 110 &  &  & 70 \\
 &  & 4 &  &  &  & 96 & 114 &  &  &  & 164 & 208 &  &  &  &  &  &  \\
 &  & 5 &  &  &  &  & 146 &  &  &  &  & 254 &  &  &  &  &  &  \\
\cmidrule(lr){2-19}
 & \multirow{5}{*}{$B_{\rm aux}^{\rm d}$} & 1 & 22 & 23 & 30 & 28 & 27 & 34 & 47 & 52 & 74 & 87 & 29 & 50 & 90 & 23 & 36 & 59 \\
 &  & 2 &  & 40 & 48 & 48 & 41 &  & 68 & 81 & 103 & 128 &  & 55 & 99 &  & 41 & 67 \\
 &  & 3 &  &  & 74 & 78 & 80 &  &  & 109 & 147 & 171 &  &  & 106 &  &  & 70 \\
 &  & 4 &  &  &  & 112 & 124 &  &  &  & 184 & 252 &  &  &  &  &  &  \\
 &  & 5 &  &  &  &  & 176 &  &  &  &  & 298 &  &  &  &  &  &  \\
\midrule\midrule
\multirow{10}{*}{$\Ldiv$, $k = 1$} & \multirow{5}{*}{$B_{\rm hom}^{\rm d}$} & 1 & 10 & 10 & 10 & 10 & 9 & 13 & 12 & 13 & 13 & 15 & 11 & 11 & 11 & 12 & 12 & 13 \\
 &  & 2 &  & 10 & 10 & 10 & 9 &  & 13 & 14 & 15 & 16 &  & 12 & 11 &  & 13 & 13 \\
 &  & 3 &  &  & 10 & 10 & 10 &  &  & 14 & 16 & 16 &  &  & 11 &  &  & 14 \\
 &  & 4 &  &  &  & 10 & 10 &  &  &  & 15 & 16 &  &  &  &  &  &  \\
 &  & 5 &  &  &  &  & 10 &  &  &  &  & 16 &  &  &  &  &  &  \\
\cmidrule(lr){2-19}
 & \multirow{5}{*}{$B_{\rm aux}^{\rm d}$} & 1 & 8 & 8 & 9 & 8 & 7 & 11 & 10 & 11 & 11 & 13 & 14 & 14 & 14 & 9 & 11 & 11 \\
 &  & 2 &  & 9 & 9 & 8 & 7 &  & 11 & 12 & 13 & 14 &  & 14 & 14 &  & 11 & 12 \\
 &  & 3 &  &  & 9 & 8 & 8 &  &  & 12 & 14 & 14 &  &  & 14 &  &  & 12 \\
 &  & 4 &  &  &  & 8 & 8 &  &  &  & 14 & 14 &  &  &  &  &  &  \\
 &  & 5 &  &  &  &  & 8 &  &  &  &  & 14 &  &  &  &  &  &  \\
\midrule\midrule
\multirow{6}{*}{$\Lcurl$, $k = 0$} & \multirow{3}{*}{$B_{\rm hom}^{\rm c}$} & 1 &  &  &  &  &  &  &  &  &  &  & 39 & 44 & 81 & 26 & 41 & 66 \\
 &  & 2 &  &  &  &  &  &  &  &  &  &  &  & 62 & 91 &  & 43 & 70 \\
 &  & 3 &  &  &  &  &  &  &  &  &  &  &  &  & 111 &  &  & 74 \\
\cmidrule(lr){2-19}
 & \multirow{3}{*}{$B_{\rm aux}^{\rm c}$} & 1 &  &  &  &  &  &  &  &  &  &  & 38 & 43 & 63 & 24 & 34 & 53 \\
 &  & 2 &  &  &  &  &  &  &  &  &  &  &  & 63 & 84 &  & 38 & 60 \\
 &  & 3 &  &  &  &  &  &  &  &  &  &  &  &  & 110 &  &  & 63 \\
\midrule\midrule
\multirow{6}{*}{$\Lcurl$, $k = 1$} & \multirow{3}{*}{$B_{\rm hom}^{\rm c}$} & 1 &  &  &  &  &  &  &  &  &  &  & 24 & 26 & 35 & 21 & 26 & 33 \\
 &  & 2 &  &  &  &  &  &  &  &  &  &  &  & 26 & 34 &  & 28 & 36 \\
 &  & 3 &  &  &  &  &  &  &  &  &  &  &  &  & 34 &  &  & 37 \\
\cmidrule(lr){2-19}
 & \multirow{3}{*}{$B_{\rm aux}^{\rm c}$} & 1 &  &  &  &  &  &  &  &  &  &  & 25 & 26 & 26 & 22 & 23 & 22 \\
 &  & 2 &  &  &  &  &  &  &  &  &  &  &  & 26 & 26 &  & 23 & 23 \\
 &  & 3 &  &  &  &  &  &  &  &  &  &  &  &  & 26 &  &  & 23 \\
\bottomrule
\end{tabular}

  \caption{Dirichlet boundary conditions (\sect{subsec:results-bcs}). FGMRES iterations of $\Bgmg$ and $\Bmg$ for the lowest-order Poisson problem and for the Riesz maps $\Ldivh$ and $\Lcurlh$ of order $k$, with Dirichlet boundary conditions on the whole boundary, on the fixed geometries, with $L+1$ levels.}
  \label{tab:gmg-dirichlet}
\end{table}

\subsection{Multiphysics preconditioners}\label{subsec:results-multiphysics}

Finally, we show how the preconditioners for the Riesz maps $\Lxh$ introduced above can be used as building blocks for a variety of multiphysics problems. Block preconditioning \cite{Mardal2011} is a widely used strategy to obtain robust and efficient solvers for coupled systems. We will here focus on its application to the mixed Hodge Laplacians of the de Rham complex. The forms of \sect{sec:formulations} for ${\rm x} = 0$ provide the $L^2$ block.

Let $({\rm x},{\rm y}) \in \{({\rm g},{\rm c}),\, ({\rm c},{\rm d}),\, ({\rm d},0)\}$ denote two consecutive spaces of the complex, so that $\mathrm{D}^{\rm x}$ maps $V_h^{\rm x}$ into $V_h^{\rm y}$. On a generic domain $\gdom$, the mixed Hodge Laplacian is given by the bilinear form
\begin{equation}\label{eq:hodge-laplacian}
  \mathcal{L}^{\rm xy}\big((\sigma,u),(\rho,v);\gdom\big) \coloneqq
  \inner{\sigma}{\rho}{\gdom} - \inner{\mathrm{D}^{\rm x} \rho}{u}{\gdom}
  - \inner{\mathrm{D}^{\rm x} \sigma}{v}{\gdom} - \inner{\mathrm{D}^{\rm y} u}{\mathrm{D}^{\rm y} v}{\gdom}.
\end{equation}
Its unfitted counterpart stabilises each term of \eqref{eq:hodge-laplacian} against the space it belongs to, that is
\begin{equation}
  \begin{aligned}
    \mathcal{L}^{\rm xy}_h\big((\sigma,u),(\rho,v)\big) &\coloneqq \mathcal{L}^{\rm xy}\big((\sigma,u),(\rho,v);\dom\big) + \tau^2 S^{\rm xy}\big((\sigma,u),(\rho,v)\big), \\
    S^{\rm xy}\big((\sigma,u),(\rho,v)\big) &\coloneqq
    S_{\rm x}^{\rm ag}(\sigma,\rho) - S_{\rm y}^{\rm ag}(\mathrm{D}^{\rm x} \rho, u)
    - S_{\rm y}^{\rm ag}(\mathrm{D}^{\rm x} \sigma, v) - S_{{\rm y}'}^{\rm ag}(\mathrm{D}^{\rm y} u, \mathrm{D}^{\rm y} v),
  \end{aligned}
\end{equation}
where ${\rm y}'$ is defined as in \sect{sec:formulations} (for ${\rm y} = 0$ the last term vanishes, since $\mathrm{D}^0 = 0$). The unfitted problem then reads: find $(\sigma, u) \in V_h^{\rm x} \times V_h^{\rm y}$ such that
\begin{equation}\label{eq:hodge-problem}
  \mathcal{L}^{\rm xy}_h\big((\sigma,u),(\rho,v)\big) = \inner{g}{\rho}{\dom} - \inner{f}{v}{\dom} \qquad \forall (\rho, v) \in V_h^{\rm x} \times V_h^{\rm y}.
\end{equation}
The pair $({\rm c},{\rm d})$ is only considered in three dimensions, and the pair $({\rm d},0)$ is the mixed form of the Poisson problem.

Following \cite{Mardal2011}, problem \eqref{eq:hodge-problem} is preconditioned by the inverse of the Riesz map of the product space $V_h^{\rm x} \times V_h^{\rm y}$, that is the block-diagonal operator
\begin{equation}\label{eq:hodge-precond}
  \mathcal{B}^{\rm xy} \coloneqq
  \begin{pmatrix}
    (A_h^{\rm x})^{-1} & 0 \\
    0 & (A_h^{\rm y})^{-1}
  \end{pmatrix},
\end{equation}
where $\Ah$ is the operator induced by $\Lxh$ on $V_h^{\rm x}$, with $\kappa = 1$. We then replace the exact inverses $(A_h^{\rm x})^{-1}$ and $(A_h^{\rm y})^{-1}$ in \eqref{eq:hodge-precond} by preconditioners $B^{\rm x}$ and $B^{\rm y}$ of the previous sections. For ${\rm y} = 0$, the block $A_h^0$ induced by $\Lambda^0_h$ on $V_h^0$ is treated by the same two preconditioners, built on the discontinuous space with the same hierarchy, transfer operators and patch smoothers. We denote by $\mathcal{B}^{\rm xy}_{\rm hom}$ and $\mathcal{B}^{\rm xy}_{\rm aux}$ the preconditioners obtained when each block is approximated by five FGMRES iterations preconditioned by $\Bgmg$ and $\Bmg$, respectively.

The exact solution of each mixed problem is the pair formed by the exact solutions of its two Riesz maps. Results are reported in \tab{tab:gmg-hodge} for the three problems on the fixed geometries, with the same layout as \tab{tab:gmg-fixed}, where along a row the finest mesh is refined at a fixed number of levels, and along a column levels are added at a fixed finest mesh.

\begin{table}[htbp]
  \centering
  \scriptsize
  \setlength{\tabcolsep}{1.5pt}
  \begin{tabular}{llc|w{c}{1.5em}w{c}{1.5em}w{c}{1.5em}w{c}{1.5em}w{c}{1.5em}|w{c}{1.5em}w{c}{1.5em}w{c}{1.5em}w{c}{1.5em}w{c}{1.5em}|w{c}{1.5em}w{c}{1.5em}w{c}{1.5em}|w{c}{1.5em}w{c}{1.5em}w{c}{1.5em}}
\toprule
 & & & \multicolumn{5}{c|}{Disk} & \multicolumn{5}{c|}{Flower} & \multicolumn{3}{c|}{Sphere} & \multicolumn{3}{c}{Popcorn} \\
\cmidrule(lr){4-8} \cmidrule(lr){9-13} \cmidrule(lr){14-16} \cmidrule(lr){17-19}
 & & $L$ \textbackslash{} $n$ & 48 & 96 & 192 & 384 & 768 & 48 & 96 & 192 & 384 & 768 & 24 & 48 & 96 & 24 & 48 & 96 \\
\midrule\midrule
\multirow{10}{*}{$\mathcal{L}^{\rm gc}$} & \multirow{5}{*}{$\mathcal{B}_{\rm hom}^{\rm gc}$} & 1 & 10 & 10 & 10 & 10 & 9 & 9 & 11 & 10 & 9 & 9 & 22 & 24 & 27 & 25 & 35 & 27 \\
 &  & 2 &  & 11 & 11 & 11 & 9 &  & 11 & 11 & 11 & 10 &  & 21 & 25 &  & 35 & 25 \\
 &  & 3 &  &  & 11 & 11 & 11 &  &  & 11 & 11 & 11 &  &  & 25 &  &  & 28 \\
 &  & 4 &  &  &  & 11 & 11 &  &  &  & 11 & 11 &  &  &  &  &  &  \\
 &  & 5 &  &  &  &  & 11 &  &  &  &  & 11 &  &  &  &  &  &  \\
\cmidrule(lr){2-19}
 & \multirow{5}{*}{$\mathcal{B}_{\rm aux}^{\rm gc}$} & 1 & 8 & 8 & 8 & 8 & 8 & 8 & 8 & 8 & 8 & 8 & 19 & 18 & 19 & 18 & 18 & 19 \\
 &  & 2 &  & 8 & 8 & 8 & 8 &  & 8 & 8 & 8 & 8 &  & 18 & 18 &  & 18 & 19 \\
 &  & 3 &  &  & 8 & 8 & 8 &  &  & 8 & 8 & 8 &  &  & 19 &  &  & 20 \\
 &  & 4 &  &  &  & 8 & 8 &  &  &  & 8 & 8 &  &  &  &  &  &  \\
 &  & 5 &  &  &  &  & 8 &  &  &  &  & 8 &  &  &  &  &  &  \\
\midrule\midrule
\multirow{6}{*}{$\mathcal{L}^{\rm cd}$} & \multirow{3}{*}{$\mathcal{B}_{\rm hom}^{\rm cd}$} & 1 &  &  &  &  &  &  &  &  &  &  & 15 & 15 & 15 & 17 & 20 & 15 \\
 &  & 2 &  &  &  &  &  &  &  &  &  &  &  & 14 & 13 &  & 20 & 15 \\
 &  & 3 &  &  &  &  &  &  &  &  &  &  &  &  & 15 &  &  & 17 \\
\cmidrule(lr){2-19}
 & \multirow{3}{*}{$\mathcal{B}_{\rm aux}^{\rm cd}$} & 1 &  &  &  &  &  &  &  &  &  &  & 14 & 13 & 13 & 13 & 13 & 13 \\
 &  & 2 &  &  &  &  &  &  &  &  &  &  &  & 13 & 13 &  & 13 & 13 \\
 &  & 3 &  &  &  &  &  &  &  &  &  &  &  &  & 13 &  &  & 13 \\
\midrule\midrule
\multirow{10}{*}{$\mathcal{L}^{\rm d0}$} & \multirow{5}{*}{$\mathcal{B}_{\rm hom}^{\rm d0}$} & 1 & 9 & 9 & 9 & 9 & 9 & 9 & 10 & 10 & 9 & 9 & 9 & 10 & 10 & 10 & 10 & 10 \\
 &  & 2 &  & 10 & 10 & 10 & 9 &  & 10 & 10 & 10 & 10 &  & 10 & 10 &  & 10 & 10 \\
 &  & 3 &  &  & 10 & 10 & 10 &  &  & 10 & 10 & 10 &  &  & 10 &  &  & 10 \\
 &  & 4 &  &  &  & 10 & 10 &  &  &  & 10 & 10 &  &  &  &  &  &  \\
 &  & 5 &  &  &  &  & 10 &  &  &  &  & 10 &  &  &  &  &  &  \\
\cmidrule(lr){2-19}
 & \multirow{5}{*}{$\mathcal{B}_{\rm aux}^{\rm d0}$} & 1 & 8 & 8 & 8 & 8 & 8 & 8 & 7 & 7 & 7 & 7 & 9 & 9 & 10 & 9 & 10 & 10 \\
 &  & 2 &  & 8 & 8 & 8 & 8 &  & 8 & 7 & 7 & 7 &  & 9 & 10 &  & 10 & 10 \\
 &  & 3 &  &  & 8 & 8 & 8 &  &  & 7 & 7 & 7 &  &  & 10 &  &  & 10 \\
 &  & 4 &  &  &  & 8 & 8 &  &  &  & 7 & 7 &  &  &  &  &  &  \\
 &  & 5 &  &  &  &  & 8 &  &  &  &  & 7 &  &  &  &  &  &  \\
\bottomrule
\end{tabular}

  \caption{Multiphysics preconditioners (\sect{subsec:results-multiphysics}). FGMRES iterations of $\mathcal{B}^{\rm xy}_{\rm hom}$ and $\mathcal{B}^{\rm xy}_{\rm aux}$ for the lowest-order mixed Hodge Laplacians $\mathcal{L}^{\rm xy}_h$ on the fixed geometries, with $L+1$ levels.}
  \label{tab:gmg-hodge}
\end{table}

\section{Conclusion}

We have designed two families of multigrid preconditioners for stabilised unfitted discretisations of the $H(\Grad)$-, $H(\Curl)$- and $H(\Div)$-conforming Riesz maps. The homogeneous \ac{gmg} preconditioner applies the unfitted discretisation on every level of a hierarchy that inherits the boundary representation of the finest mesh, and relies on star patches built on the aggregates of the stabilisation. The auxiliary-space \ac{gmg} preconditioner instead confines the unfitted operator to the finest level, where it is treated by the same patch smoother, and delegates the remaining work to a body-fitted problem on the active mesh, solved by a non-nested \ac{gmg} on a hierarchy of active meshes. Numerical experiments in two and three dimensions show that both preconditioners are robust with respect to the mesh size, the number of levels, the cut ratio $\eps$ and the weight $\kappa$ of the differential term, for the aggregate-wise projection stabilisation as well as for the face-based one, and that this robustness carries over to higher polynomial degrees. For Dirichlet boundary conditions imposed by Nitsche's method, the preconditioners retain the same robustness for high-order discretisations, but are not robust for the lowest-order case. This is in line with the observed degraded approximability of the lowest-order methods. The same preconditioners serve as the diagonal blocks of operator preconditioners for the mixed Hodge Laplacians of the de Rham complex, with the same robustness, which makes them building blocks for multiphysics problems beyond the Riesz maps they were designed for.

Several questions remain open. The aggregates couple cells across process boundaries, so an efficient distributed-memory implementation of the aggregate-based patch smoothers, and of the aggregation itself on coarse levels, is not immediate. Also, some multiphysics problems of interest, such as the incompressible elasticity and Navier--Stokes equations, do not fit in the de Rham complex. In the body-fitted setting, robust multigrid strategies exist that revolve around the same principles as those discussed here \cite{Schöberl1999,10.1093/imanum/drab083,doi:10.1137/18M1219370,osti_2424703}. Thus, the methods designed in this work should extend to these problems in the unfitted setting as well.
Then, recovering robustness for the lowest-order $H(\Curl)$- and $H(\Div)$-conforming discretisations under Dirichlet boundary conditions remains an open question.
Finally, the robustness observed here is supported by numerical evidence only. A rigorous analysis has to deal with hierarchies that are non-nested, which requires a substantially different approach from the body-fitted theory, and is left for future work.

\section*{Acknowledgements}

This research was partially funded by the Australian Government through the Australian Research Council (project numbers DP210103092 and DP220103160). This work was also supported by computational resources provided by the Australian Government through NCI and Pawsey under the NCMAS Merit Allocation Schemes, as well as computational resources provided by Monash University through NCI under the Monash-NCI allocation scheme.

\printbibliography

\end{document}